 \documentclass{amsart}

\usepackage{epsfig}  		
\usepackage{epic,eepic}       

\usepackage{a4wide,amssymb,latexsym,amsthm,amsmath}
\usepackage{eucal}
\usepackage{graphicx}
\usepackage{algorithm} 
\usepackage{algorithmic}
\usepackage{graphicx}
\usepackage{subcaption}
\usepackage{bbm}
\usepackage{lipsum}
\usepackage{placeins}

\newcommand{\argmin}{\mathrm{arg\,min}}

\newtheorem{thm}{Theorem}[section]

\theoremstyle{definition}
\newtheorem{definition}[thm]{Definition}

\theoremstyle{remark}

\newtheorem{remark}[thm]{Remark}

\numberwithin{equation}{section}

\begin{document}


\title[Fast algorithm of the shear-compression damage model for block caving]{A fast algorithm of the shear-compression damage model for the simulation of block caving}


\author[Gaete]{Sergio Gaete}

\author[Jofre]{Alejandro Jofre}

\author[Lecaros]{Rodrigo Lecaros}

\author[Montecinos]{Gino Montecinos}

\author[Ortega]{Jaime H. Ortega}

\author[Ram\'irez-Ganga]{Javier Ram\'irez-Ganga}

\author[San Mart\'in]{Jorge San Mart\'in}



\begin{abstract}
For underground mine, the current usual technique for ore extraction is block caving, which generates and induces seismic activity in the mine. To understand block caving method is one of the most challenging problems in underground mining. This method relies on gravity to break and transport large amounts of ore and waste. The state of art in damage models is not able to represent the real effect of the mining in the rock mass since for example the damage appears in the bottom of the domain under consideration and with this is not possible recover the subsidence sees in the mine. In this paper we present the analysis and implementation of the shear-compression damage model applied to underground mining proposed in \cite{bonnetier2020shearcompression}. We propose a fast algorithm based in the usual alternated algorithm used in gradient damage models \cite{marigo2016overview} and show that this new algorithm is faster than the usual algorithm. We show some numerical tests in 3D and present interested simulations for different damage laws producing realistic damage behavior. 
\end{abstract}


 \maketitle



\section{Introduction}

Worldwide the mining activity is facing lower grades, deeper ore bodies and stronger stresses within increasing massive operations. These challenges have several technical difficulties, uncertainties and associated systemic risks, which can affect the business survival and the long-term sustainability of the operation. The block caving method consists in that ore blocks are undermined causing the rocks to cave, and thus allowing broken ore to be removed at draw-points (see \cite{hartman2002introductory}). This method is based in the sinking principle and mineral breakage due to the removal of a large supporting area of rock and the subsequent extraction of this  by mean of tunnels or ditches collectors. The vacuum generated by the extraction of this material from the basis is filled by the material falling by the action of gravity, which added to the process of attrition given by the friction during its falling defines the size of the mineral at the extraction point.

In this paper we study the effect that occurs in the rock mass as mining activity develops by the block caving process, where we analize a model of rock mechanics considering damage proposed in \cite{bonnetier2020shearcompression}, in particular, we will consider a variational formulation where gradient damage models appear as an elliptical approach to the problem of  variational fracture mechanics  \cite{pham2010approche1,pham2010approche2,pham2011gradient}. To model this, we consider a homogeneous body $\Omega_0 \subseteq \mathbb{R}^3$, that represents the reference rock mass, and consider subdomains $S_i \subseteq \Omega_0$, $i=1,\cdots,M$, that represent the interior cavities with $S_i \subseteq S_{i+1}$, where our domain to model the block caving process will be $\Omega_i = \Omega_0 \setminus S_i$. In each time-step domain $\Omega_i$, we develop the damage model mentioned above of rock mechanics, considering the irreversibility of the damage in the domain.

The key issue of models predicting fracture is Griffith's criterion \cite{griffith1921vi}. This criterion supposes that, as a crack grows, the displacement field is instantly in a new equilibrium, since the displacement may be discontinuous across the crack increment. The resulting decrease in stored elastic energy can then be balanced with the work required to create the crack increment, postulated to be proportional to the newly created area. The proportionality constant is usually known as fracture toughness. In other words, the rate of elastic energy decreases per unit area, the energy release rate, is proportional to the fracture toughness. Griffith's criterion stipulates that the crack grows only if the energy release rate is equal to the fracture toughness. Traditionally, these ideas could be formalized only for relatively simple crack topologies and often only for a pre-defined crack path. Only recently was the theory of brittle fracture freed from this restriction \cite{ambrosio1997energies,francfort1998revisiting}. Ambrosio and Braides \cite{ambrosio1997energies} propose minimizing the sum of stored elastic energy and surface energy of discontinuity sets, to obtain displacements that are stable in the sense of Griffith. The first well-possed mathematical models of quasi-static fracture can be found in Dal Maso and Toader \cite{dal2002model}, Francfort and Larsen \cite{francfort2003existence}, and Francfort and Marigo \cite{francfort1998revisiting}. In these references, the Dirichlet data $u_D$ is varying in time and an evolution $u$ is sought such that, at each time $t$, $u(t)$ minimizes the potential energy subject to the Dirichlet boundary condition, and subject to an irreversibility constraint on the crack set. In Bourdin, Larsen and Richardson \cite{Bourdin2011}, a discrete-time candidate for such a model is proposed, based on the Ambrosio-Tortorelli approximation, which $\Gamma$-converges, to the Griffith energy, see Ambrosio and Tortorelli \cite{ambrosio1990approximation}. The Ambrosio-Tortorelli approximation is particularly convenient for numerical implementation and was proposed in Bourdin, Francfort and Marigo \cite{bourdin2000numerical} and Bourdin \cite{bourdin2007numerical} for the simulation of the quasi-static model.

The fundamental link between gradient damage models and Griffith model of fracture mechanics \cite{griffith1921vi,ambrosio1997energies,francfort1998revisiting} relies strongly on the variational structure of both models. Indeed, gradient damage models appear as an elliptic approximation of the variational fracture mechanics problem. The variational approach of brittle fracture recasts the evolution problem for the cracked state body as a minimality principle for an energy functional sum of the elastic energy and the energy dissipated to create the crack \cite{francfort1998revisiting}. Mathematical results based on $\Gamma$-Convergence theory show that when the internal length of a large class of gradient damage models tends to zero, the global minimum of the damage energy functional tends towards the global minima of the energy functional of Griffith brittle fracture \cite{braides1998approximation}. 

The aim of this work is to give numerical tools to the study of the mathematical model which describes the block caving process. For this, we consider the damage model proposed in \cite{bonnetier2020shearcompression} where, for the damage criterion, the authors decompose the stress tensor in its deviatoric and spherical part in order to be able to control the damage produce by the normal and shear components of the stress. Indeed, one would like to understand the impact of creating an extraction cavity on the state of the material above it, so as to assess security issues for the mining operation or so as to optimize extraction. We propose a fast algorithm to model this approach in order to improve the velocity of the convergence in the anternate algortithm and use this new algorithm with different damage laws.

This paper is organized as follows. Section \ref{SCDModel} concerns a brief presentations of the damage model to solve the block caving process. Different damage laws are used considering the damage model in Section \ref{DamLaws} and report 3D test cases. In section \ref{Cav}, we report the cavity test cases in order to model the block caving proces. Finally the Section \ref{NewAlg} describes, briefly, the numerical strategy to attack the underground mining problem and shows 2D test cases where our new damage model is compared with other existing models, and present different interesting simulations in the sense of the underground mining. In Section \ref{Con}, the conclusions and comments are drawn.

In all the work we will use the summation convention on repeated indices: Vectors and second order tensors are indicated by a lowercase letter, such as $u$ and $\sigma$ for the displacement field and the stress field. Their components are denoted by $u_i$ and $\sigma_{ij}$. Third or fourth order tensors as well as their components are indicated by a capital letter, such as $A$ or $A_{ijkl}$ for the stiffness tensor. Such tensors are considered to be linear maps acting on vectors or second order tensors and this action is denoted without dots, such as $A\varepsilon$ whose $ij$-component is $A_{ijkl}\varepsilon_{kl}$. The inner product between two vectors or two tensors of the same order is indicated by a colon. For instance $a:b$ stands for $a_ib_i$ and $\sigma : \varepsilon$ for $\sigma_{ij}\varepsilon_{ij}$. We use the notation $A>0$ to denote a positive definite tensor. The time derivative is indicated by a dot, for instance $\dot{\alpha}:=\frac{\partial \alpha}{\partial t}$.

\section{Shear-compression damage model}\label{SCDModel}

\subsection{Setting of the model}

We introduce the main ideas of the shear-compression damage model by a variational approach in a quasi-static setting, the readers interested by more details should refer \cite{bonnetier2020shearcompression}. We consider a homogeneous $n$-dimensional body whose reference configuration is the open connected bounded set $\Omega \subseteq \mathbb{R}^n$ and we assume that the local elastic material behavior can be characterized by two constants: the Young's modulus $E$ and the Poisson's ratio $\nu$. This body is made of a damaging material whose behavior is defined as follow.

The damage can be represented by a scalar parameter $\alpha \in [0,1]$, where $\alpha=0$ is the undamaged state and $\alpha=1$ is the fully damage state. For this damage parameter we assume the irreversibility condition to prevent the unphysical healing of the material, that is, if an area of the material is damaged, it is not recovered. The behavior of the material is characterized by the state function $W_{\ell}$ which gives the energy density at each point $x$. It depends on the local strain $\varepsilon(u)$, defined by $\varepsilon(u)=\frac{\nabla u + \nabla u^T}{2}$, the local damage $\alpha(x)$ and the local gradient $\nabla \alpha(x)$ of the damage field al $x$. Specifically, we assume that $W_{\ell}$ takes the following form
\begin{equation}\label{statefunc}
W_{\ell}(\varepsilon(u),\alpha,\nabla \alpha):= \frac{1}{2}A(\alpha)\varepsilon(u):\varepsilon(u)+w(\alpha)+\frac{1}{2}w_1\ell^2\vert \nabla \alpha \vert^2,
\end{equation}
where $\alpha \mapsto w(\alpha)$ describes the damage dissipation during a homogeneous damage evolution and its maximal value $w(1)=w_1$, with $0<w_1<\infty$, is the energy completely dissipated during such process when damage attains 1. $A(\alpha)$ represents the rigidity tensor of the material in the damage state $\alpha$, the material becomes less rigid when $\alpha$ increases, that is, 
\begin{equation}
    A(\alpha)= \mathrm{a}(\alpha) A_0,
\end{equation}
where $\alpha \mapsto \mathrm{a}(\alpha)$ is a decreasing scalar function such that $\mathrm{a}(0)=1$ and $\mathrm{a}(1)=0$, and where $A_0$ is the isotropic elasticity tensor. The last term in the right-hand side of (\ref{statefunc}) is the non-local part of the energy, where, for simplicity, is assumed to be a quadratic function of the gradient of damage. This term plays a regularizing role by limiting the possibilities of localization of the damage field, where the parameter $\ell >0$  can be considered as an internal length characteristic of the material, which controls the thickness of the damage localization zones.

We assume that the body $\Omega$ is subjected to a time dependent loading $U(\cdot,t)$, $F(\cdot,t)$ and $f(\cdot,t)$ which consists in an imposed displacement on the part of the boundary $\partial \Omega_U$, surface forces on the complementary part of the boundary $\partial \Omega_F$, and volume forces, $t$ denoting the time parameter. The space of kinematically admissible displacement fields at time $t$ is the set
\begin{equation}
\mathcal{C}_t=\left\{ v \in (H^1(\Omega))^n \mbox{ : } v=U(\cdot,t) \mbox{ on } \partial \Omega_U \right\},
\end{equation}
whereas the set of accessible damage state from the current damage set $\alpha_t$ is, by virtue of the irreversibility condition,
\begin{equation}
  \mathcal{D}_t(\alpha)=\left\{ \beta \mbox{ : } \alpha(x,t) \leq \beta(x) \leq 1 \mbox{ in } \Omega\right\}.
\end{equation}

In the framework of the simulation of the block caving process, we consider the constitutive equations given in \cite{bonnetier2020shearcompression}, where we have the following relationships
\begin{eqnarray}
\sigma(u,\alpha)&=& \mathrm{a}(\alpha)A_0\varepsilon(u),\\ 
\frac{1}{2}H(\varepsilon,\alpha)+w'(\alpha)-w_1\ell^2\Delta \alpha &=&-R(\dot{\alpha}),
\end{eqnarray}
where, $R(\dot{\alpha}) = 0$ if $\dot{\alpha}\geq 0$ and $R(\dot{\alpha})\in (-\infty,0]$ if $\dot{\alpha}= 0$ and $H(\varepsilon,\alpha)$ is a function representing the damage generator defined by
\begin{equation}
H(\varepsilon,\alpha) = \frac{\partial}{\partial \alpha}\left( \frac{1}{E}\left(\sigma^d:\sigma^d - \frac{2}{3} \kappa \sigma^s:\sigma^s\right)\right).
\end{equation}
with $\sigma^d$ and $\sigma^s$ the deviatoric and spheric part of the tensor stress $\sigma=\sigma(u,\alpha)=\mathrm{a}(\alpha)A_0\varepsilon(u)$. To simplify the presentation, we only consider $\kappa =1$, the readers interested by the case where $\kappa >0$ should refer to \cite{bonnetier2020shearcompression} where a complete analysis is made. Then, in order to model the damage in the underground mining process, the evolution model for block caving satisfies the following conditions
\begin{enumerate} 
	\item The stress tensor $\sigma(x,t)=\sigma(u(x,t),\alpha(x,t))=\mathrm{a}(\alpha(x,t))A_0\varepsilon(u(x,t))$ satisfies the equilibrium equations
	\begin{equation}
	\begin{array}{r c l c l}
	\mathrm{div}(\sigma(x,t)) +f(x,t)  &= &0 & \mbox{ in } &\Omega, \\
	\sigma(x,t) \cdot n &= &F(x,t) &\mbox{ on } &\partial \Omega_F, \\
	u(x,t) &= &U(x,t) &\mbox{ on } &\partial \Omega_U.
	\end{array}
	\end{equation}
	
	\item The damage field $\alpha(x,t)$ satisfies the nonlocal damage criterion in $\Omega$
	\begin{equation}\label{SCdamcrit}
	\frac{(\mathrm{a}^2(\alpha(x,t)))'}{2E}\left( (A_0\varepsilon)^d:(A_0\varepsilon)^d - \frac{2}{3}(A_0\varepsilon)^s:(A_0\varepsilon)^s \right)+ w'(\alpha(x,t))-w_1\ell^2 \Delta \alpha(x,t) \geq 0,
	\end{equation}
	and the nonlocal consistency condition in $\Omega$
	\begin{equation}
	\left(\frac{(\mathrm{a}^2(\alpha(x,t)))'}{2E}\left( (A_0\varepsilon)^d:(A_0\varepsilon)^d - \frac{2}{3} (A_0\varepsilon)^s:(A_0\varepsilon)^s \right)+ w'(\alpha(x,t))-w_1\ell^2 \Delta \alpha(x,t) \right) \dot{\alpha}(x,t) = 0, 
	\end{equation}
\end{enumerate}

\subsection{Numerical implementation}

The numerical solution considers a discretization version of the evolution problem. We use a alternate algorithm that consist in solving a series of subproblems on $u$ at fixed $\alpha$, and vice versa on $\alpha$ at fixed $u$, until convergence. Given the displacement and the damage field $(u_{i-1},\alpha_{i-1})$ at time step $t_{i-1}$, the solution at time step $t_i$ is obtained by solving the following alternate problem.  First, for the displacement $u_i$, we solve the following variational problem: Find $u$ such that
\begin{eqnarray}\label{elasteq} 
\forall  v \in \mathcal{C}_{t_i}, & \displaystyle \int_{\Omega}\sigma(u,\alpha_{i-1}):\varepsilon(v) dx = \int_{\Omega} f_i \cdot v dx + \int_{\partial \Omega_F} F_i \cdot v dS,
\end{eqnarray}
Later, with $u_i$ calculated in the previous variational problem, we obtain $\alpha_i$ solving the following bound-constrained minimization problem
\begin{equation}\label{moddamageeq}
\inf\left\{ \overline{\mathcal{P}}(u_i,\alpha) \mbox{ : } \alpha \in \mathcal{D}_t(\alpha_{i-1})\right\},
\end{equation}
where
\begin{equation}\label{moddamageeqmod3}
\overline{\mathcal{P}}(u,\alpha)=\int_{\Omega} \frac{\mathrm{a}^2(\alpha)}{2E} \left( \left(A_0\varepsilon(u)\right)^d:\left(A_0\varepsilon(u)\right)^d-\frac{2}{3}\kappa \left(A_0\varepsilon(u)\right)^s:\left(A_0\varepsilon(u)\right)^s \right)+w(\alpha)+w_1\ell^2 \vert \nabla \alpha\vert^2.
\end{equation}

The unilateral constraint $\alpha(x)\geq \alpha_{i-1}$ is the time-discrete version of the irreversibility of damage. The solution strategy is summarized in Algorithm \ref{alg:damage}. The unilateral constraint included in the damage problem requires the use of variational inequalities solvers, here we use the open-source library PETSc \cite{balay2004petsc}, where this capability is available. The problem is discretized in space with standard triangular finite elements with piecewise linear approximation for $u$ and $\alpha$. Finite element implementations based on the open-source FEniCS library \cite{langtangen2016solving, logg2012automated} are used for this problems. 
\begin{algorithm}[h!]
	\begin{algorithmic}[1]
		\RETURN Solution of time step $t_i$.
		\STATE Given $(u_{i-1},\alpha_{i-1})$, the sate at the previous loading step.
		\STATE Set $(u^{(0)},\alpha^{(0)}):=(u_{i-1},\alpha_{i-1})$ and error$^{(0)} = 1.0$ 
		\WHILE {error$^{(p)}>$tolerance}
		\STATE Solve $u^{(p)}$ from (\ref{elasteq}) with $\alpha^{(p-1)}$.
		\STATE Find $\displaystyle\alpha^{(p)}:=  \argmin_{\alpha \in \mathcal{D}(\alpha_{i-1})} \overline{\mathcal{P}}(u^{(p)},\alpha)$.
		\STATE error$^{(p)} = \Vert\alpha^{(p-1)} -\alpha^{(p)}\Vert_{\infty}$.
		\ENDWHILE
		\STATE  Set $(\mathbf{u}_{i}, \alpha_{i}) = (u^{p}, \alpha^{p}).$
	\end{algorithmic}
	\caption{Numerical algorithm to solve the damage problem}\label{alg:damage}
\end{algorithm}

\section{Study of the different damage laws}\label{DamLaws}

In order to have a better approximation of impact of the underground mining in the rock mass, it is necessary to model different materials that can represent the different types of rocks found in mining deposits. In this section we will study numerical results to analyze the different combinations of $\mathrm{a}(\alpha)$ and $w(\alpha)$.

\subsection{Hardening properties}

For a homogeneous damage distribution $\alpha$, let us define the elastic domains in the strain and stress spaces by 
\begin{eqnarray}\label{elasticspaces}
\mathcal{R}(\alpha) =\left\{ \varepsilon \in \mathbb{M}_{s} \mbox{ : } \frac{1}{2}\mathrm{a}'(\alpha)A_0\varepsilon : \varepsilon  \leq w'(\alpha)\right\}, & \displaystyle \mathcal{R}^*(\alpha)=\left\{ \sigma \in \mathbb{M}_{s} \mbox{ : } \frac{1}{2} S'(\alpha)\sigma : \sigma  \leq w'(\alpha) \right\},
\end{eqnarray}
where $\mathbb{M}_s$ is the space of symmetric tensors and $S(\alpha)=A^{-1}(\alpha)$.

The evolution of the sizes of (\ref{elasticspaces}) play a important role in the qualitative properties in the evolution damage problems. To see this, first we introduce the following definition
\begin{definition}[Hardening properties]
	We say that the material behavior is: {\bf Strain Hardening} if $\alpha\mapsto \mathcal{R}(\alpha)$ is increasing, {\bf Stress Hardening} if $\alpha \mapsto \mathcal{R}^*(\alpha)$ is increasing and {\bf  Stress Softening} if $\alpha \mapsto \mathcal{R}^*(\alpha)$ is decreasing.
\end{definition}

\begin{remark}
	The monoticity properties mut be understood in the sense of the set inclusion. For example, the Strain Hardening propertie means that
	\begin{equation}
	\alpha_1 < \alpha_2 \Rightarrow \mathcal{R}(\alpha_1)\subseteq \mathcal{R}(\alpha_2).
	\end{equation}
\end{remark}

Depending on the behavior of the quadratic forms defining these domains, the material is said to be
\begin{itemize}
	\item {\it Strain Hardening} when $\alpha \mapsto \left(-A'(\alpha)/w'(\alpha)\right)$ is decreasing with respect to $\alpha$, i.e.,
	\begin{equation}
	w'(\alpha) A''(\alpha) - w''(\alpha) A'(\alpha) >0.
	\end{equation}
	
	\item {\it Stress Hardening} (resp. Softening) when $\alpha \mapsto \left(S'(\alpha)/w'(\alpha)\right)$ is decreasing (resp. increasing) with respect to $\alpha$, i.e.,
	\begin{equation}
	w'(\alpha) S''(\alpha) - w''(\alpha) S'(\alpha) <0 \mbox{ (resp. $>$ 0 )}.
	\end{equation}
\end{itemize}

\subsection{Models for damage laws}\label{Models}

For the study of the different behavior, following \cite{marigo2016overview}, we consider four different models, with its damage laws and hardening properties respectively. This models are defined as follows:
\begin{itemize}
	\item {\bf Model 1:} We consider a model with an elastic phase. For this, we assume a linear function $w(\alpha)$ and quadratic function $\mathrm{a}(\alpha)$, that is, we take the following damage law
	\begin{eqnarray}
	w(\alpha) = w_1\alpha, & \mathrm{a}(\alpha)  = (1-\alpha)^2.
	\end{eqnarray}	
	This law satisfies both the strain hardening and stresss softening conditions for $\alpha \in [0,1)$.
	
	\item {\bf Model 2:} We consider the original Ambrosio-Tortorelli regularization model. For this, we assume a quadratic function $w(\alpha)$ and quadratic function $\mathrm{a}(\alpha)$, that is, we take the following damage law
	\begin{eqnarray}
	w(\alpha) = w_1\alpha^2, & \mathrm{a}(\alpha)  = (1-\alpha)^2.
	\end{eqnarray}
	This law satisfies the condition strain hardening for $\alpha \in [0,1)$ and stress softening only for $\alpha \geq 1/4$.
	
	\item {\bf Model 3:} We consider a family of models with the same homogeneous strain-stress response. For this, we consider the following family of damage models indexed by parameter $p>0$
	\begin{eqnarray}
	w(\alpha) = w_1\left(1-(1-\alpha)^{p/2}\right), & \mathrm{a}(\alpha)  = (1-\alpha)^p.
	\end{eqnarray}
	This case is a generalization of the Model 1 which is recovered for $p=2$. It satisfies both the strain hardening and stress softening conditions for $\alpha \in [0,1)$ and any $p>0$.
	
	\item {\bf Model 4:} We consider a model where ultimate fracture occurs at finite strain. For this, we assume that it is defined by the following material functions parametrized by a scalar parameter $k>1$
	\begin{eqnarray}
	w(\alpha) = w_1\alpha, & \mathrm{a}(\alpha) = \displaystyle \frac{1-\alpha}{1+(k-1)\alpha}.
	\end{eqnarray}	
\end{itemize}

\subsection{Numerical results}

The experiment to be modeled consists of a homogeneous cylindrical domain of $L=0.2$ meters length and radius $R= 0.06$ meters. For the lower base of the cylinder we will consider that the displacement is null in all directions, while in the upper base we will consider controlled displacement in the axis $z$, that is, $u_z=-0.005t$ meters. In all these examples we will consider free condition for the displacement in the lateral boundary, that is, $\sigma\cdot n=0$. For the case of damage, the boundary conditions will be considered homogeneous Neumann type, that is, $\frac{\partial\alpha}{\partial n}=0$. Finally, we consider the following values for the material parameters
\begin{equation}
\begin{array}{c c c c c c}
E =2.9 \cdot 10^{10}[Pa], & \nu = 0.3, & p=4,  & k=2, &\mbox{ and }  & w_1=10^6\left[\frac{N}{m^3}\right].
\end{array}
\end{equation}

\begin{figure}[h!]
\centering
\begin{subfigure}[b]{0.4\textwidth}
	\includegraphics[width=\textwidth]{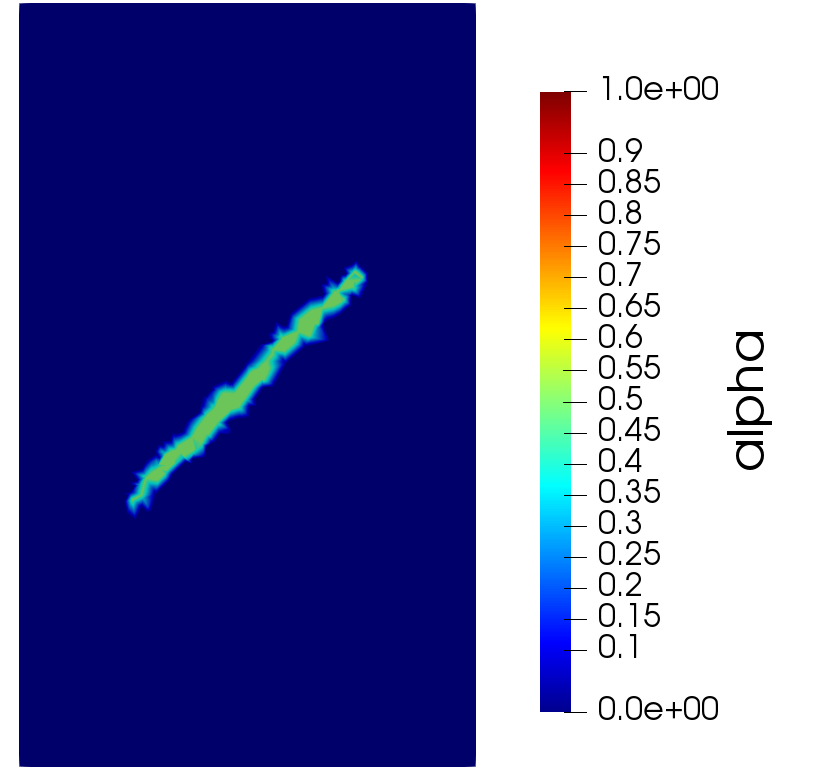}
	\caption{t=0.0 seg.}
\end{subfigure}
\begin{subfigure}[b]{0.4\textwidth}
	\includegraphics[width=\textwidth]{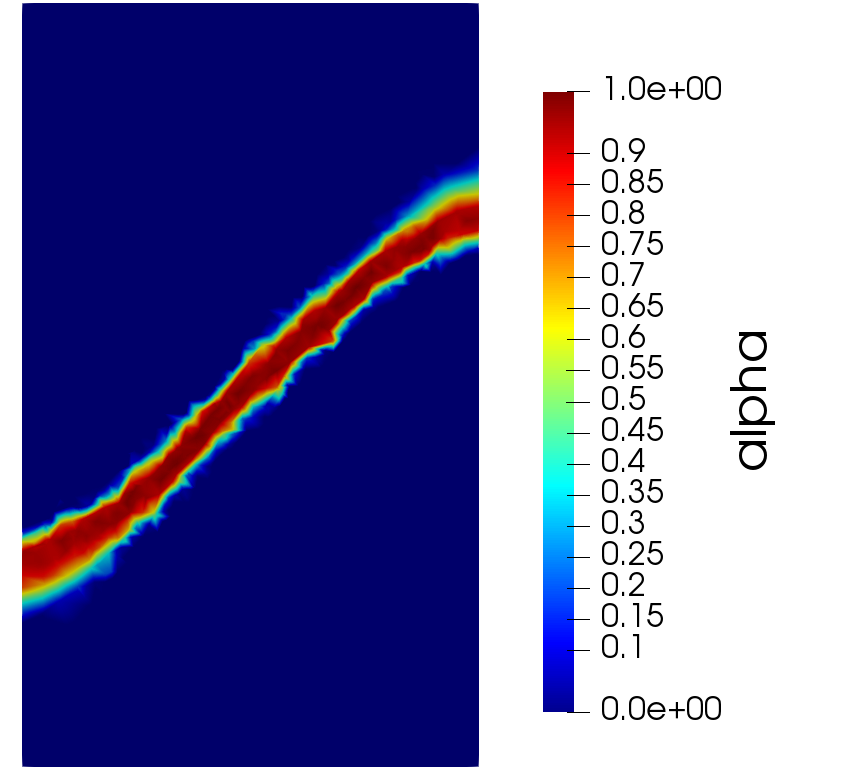}
	\caption{t=1.0 seg.}
\end{subfigure}
\caption{Damage field distribution in the cylinder for the model 1.}\label{cyl_mod1}
\end{figure}

\begin{figure}[h!]
\centering
\begin{subfigure}[b]{0.4\textwidth}
	\includegraphics[width=\textwidth]{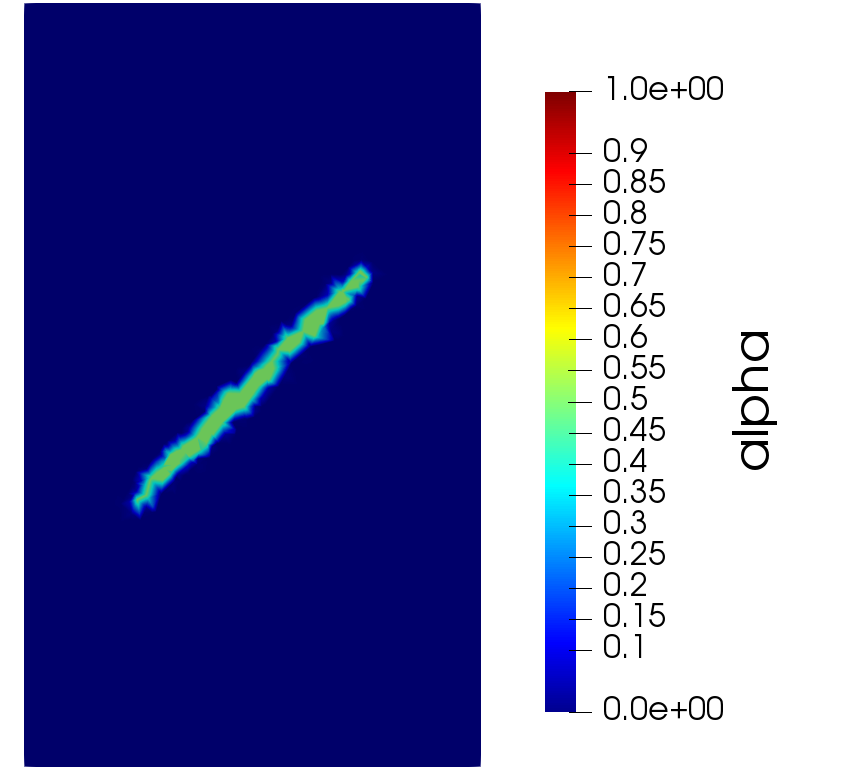}
	\caption{t=0.0 seg.}
\end{subfigure}
\begin{subfigure}[b]{0.4\textwidth}
	\includegraphics[width=\textwidth]{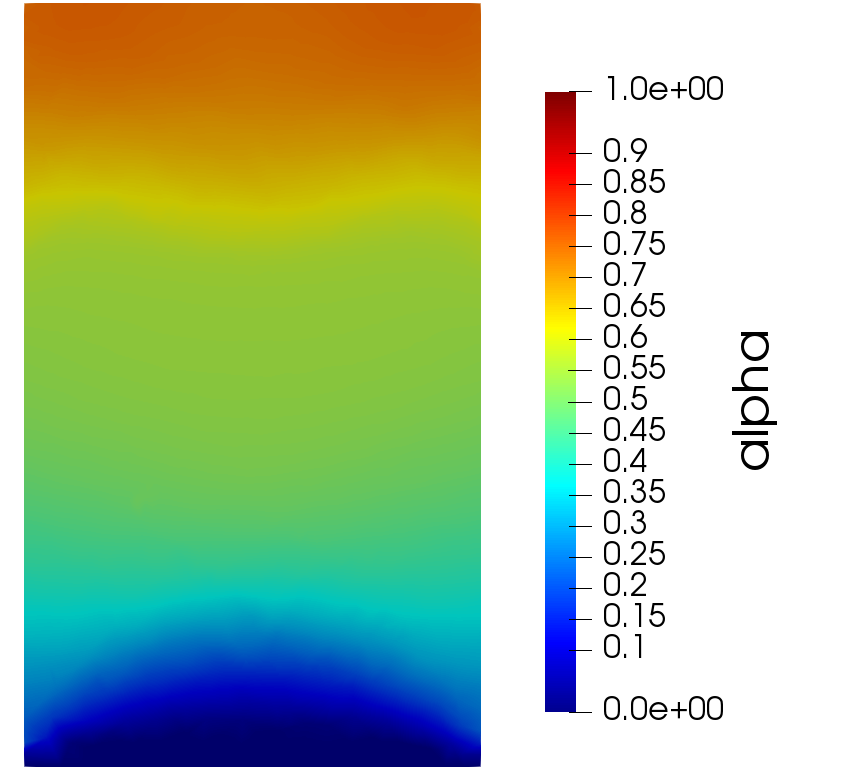}
	\caption{t=1.0 seg.}
\end{subfigure}
\caption{Damage field distribution in the cylinder for the model 2.}\label{cyl_mod2}
\end{figure}

\begin{figure}[h!]
\centering
\begin{subfigure}[b]{0.4\textwidth}
	\includegraphics[width=\textwidth]{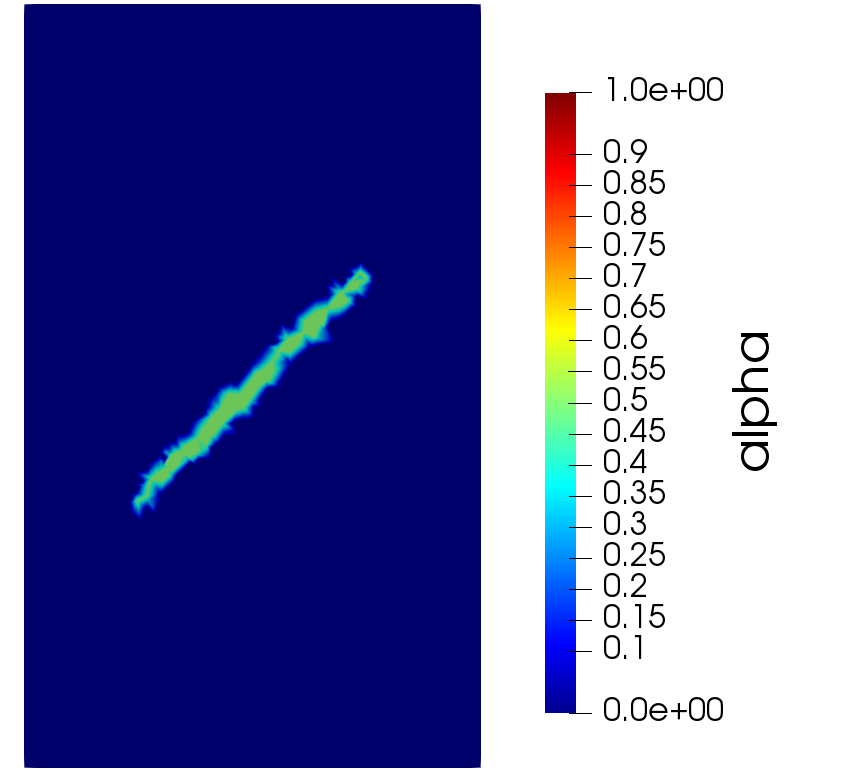}
	\caption{t=0.0 seg.}
\end{subfigure}
\begin{subfigure}[b]{0.4\textwidth}
	\includegraphics[width=\textwidth]{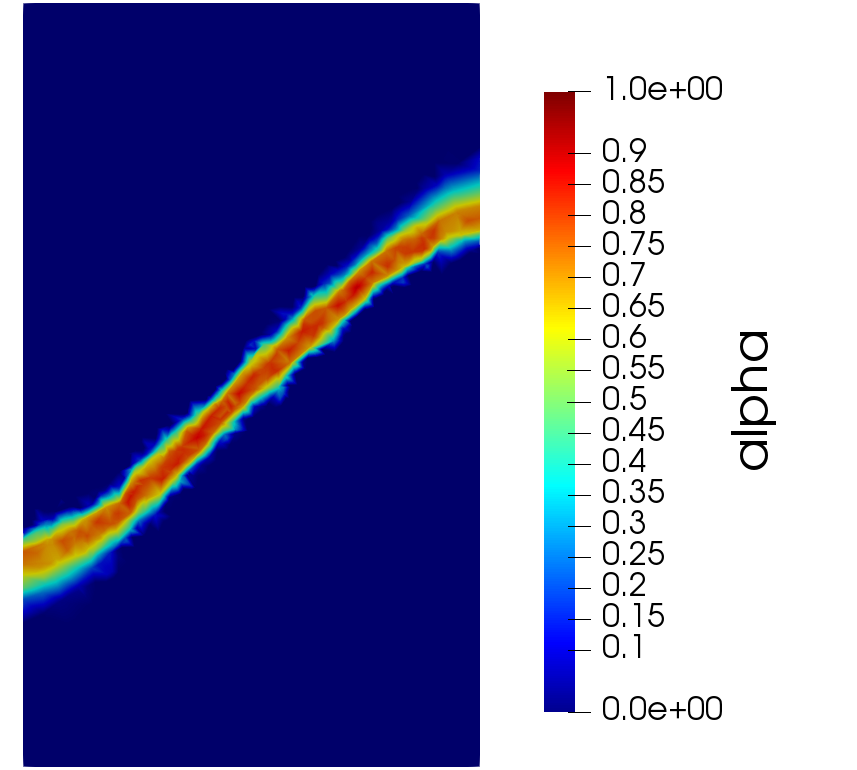}
	\caption{t=1.0 seg.}
\end{subfigure}
\caption{Damage field distribution in the cylinder for the model 3.}\label{cyl_mod3}
\end{figure}

\begin{figure}[h!]
\centering
\begin{subfigure}[b]{0.4\textwidth}
	\includegraphics[width=\textwidth]{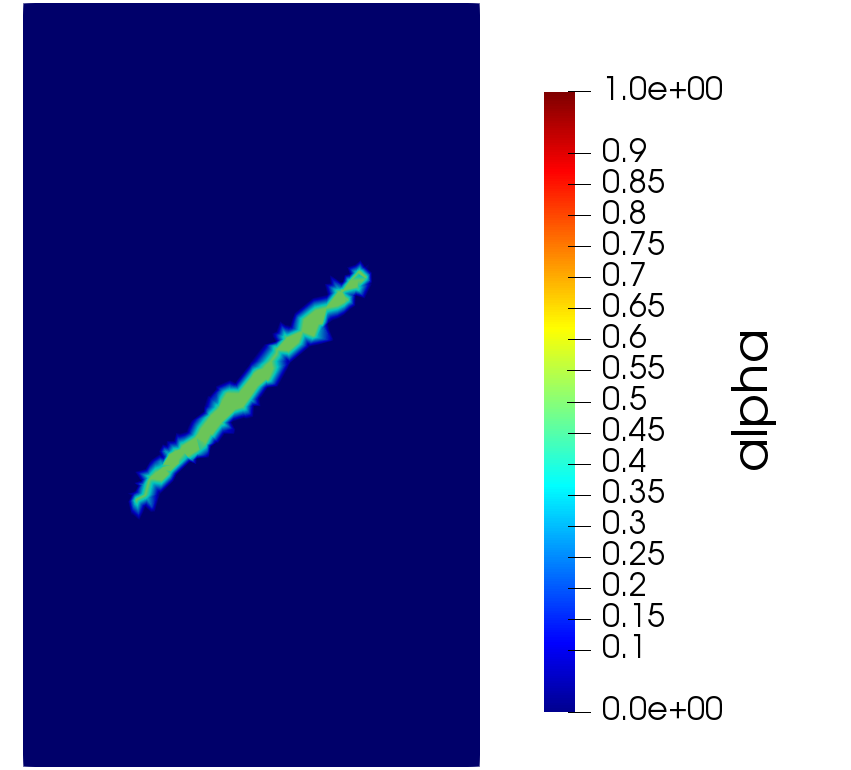}
	\caption{t=0.0 seg.}
\end{subfigure}
\begin{subfigure}[b]{0.4\textwidth}
	\includegraphics[width=\textwidth]{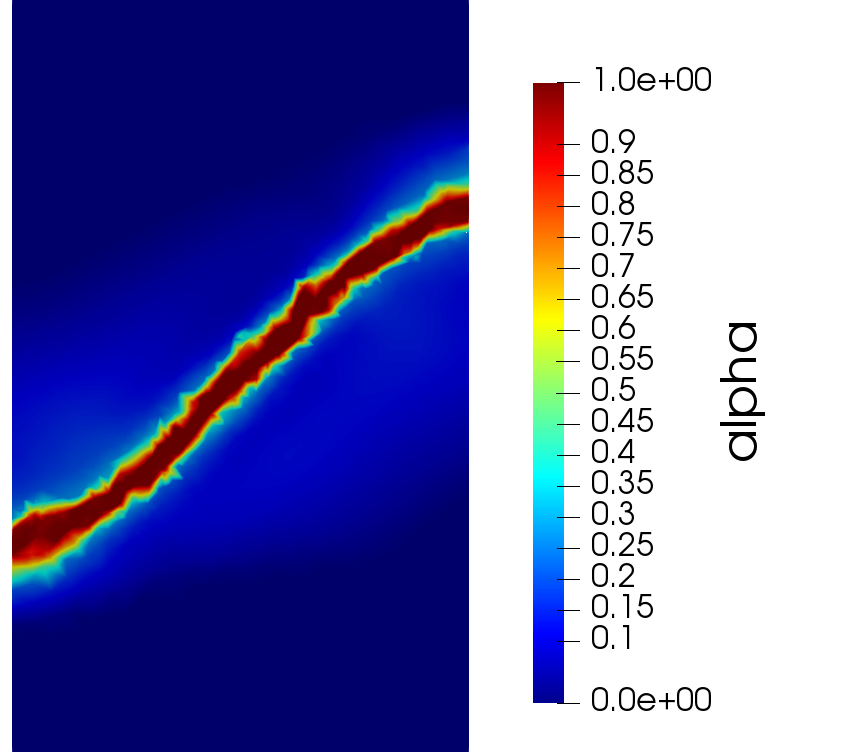}
	\caption{t=1.0 seg.}
\end{subfigure}
\caption{Damage field distribution in the cylinder for the model 4.}\label{cyl_mod4}
\end{figure}

The Figures \ref{cyl_mod1}-\ref{cyl_mod4} display the distribution of the damage for the different models proposed previously and for the times $t=0.0$ seconds and $t=1.0$ seconds with an initial diagonal fracture in the center of the cylinder with magnitude $\alpha =0.5$. Is possible see that, the different models have different behaviour for the distribution of the damage. In the Models 1,3 and 4 the damage only appear in the path of the diagonal fracture, where the values of the damage grows attaining close values of $\alpha=1.0$ (or even $\alpha=1.0$) and extends following this path. In the Model 2 is possible see that the distribution of the damage does not consider the initial fracture and the damage is distributed more homogeneously throughout the domain.  

Through these results, we can differentiate two kinds of materials: The first kind of materials are those that, from a fracture (damaged area), follow the path of the fracture getting the material to be damaged in a particular area, for example the Models 1, 3 and 4. The second kind of materials, are those that do not consider the fracture and are damaged homogeneously (Model 2). After seeing the results, we can conclude that our model is valid to represent the behavior of fractures, since it is capable of replicating the both classes of materials.

\section{influence of the cavity in the damage model}\label{Cav}

\subsection{Formulation and boundary conditions}\label{FormImp}

For simulate the block caving process, we consider a domain that simulates a rock mass in three dimensions $\Omega_0 \subseteq \mathbb{R}^3$ with $\partial \Omega = \Gamma_{lat}\cup \Gamma_{up} \cup \Gamma_{down}$ denoting the lateral, upper and lower bounded of the domain respectively.  We consider a time depending interior domain $S(t_i) \subseteq \Omega_0$, with $\partial S(t_i) = \Gamma_{cav}(t_i)$ and $S(t_i) \subseteq S(t_{i+1})$, representing the interior cavity of this rock mass. Then in each time step, we consider the evolution domain $\Omega(t_i) = \Omega_0 \setminus S(t_i)$, that is, our domain where we calculate the damage model will be the initial domain $\Omega_0$ minus the internal cavity in each time step.

\begin{figure}[h!]
	\centering
	\includegraphics[scale=0.8]{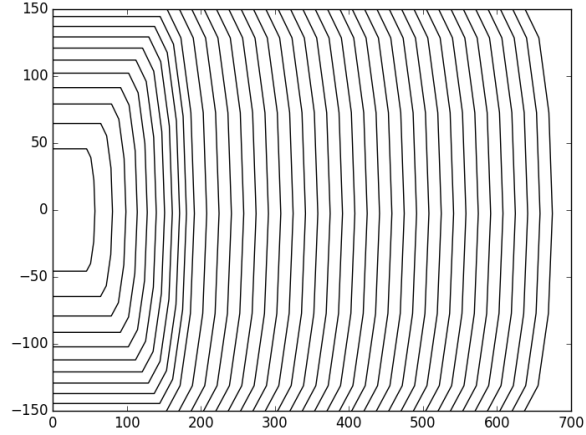} 
	\caption{Interior flat surface that represent the cavity sequences.}
	\label{soc_plane}
\end{figure}

To construct the mesh that represent our domain, first we need define the interior cavity in each time step $t_i$. For this, we consider a flat surface that represent the base of the cavity, which is projected on the $z$ axis. More specifically, the base is repeated and scaling at different levels of $z$ in the range $z \in [0,z_{max}]$ to build the cavity surface in 3D. Figure \ref{soc_plane} displays the forty interior flat surface that later represent the cavity sequences that we use in our tests.

\begin{figure}[h!]
	\centering
	\includegraphics[width=0.7\textwidth]{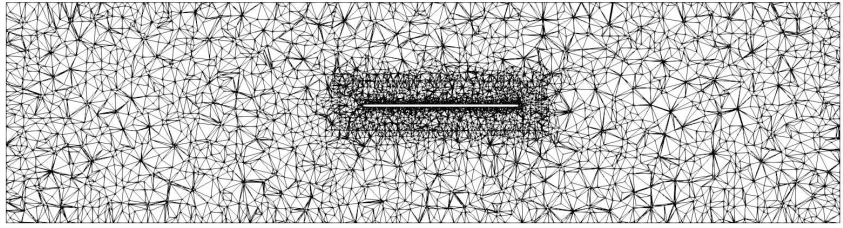} 
	\caption{Lateral cut of the mesh in the last time step.}
	\label{mesh_lat}
\end{figure}

Once the cavity in each time step is known, for each cavity,  a rectangular box mesh is built around it where near the boundary of the cavity a more refine mesh is defined in order to better reflect the effects produced by each advance of the cavity. Figure \ref{mesh_lat} shows a lateral view of the mesh in the last cavity advance, where we can be observed the greater refinement around the cavity, for our simulations.

To define the boundary conditions for the elasticity equation in our problems, let us consider the  equilibrium equation for fixed value of $\alpha$ that, no loss of generality, we assume that $\alpha=0$
\begin{equation}
\begin{array}{r c l r l}
-\mathrm{div}(\sigma) & = & f, & \mbox{ in } & \Omega,\\
\sigma \cdot n & = & 0, & \mbox{ on } & \Gamma_{cav},\\
\sigma \cdot n & = & 0, & \mbox{ on } & \Gamma_{up},\\
u \cdot n & = & 0, & \mbox{ on } & \Gamma_{down},\\
\sigma \cdot n & = & g, & \mbox{ on } & \Gamma_{lat},
\end{array}
\end{equation}
where $\sigma= 2 \mu \varepsilon(u) + \lambda \mathrm{tr}(\varepsilon(u))I$ and $f = (0,0,- \rho g)^T$. We will assume that, in the rock mass, we do not have lateral displacement, that is, 
\begin{equation}
u= \left(\begin{array}{c} 0\\ 0 \\u_{z}(z) \end{array}\right) \mbox{ and } \varepsilon(u) = \left(\begin{array}{c c c} 0 & 0 & 0\\ 0 & 0 & 0 \\0 & 0 & u'_{z}(z) \end{array}\right), 
\end{equation}
so we will have
\begin{equation}
\sigma = \left(\begin{array}{c c c} \lambda & 0 & 0\\ 0 & \lambda & 0 \\0 & 0 & \lambda + 2 \mu  \end{array}\right) u'_{z}(z) \mbox{ and } \mathrm{div}(\sigma)=\left(\begin{array}{c} 0\\ 0 \\(\lambda+2\mu)u''_{z}(z) \end{array}\right)=\left(\begin{array}{c} 0\\ 0 \\\rho g \end{array}\right).
\end{equation}
Then
\begin{equation}
u'_{z}(z) = \frac{ \rho g}{\lambda + 2 \mu} z + C. 
\end{equation}
Now, denoting by $H_{max}$, the maximum height of the domain and considering that $\sigma_{zz}(H_{max})=0$, we have
\begin{equation}
C = - \frac{\rho g H_{max}}{\lambda + 2 \mu }.
\end{equation}
Then
\begin{equation}
u'_{z}(z) = \frac{ \rho g}{\lambda + 2 \mu} (z-H_{max}),
\end{equation}
and, thus
\begin{equation}
\sigma= \left(\begin{array}{c c c} \frac{\lambda}{\lambda + 2 \mu} & 0 & 0\\ 0 & \frac{\lambda}{\lambda + 2 \mu} & 0 \\0 & 0 & 1 \end{array}\right)\rho g  (z-H_{max}).
\end{equation} 

This gives us a condition that depends on the height $z$ and the maximum height of the domain $H_{max}$, so for non-homogeneous domains, for example a real mine, this condition gives us problems, because the border conditions have different values for each side face, since these would not have the same height, causing that the damage appears in the corners of the domain. To solve this problem we will include a Robin boundary condition of the form $\sigma \cdot n + \overline{k}(u \cdot n) \cdot n =0$, with $\overline{k}>0$ a real constant. Finally the boundary condition is of the form
\begin{equation}
g=\frac{\lambda}{\lambda + 2 \mu}\rho g (z-H_{max}) - \overline{k} (u\cdot n) \cdot n.
\end{equation}

\subsection{Numerical results}

In all our tests, we consider $\Omega_0 \in \mathbb{R}^3$ such that $x\in [-1540,2060]$, $y \in [-1050,1050]$ and $z \in [ -500,450]$. Here the meshes are defined following the Section \ref{FormImp} and we use the synthetic cavities defined in the Figures \ref{soc_plane}, where it is possible to see the different configurations of each cavity. All simulations presented here use the following values for the parameters in the damage model
\begin{equation}
\begin{array}{c c c c c c c}
E =2.9 \cdot 10^{10}[Pa], & \nu = 0.3, & p=4, & k=2, & w_1=10^6 \left[\frac{N}{m^3}\right]& \mbox{ and }  & \overline{k}=10^{9}[Pa],
\end{array}
\end{equation}

For the numerical test, first, we consider that the constant in the function $w(\alpha)$ is independent of the constant in the elliptic regularization term in (\ref{moddamageeq}), this in order to independently control the variable the energy dissipated without affecting the term of regularization. Thus, for the models shown in Section \ref{Models}, we take
\begin{equation}
\begin{array}{c c c c c}
w(\alpha)=w_{11}\alpha, & w(\alpha)=w_{11}\alpha^2, & w(\alpha)=w_{11}\left(1-(1-\alpha)^{p/2}\right), & \mbox{ with }& w_{11}\neq w_1.
\end{array}
\end{equation}

\begin{figure}[h!]
	\begin{subfigure}[b]{0.49\textwidth}
		\includegraphics[width=\textwidth]{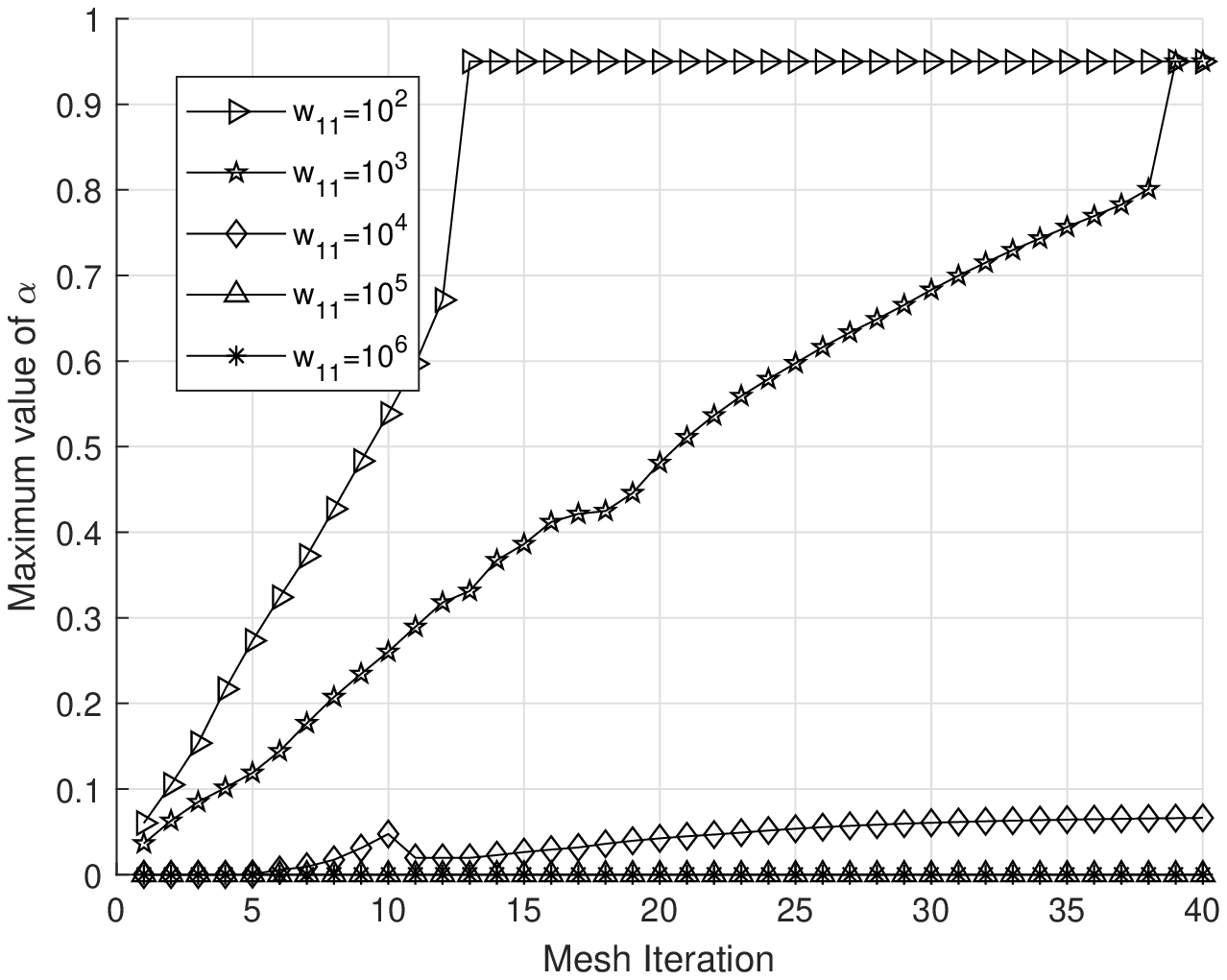}
		\caption{Model 1.}
	\end{subfigure}
		\begin{subfigure}[b]{0.49\textwidth}
		\includegraphics[width=\textwidth]{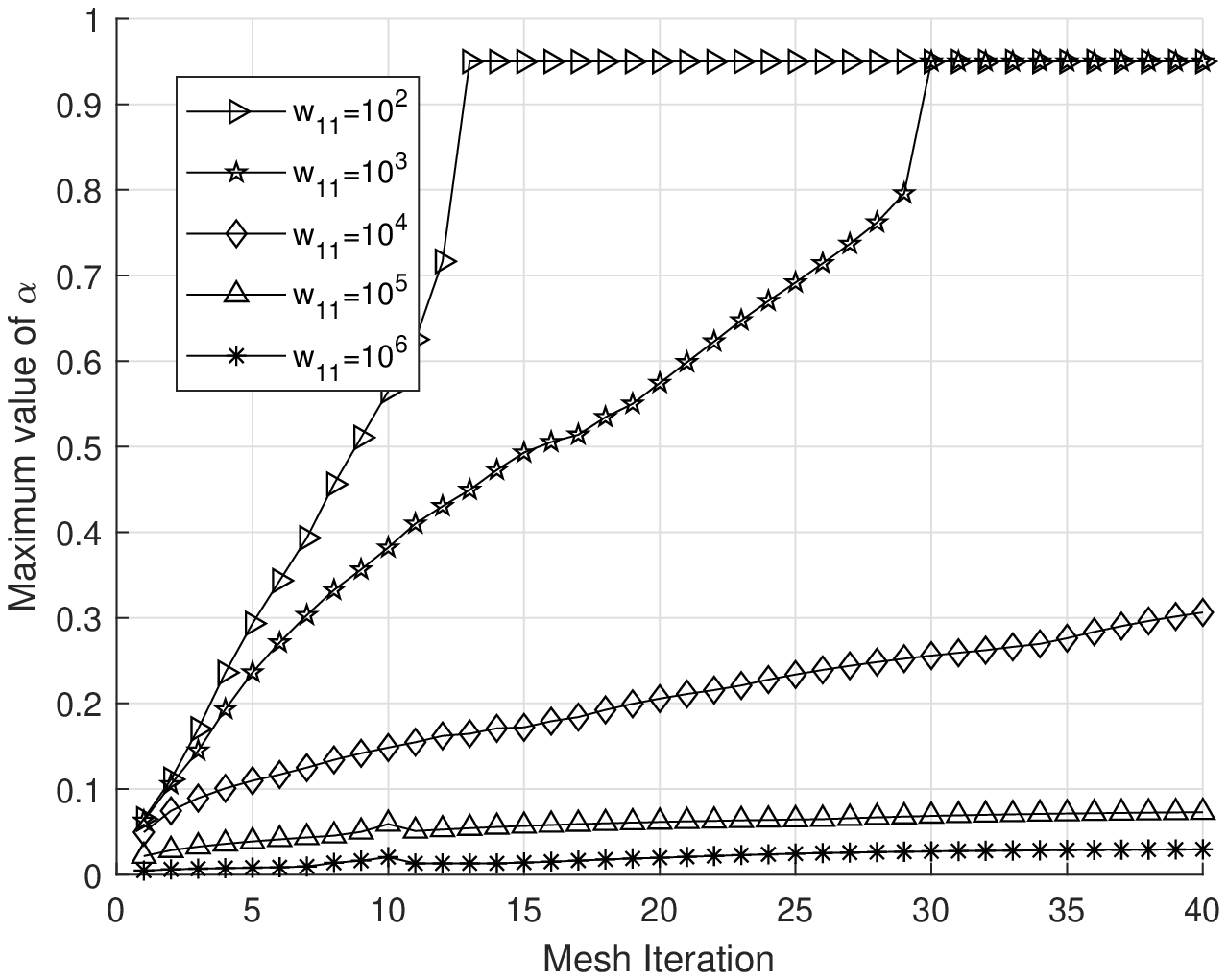}
		\caption{Model 2.}
	\end{subfigure}
		\begin{subfigure}[b]{0.49\textwidth}
		\includegraphics[width=\textwidth]{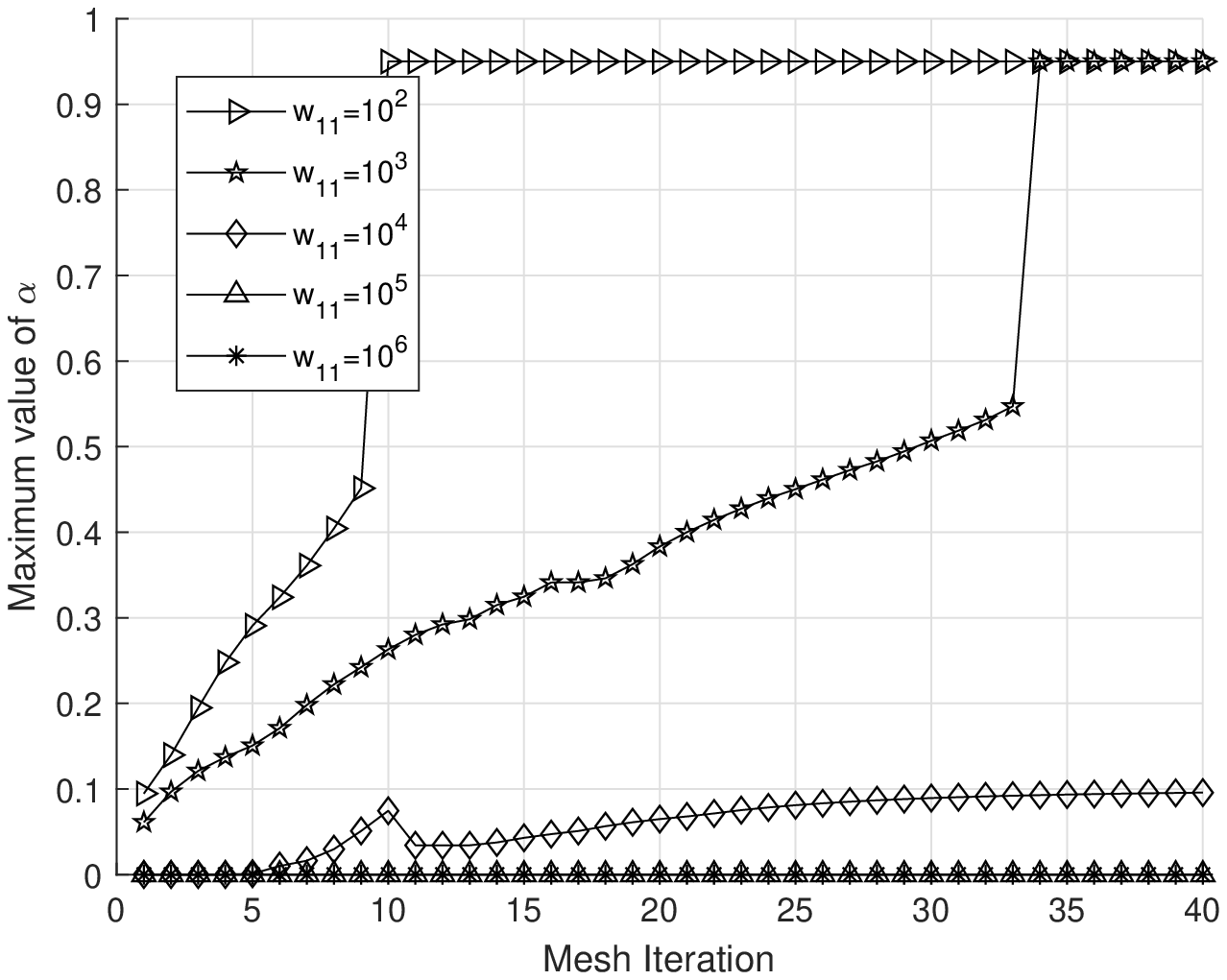}
		\caption{Model 3.}
	\end{subfigure}
	\begin{subfigure}[b]{0.49\textwidth}
		\includegraphics[width=\textwidth]{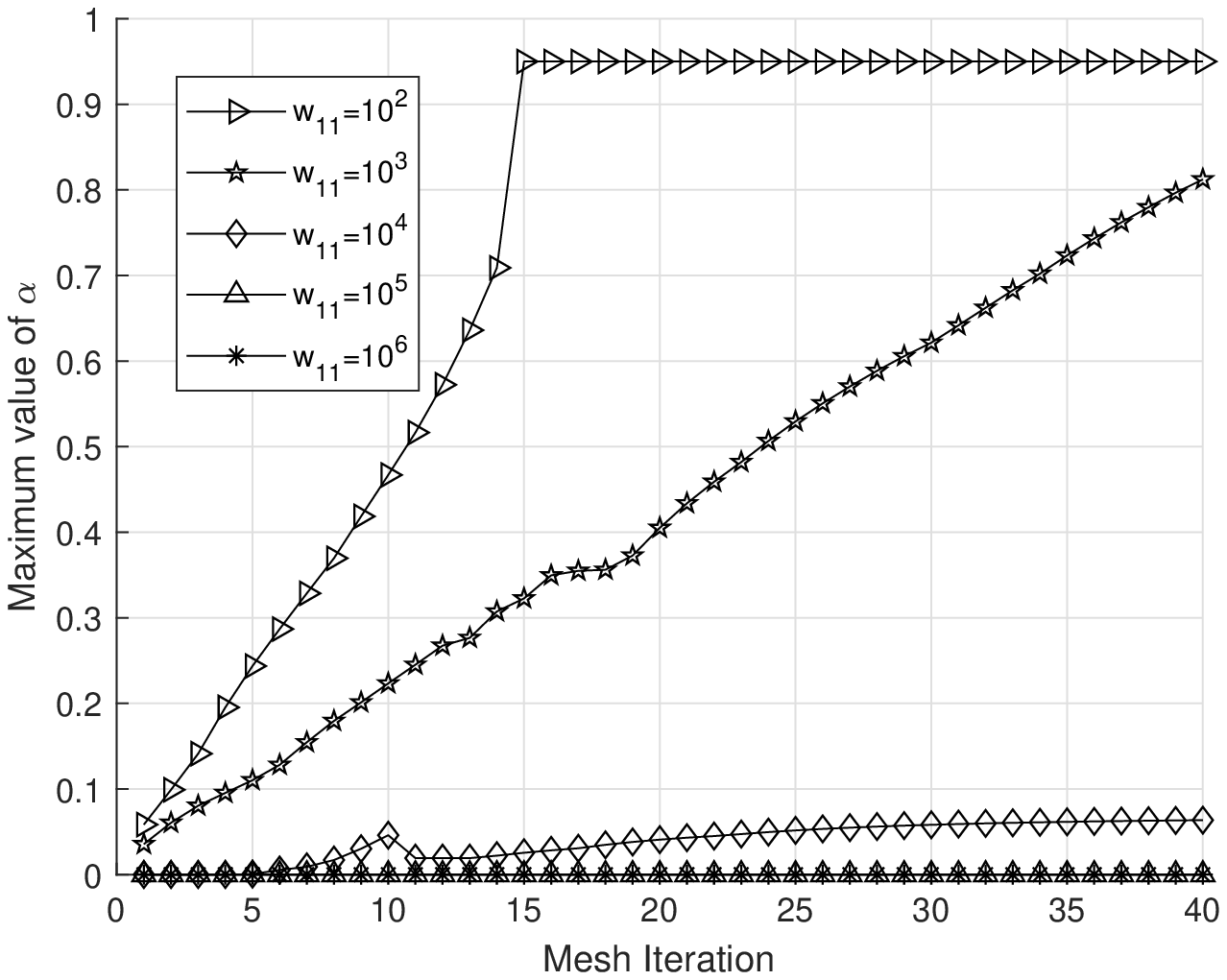}
		\caption{Model 4.}
	\end{subfigure}
	\caption{Maximum damage value for each cavity advance for different values of $w_{11}$ and $K_w$.}\label{alphamaxmods}
\end{figure}

The Figures \ref{alphamaxmods} summarize the maximum values of $\alpha$ in each mesh iteration. The images display the evolution of this maximum value and we can be see that for little values $w_{11}$ the damage grows faster than bigger values of this constant. For large $w_{11}$ values the  damage does not achieve the maximum value to fully damage state.

\begin{figure}[h!]
	\centering
	\begin{subfigure}[b]{0.49\textwidth}
		\centering
		\includegraphics[width=\textwidth]{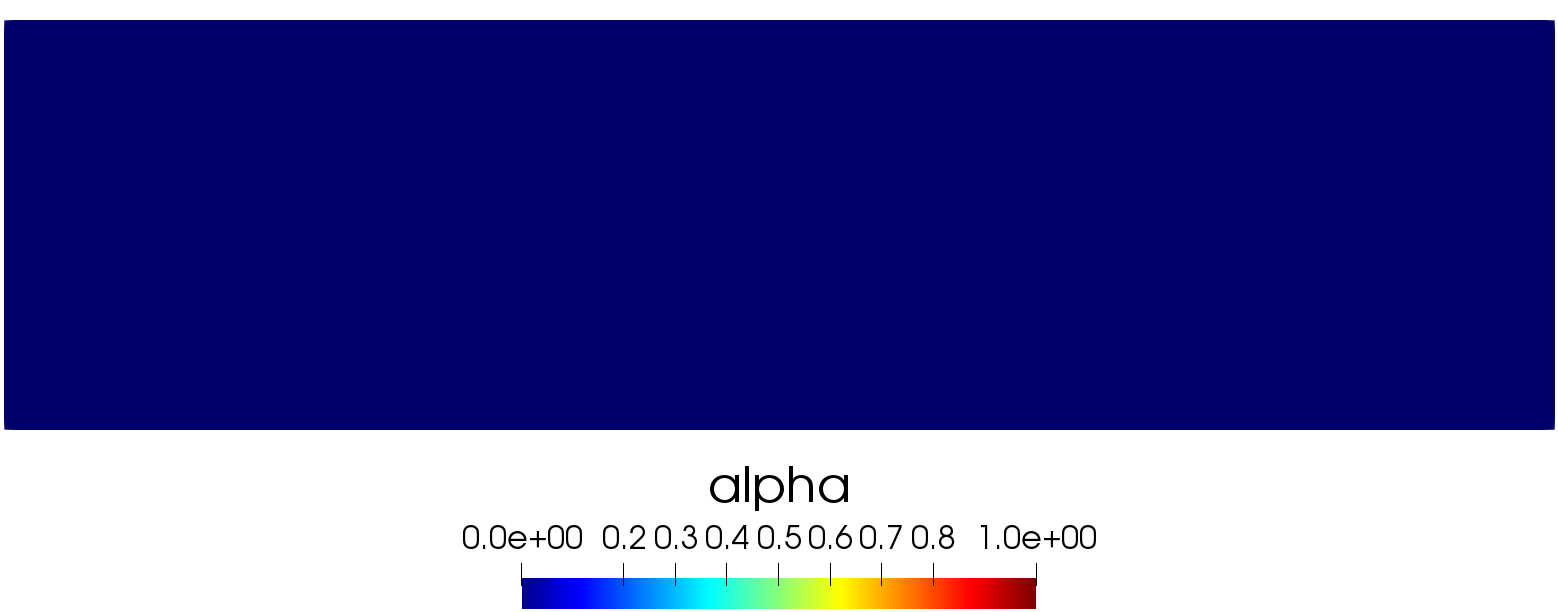}
		\caption{$\Omega(t_{0})$}
	\end{subfigure}
	\begin{subfigure}[b]{0.49\textwidth}
		\centering
		\includegraphics[width=\textwidth]{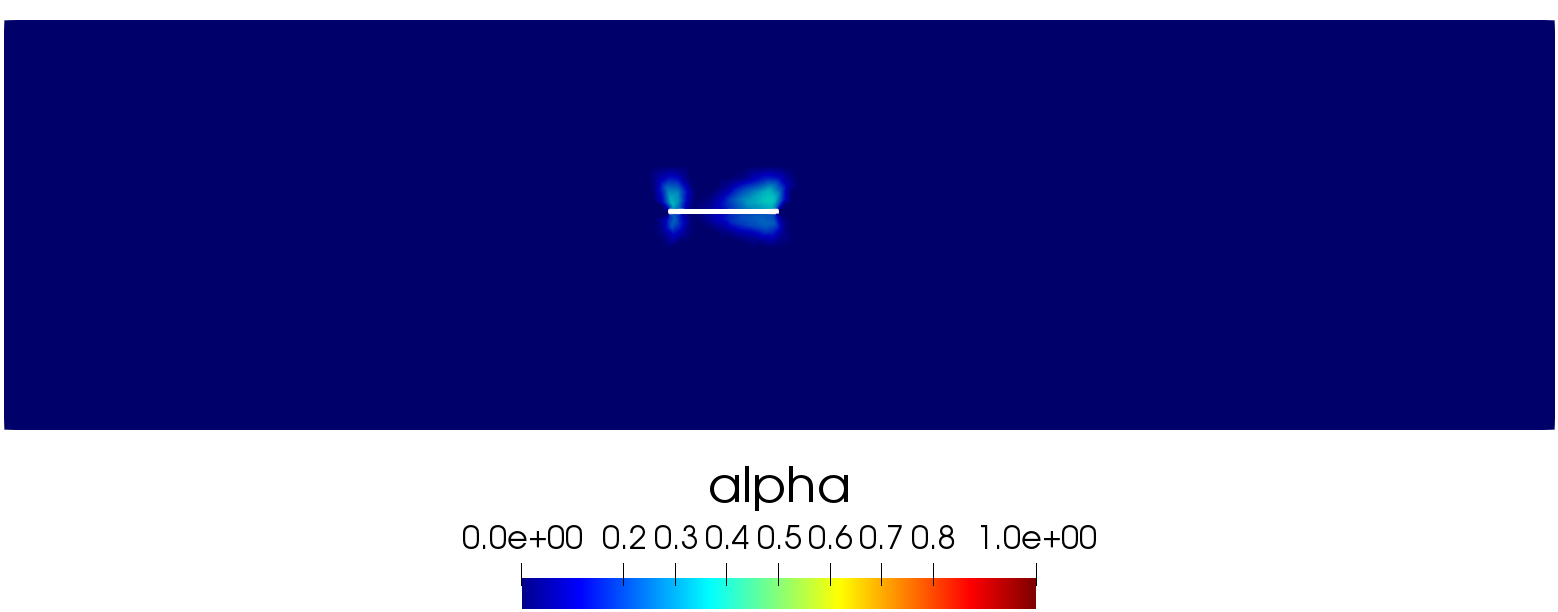}
		\caption{$\Omega(t_{15})$}
	\end{subfigure}
	\begin{subfigure}[b]{0.49\textwidth}
		\centering
		\includegraphics[width=\textwidth]{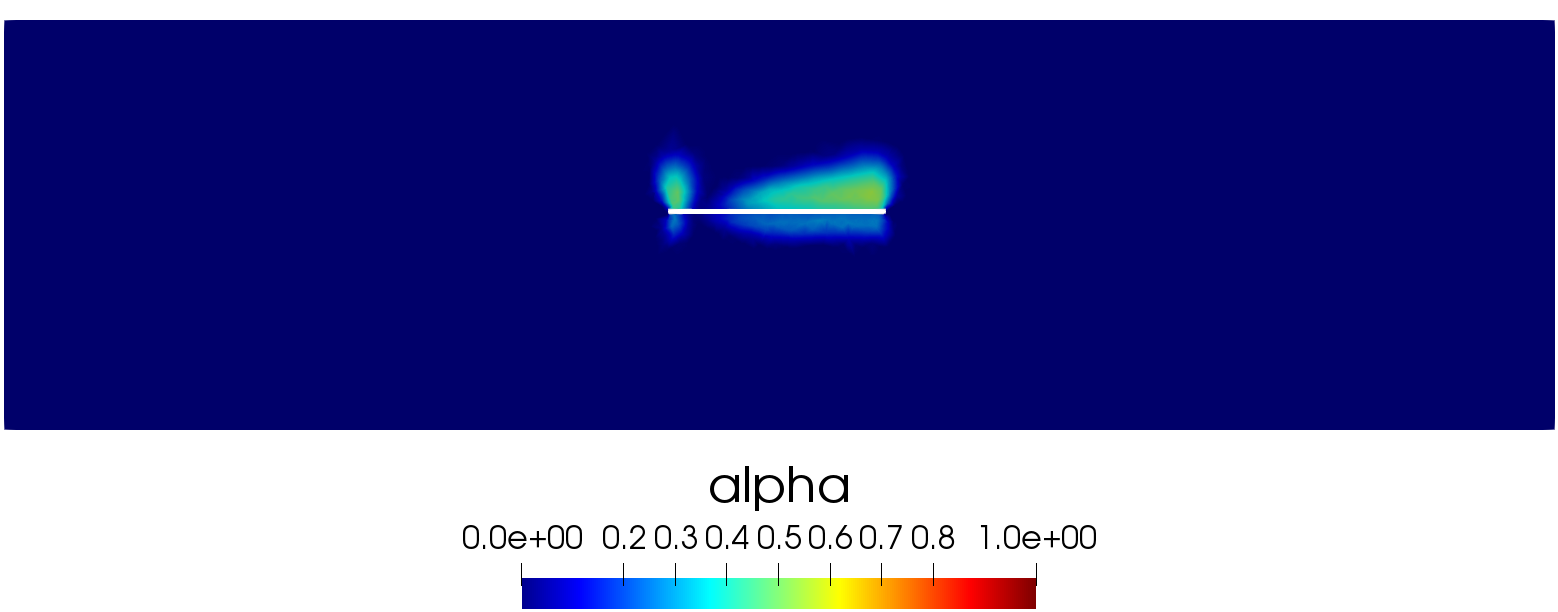}
		\caption{$\Omega(t_{30})$}
	\end{subfigure}	
	\begin{subfigure}[b]{0.49\textwidth}
		\centering
		\includegraphics[width=\textwidth]{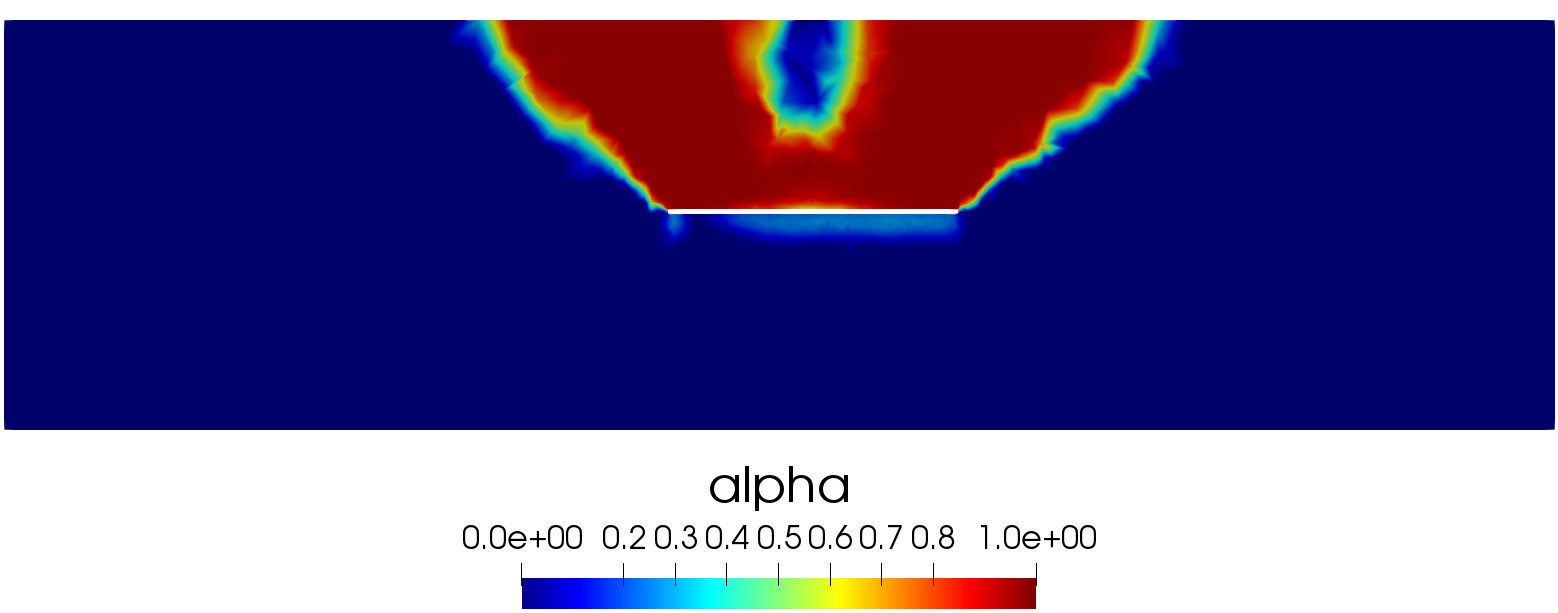}
		\caption{$\Omega(t_{40})$}
	\end{subfigure}
	\caption{Damage field distribution in the rock mass for Model 1.}\label{mod1latcut}
\end{figure}

\begin{figure}[h!]
	\centering
	\begin{subfigure}[b]{0.49\textwidth}
		\centering
		\includegraphics[width=\textwidth]{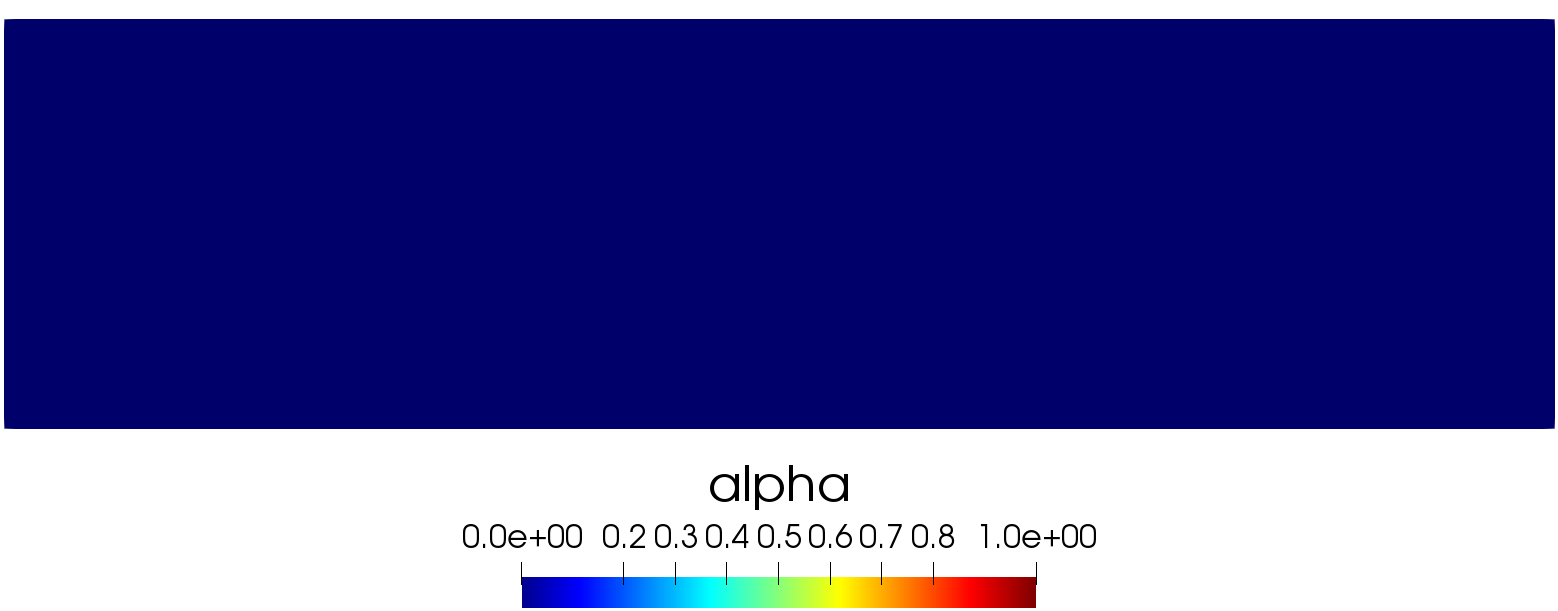}
		\caption{$\Omega(t_{0})$}
	\end{subfigure}
	\begin{subfigure}[b]{0.49\textwidth}
		\centering
		\includegraphics[width=\textwidth]{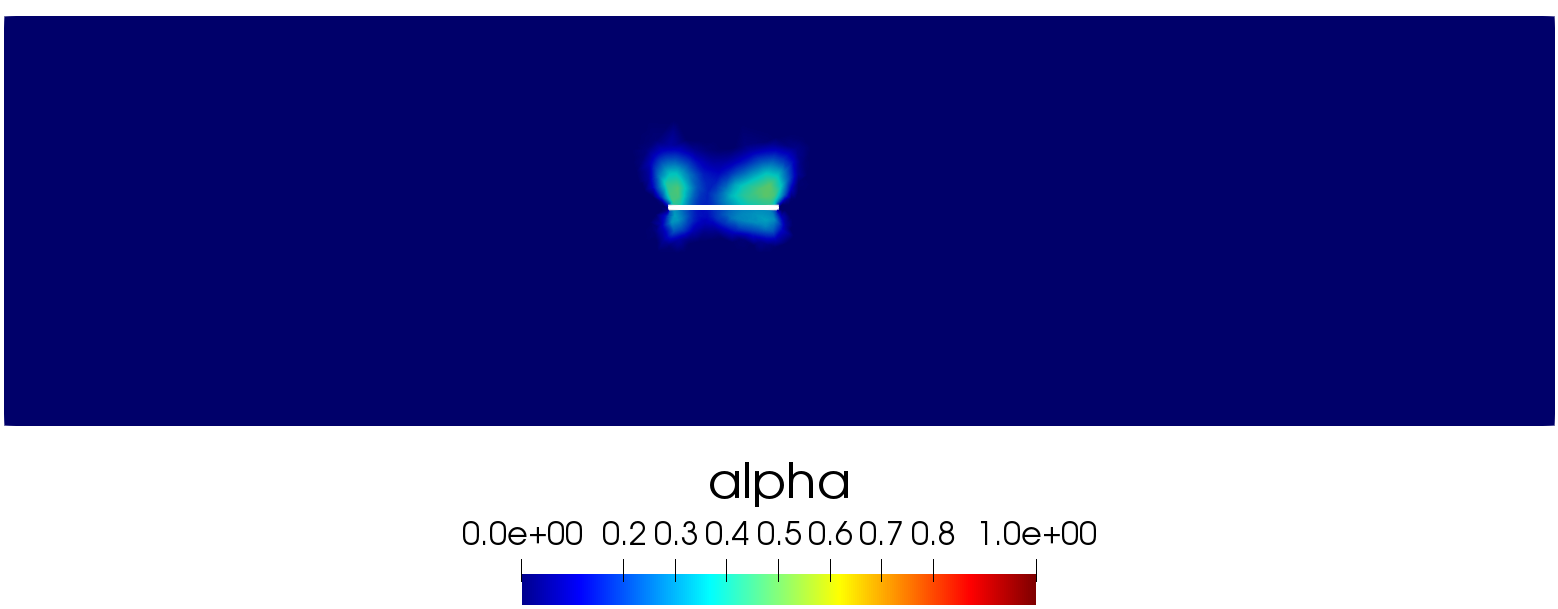}
		\caption{$\Omega(t_{15})$}
	\end{subfigure}
	\begin{subfigure}[b]{0.49\textwidth}
		\centering
		\includegraphics[width=\textwidth]{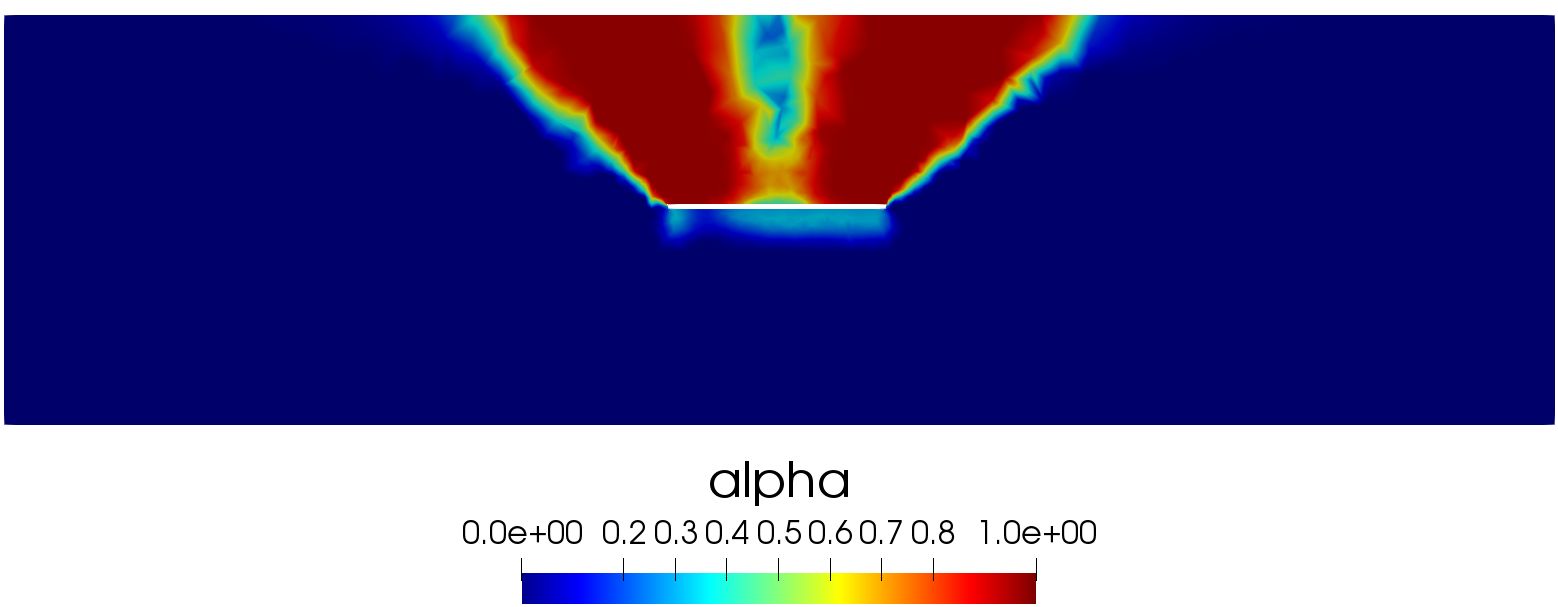}
		\caption{$\Omega(t_{30})$}
	\end{subfigure}
	\begin{subfigure}[b]{0.49\textwidth}
		\centering
		\includegraphics[width=\textwidth]{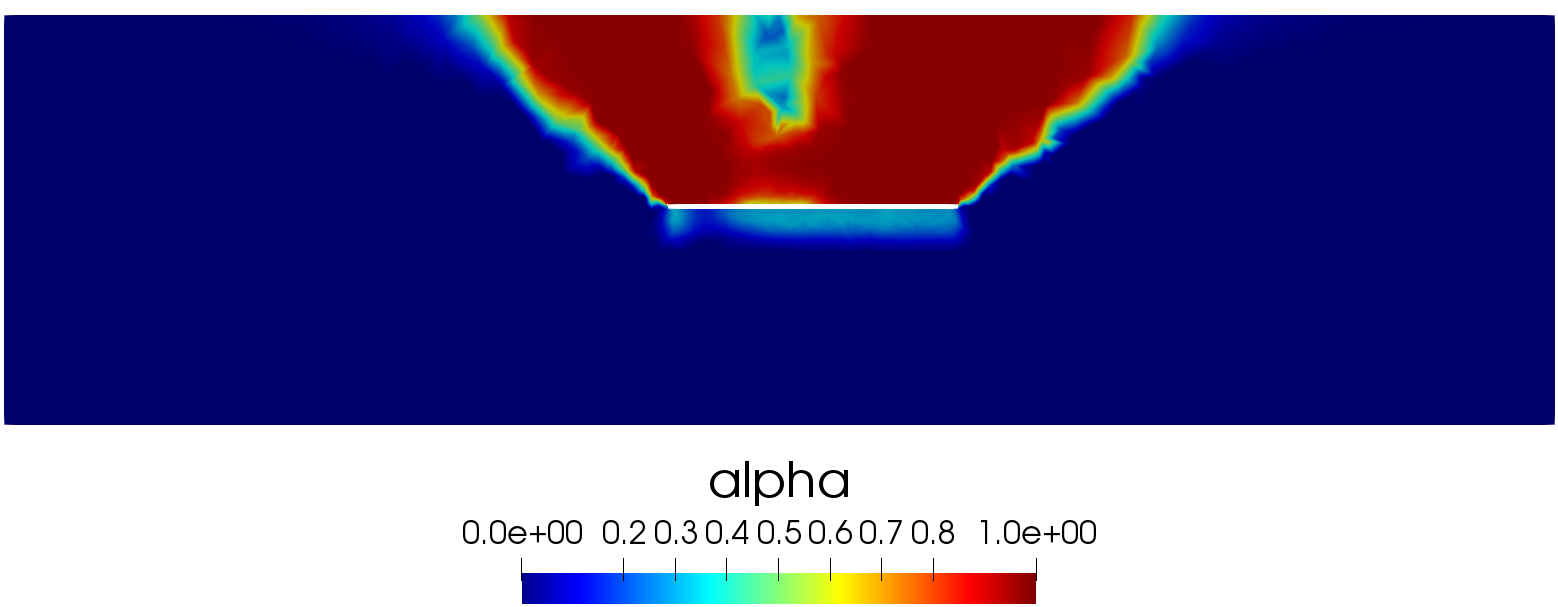}
		\caption{$\Omega(t_{40})$}
	\end{subfigure}
	\caption{Damage field distribution in the rock mass for Model 2.}\label{mod2latcut}
\end{figure}

\begin{figure}[h!]
	\centering
	\begin{subfigure}[b]{0.49\textwidth}
		\centering
		\includegraphics[width=\textwidth]{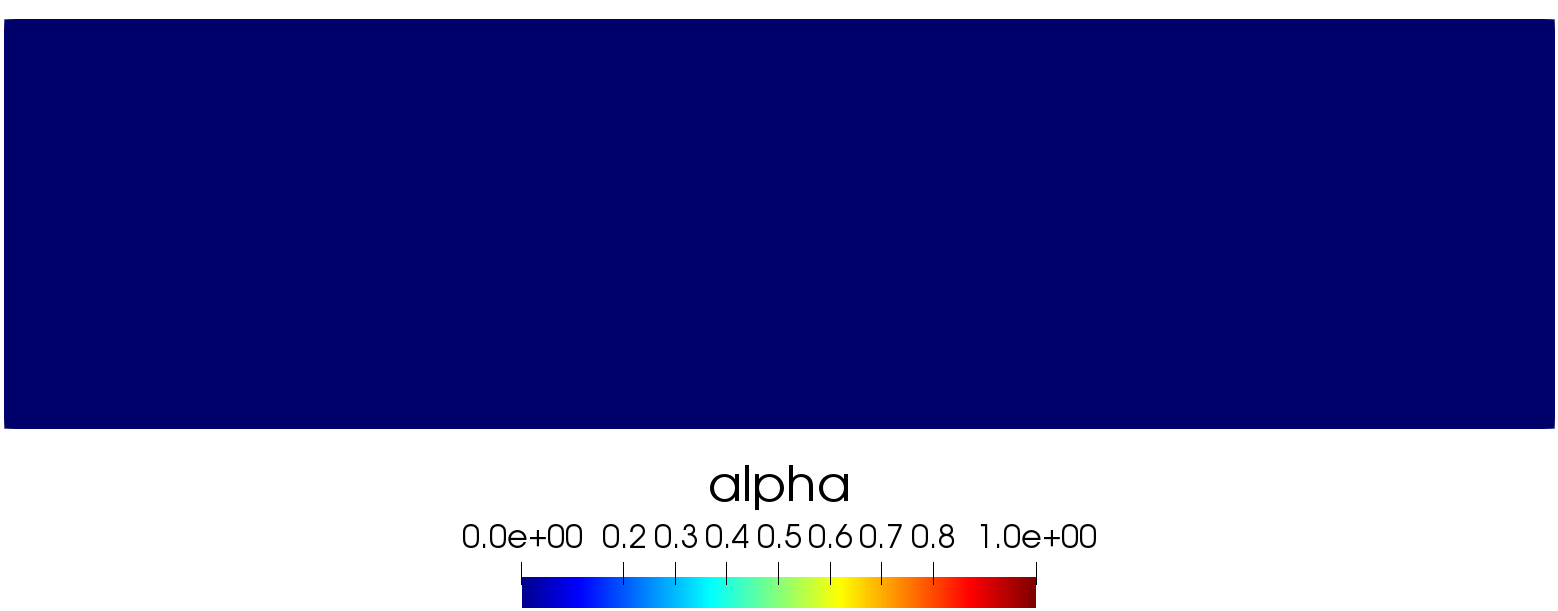}
		\caption{$\Omega(t_{0})$}
	\end{subfigure}
	\begin{subfigure}[b]{0.49\textwidth}
		\centering
		\includegraphics[width=\textwidth]{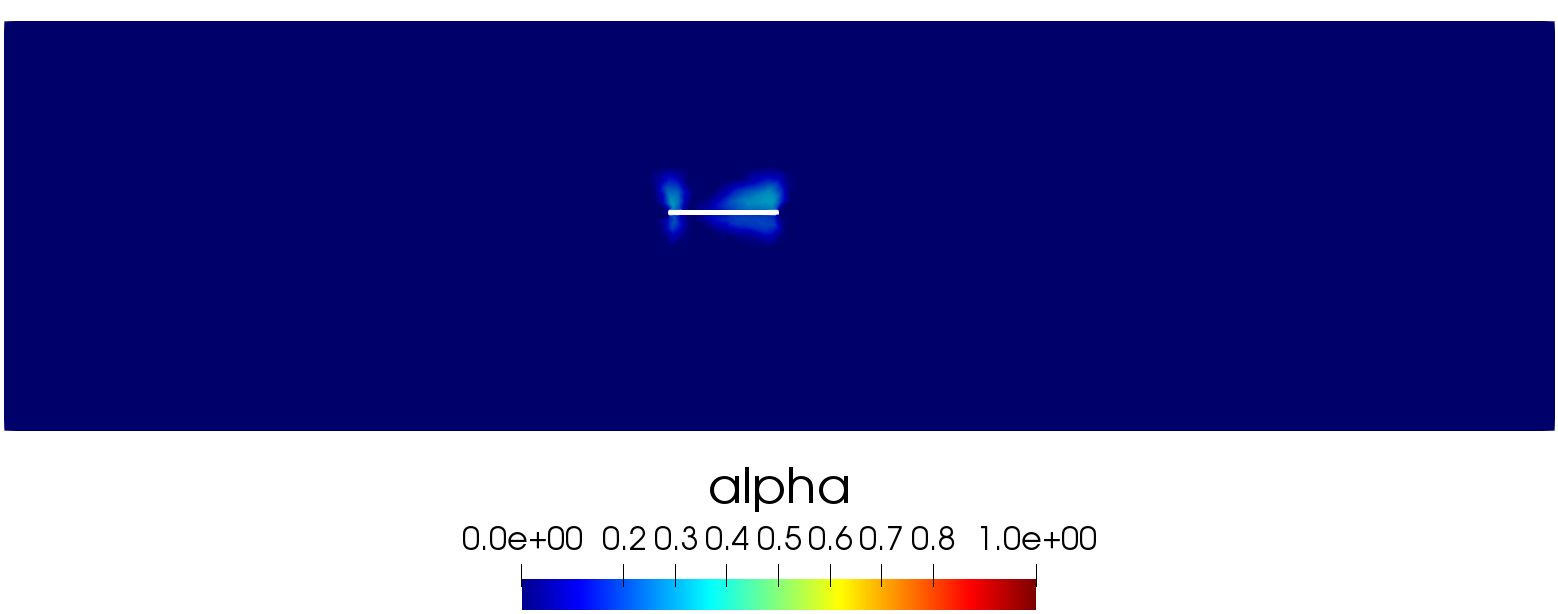}
		\caption{$\Omega(t_{15})$}
	\end{subfigure}
	\begin{subfigure}[b]{0.49\textwidth}
		\centering
		\includegraphics[width=\textwidth]{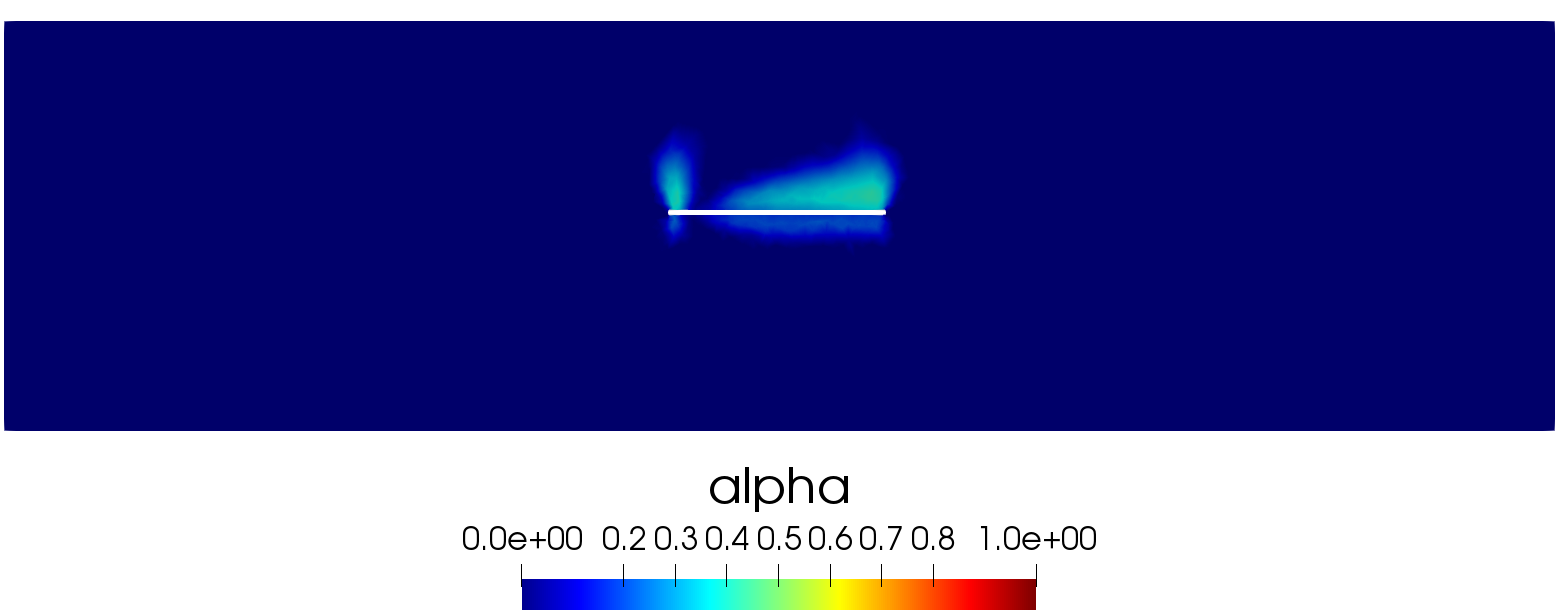}
		\caption{$\Omega(t_{30})$}
	\end{subfigure}	
	\begin{subfigure}[b]{0.49\textwidth}
		\centering
		\includegraphics[width=\textwidth]{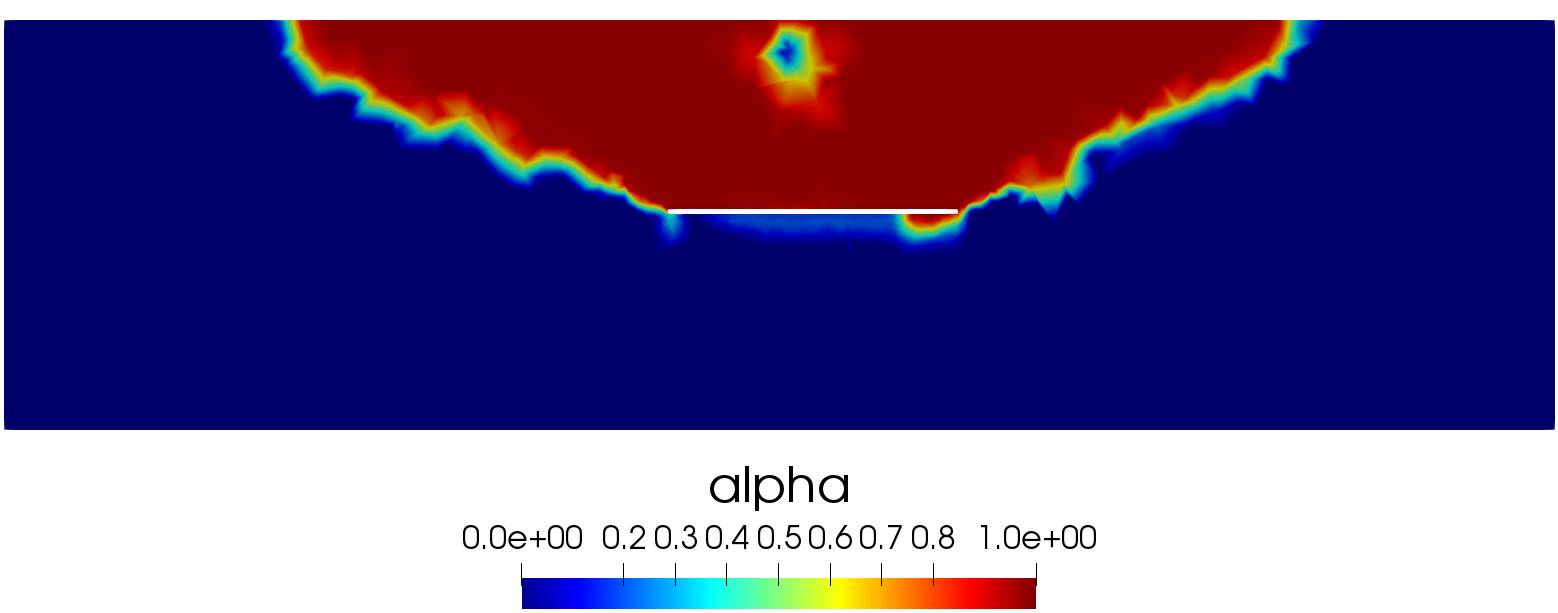}
		\caption{$\Omega(t_{40})$}
	\end{subfigure}
	\caption{Damage field distribution in the rock mass for Model 3.}\label{mod3latcut}
\end{figure}

\begin{figure}[h!]
	\centering
	\begin{subfigure}[b]{0.49\textwidth}
		\centering
		\includegraphics[width=\textwidth]{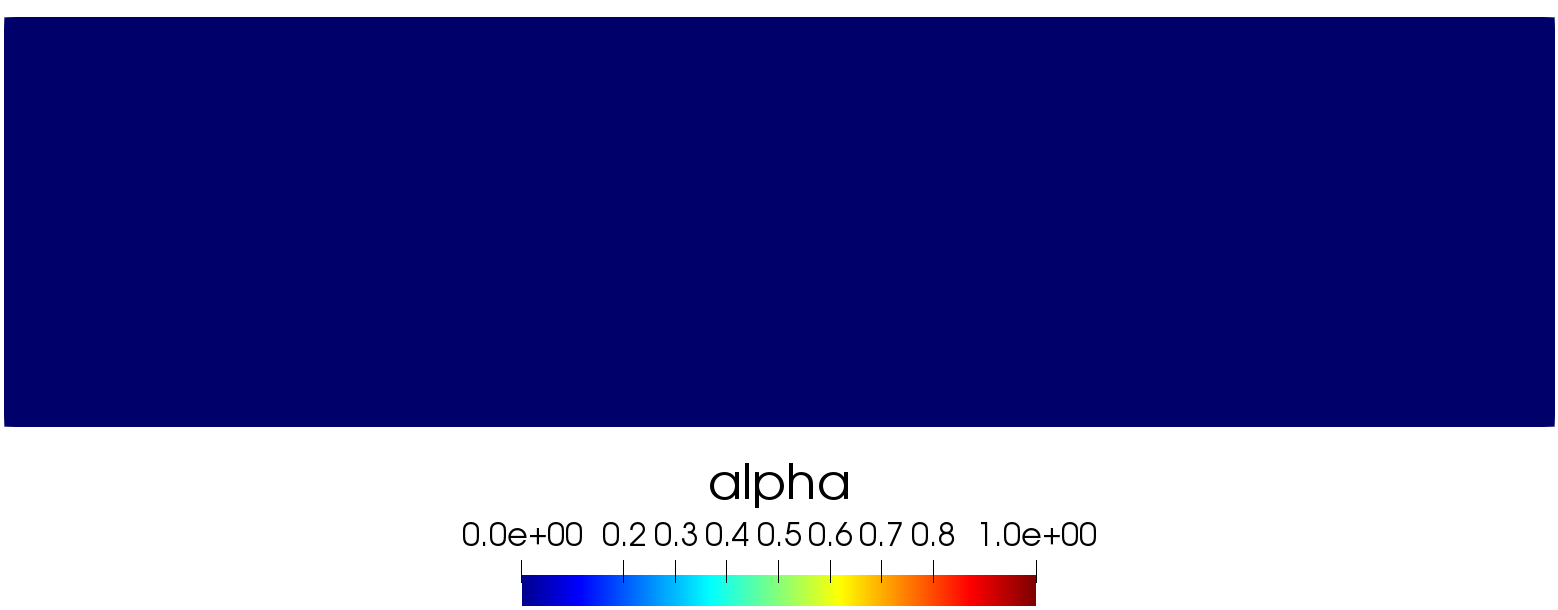}
		\caption{$\Omega(t_{0})$}
	\end{subfigure}
	\begin{subfigure}[b]{0.49\textwidth}
		\centering
		\includegraphics[width=\textwidth]{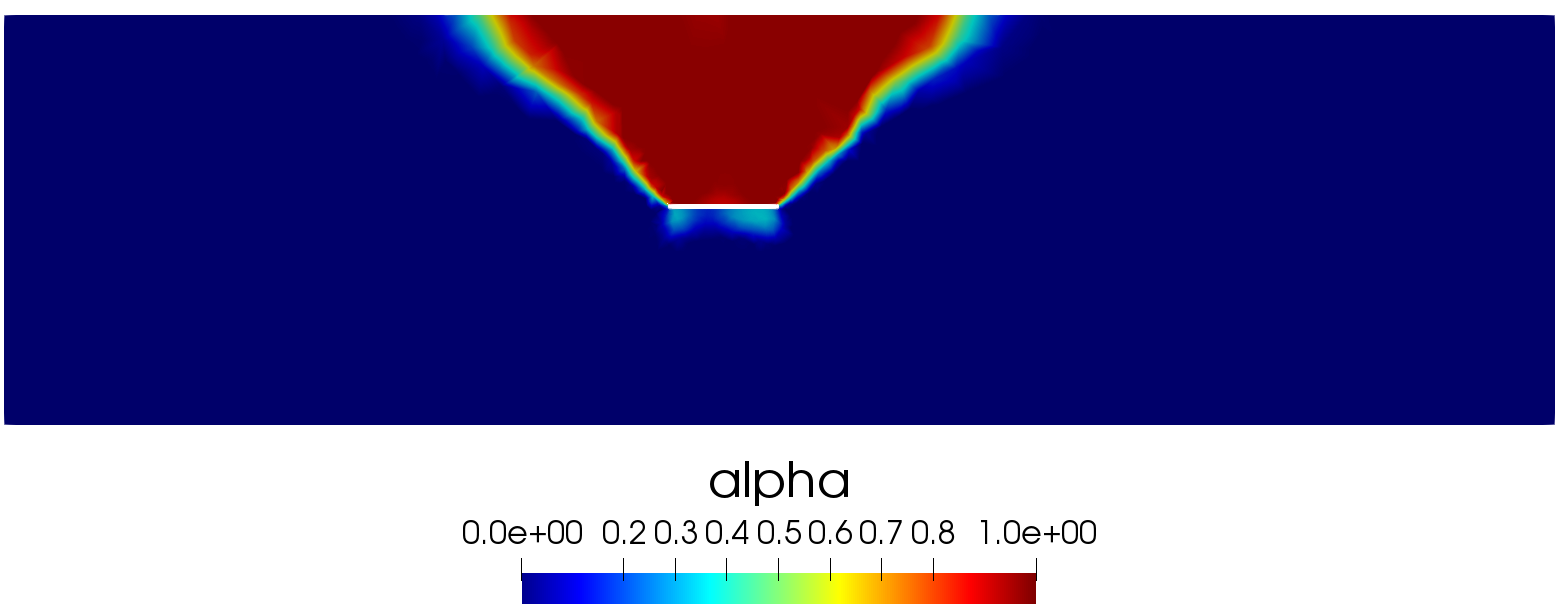}
		\caption{$\Omega(t_{15})$}
	\end{subfigure}
	\begin{subfigure}[b]{0.49\textwidth}
		\centering
		\includegraphics[width=\textwidth]{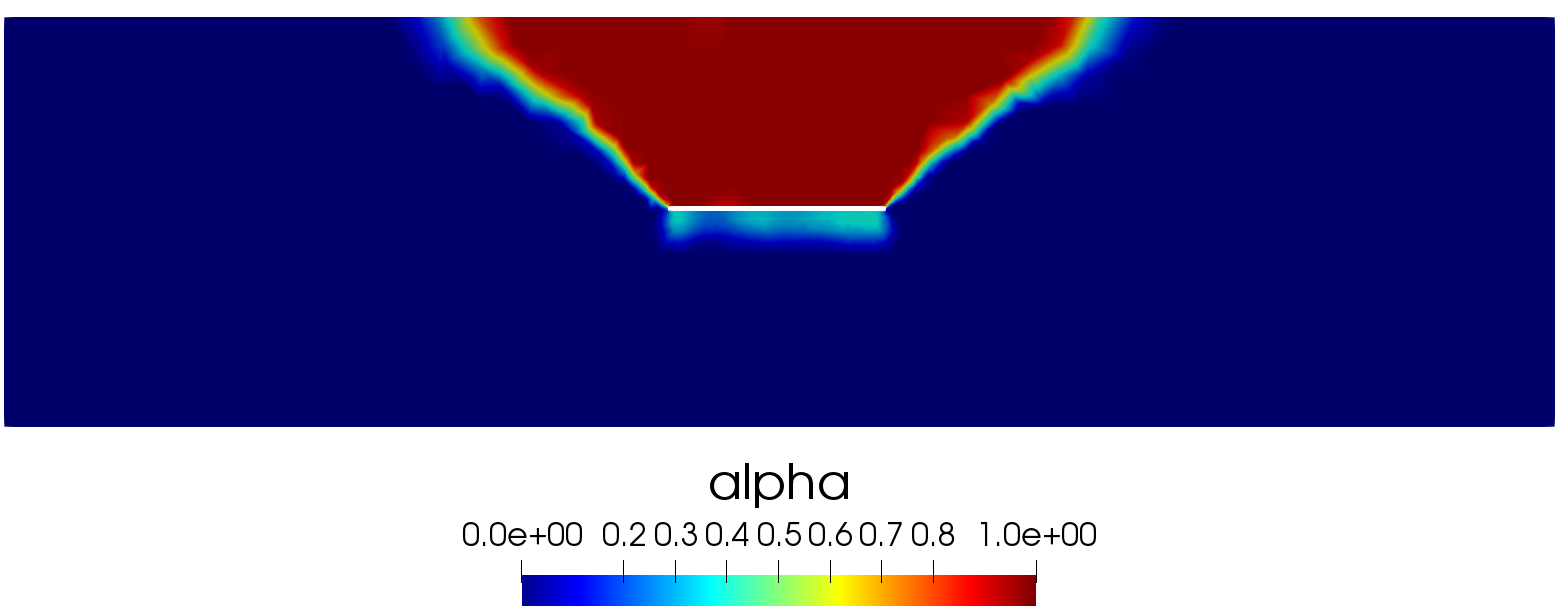}
		\caption{$\Omega(t_{30})$}
	\end{subfigure}	
	\begin{subfigure}[b]{0.49\textwidth}
		\centering
		\includegraphics[width=\textwidth]{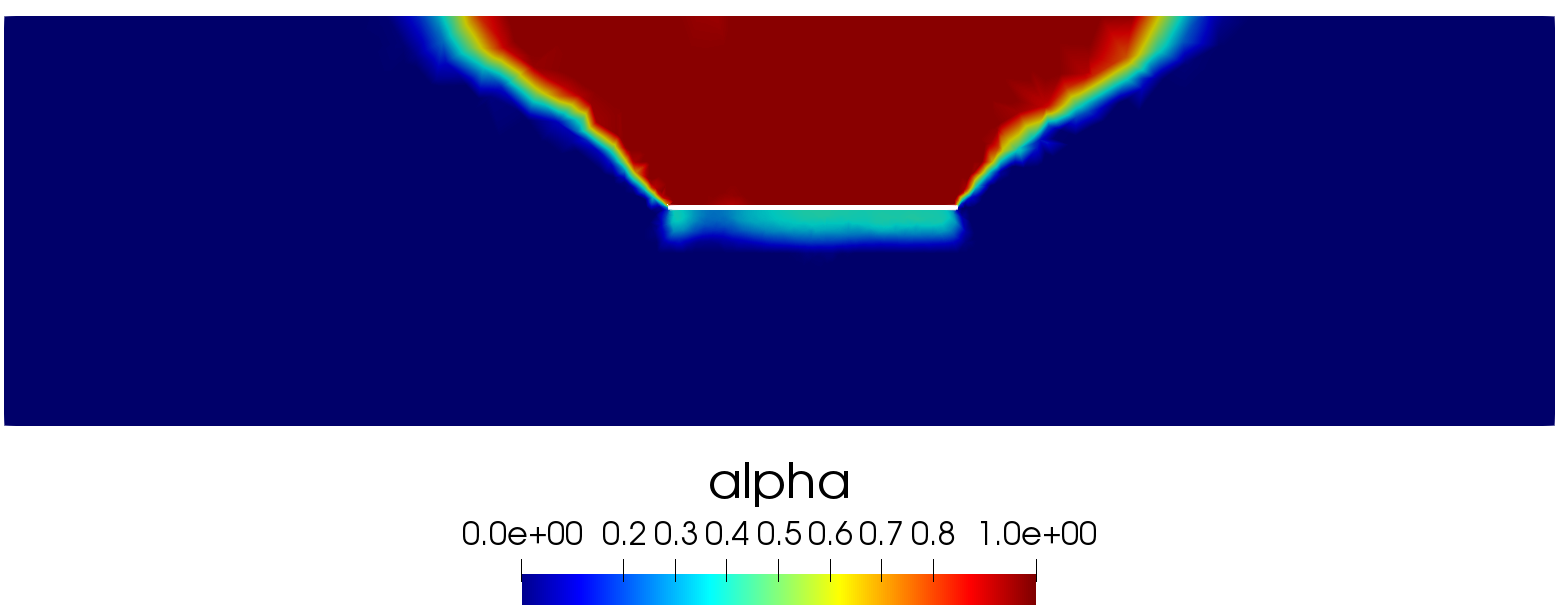}
		\caption{$\Omega(t_{40})$}
	\end{subfigure}
	\caption{Damage field distribution in the rock mass for Model 4.}\label{mod4latcut}
\end{figure}

Following the previous results, we consider that, for the Model 1, 2 and 3, $w_{11}=10^3$ and, for Model 4, $w_{11}=10^2$. Figures \ref{mod1latcut}-\ref{mod4latcut} display the evolution of the damage when the cavity advance in time, for the different models. As it showed in \cite{bonnetier2020shearcompression}, the damage is close to zero almost everywhere except around the cavity, it is possible see that the distribution of the damage reaching the upper boundary of the domain. For the Model 1, 2 and 3, the material is fully damage over the cavity whereas, for the Model 4, the damage does not reach its maximum value on this cavity.

\begin{figure}[h!]
	\centering
	\begin{subfigure}[b]{0.49\textwidth}
		\centering
		\includegraphics[width=\textwidth]{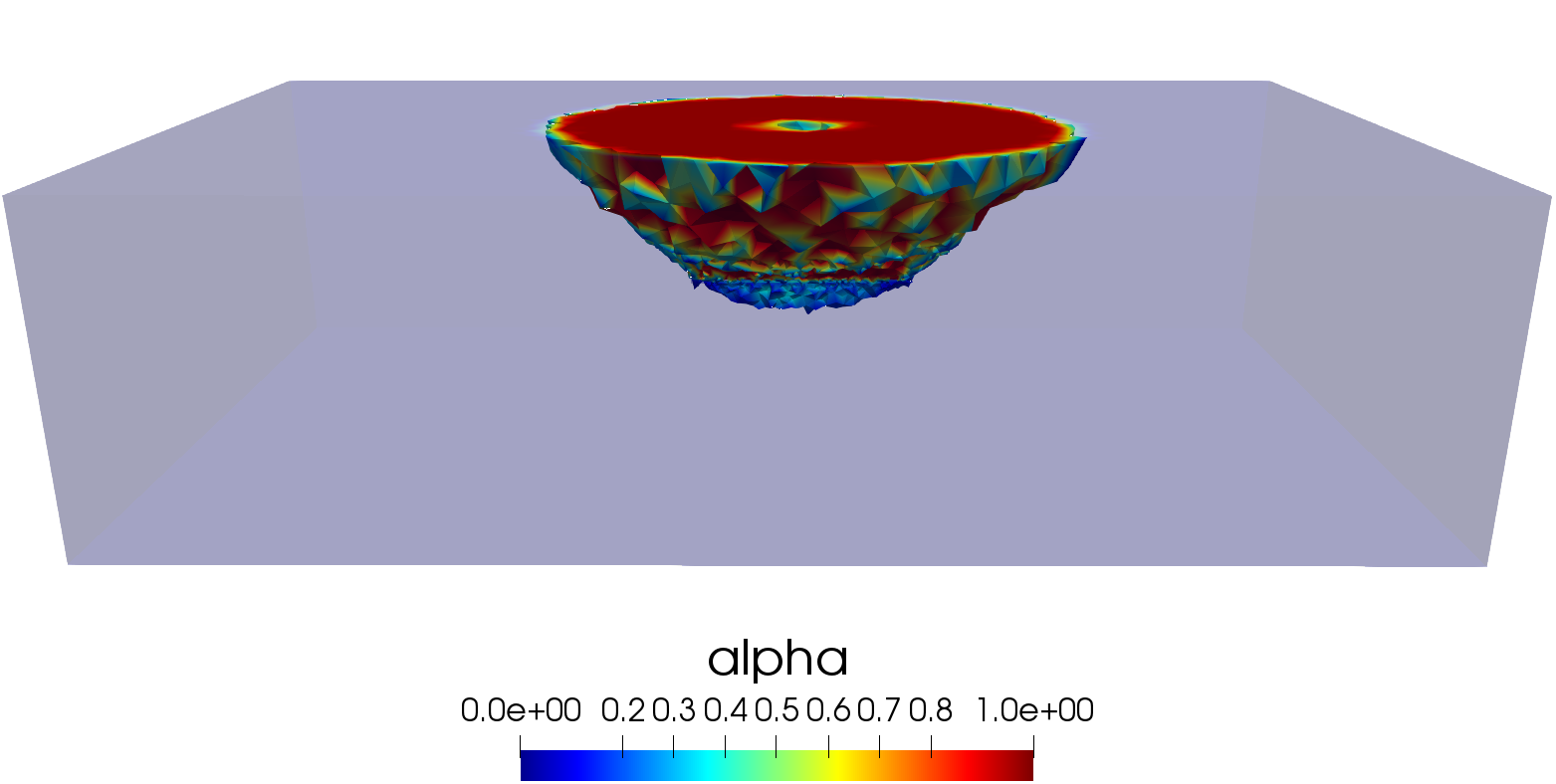}
		\caption{Model 1.}
	\end{subfigure}
	~
	\begin{subfigure}[b]{0.49\textwidth}
		\centering
		\includegraphics[width=\textwidth]{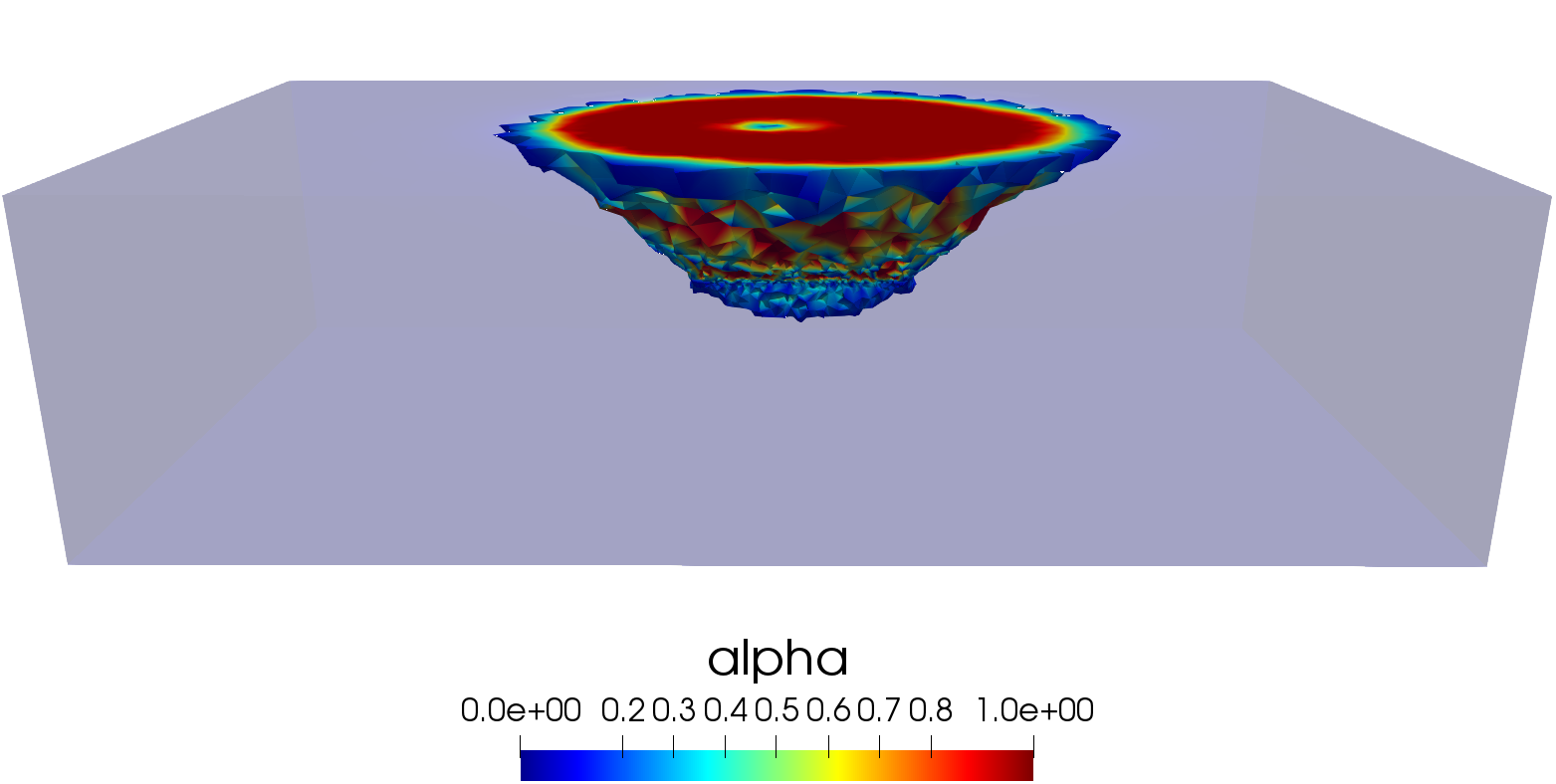}
		\caption{Model 2.}
	\end{subfigure}
	\qquad
	\begin{subfigure}[b]{0.49\textwidth}
		\centering
		\includegraphics[width=\textwidth]{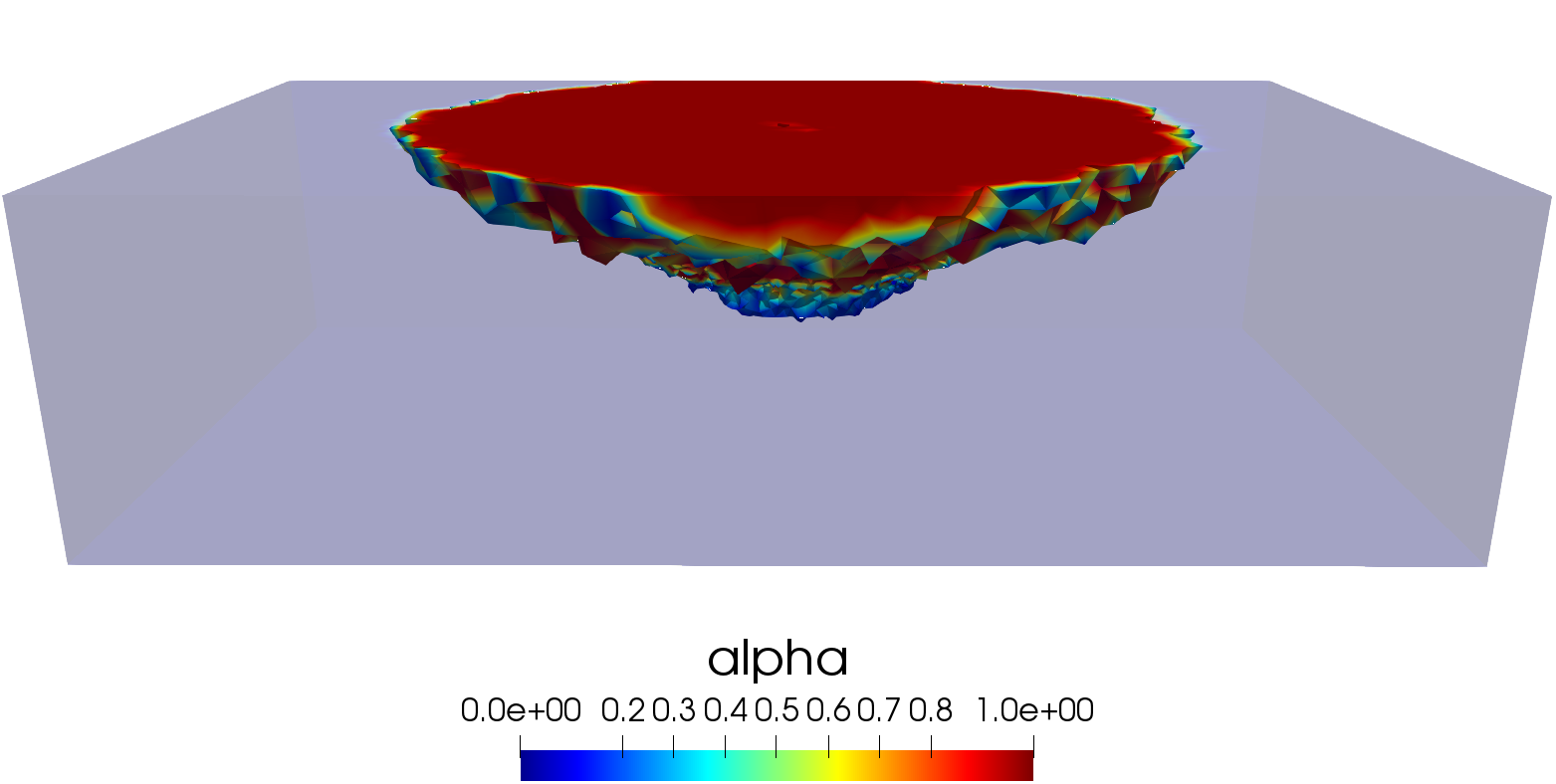}
		\caption{Model 3.}
	\end{subfigure}
	~
	\begin{subfigure}[b]{0.49\textwidth}
		\centering
		\includegraphics[width=\textwidth]{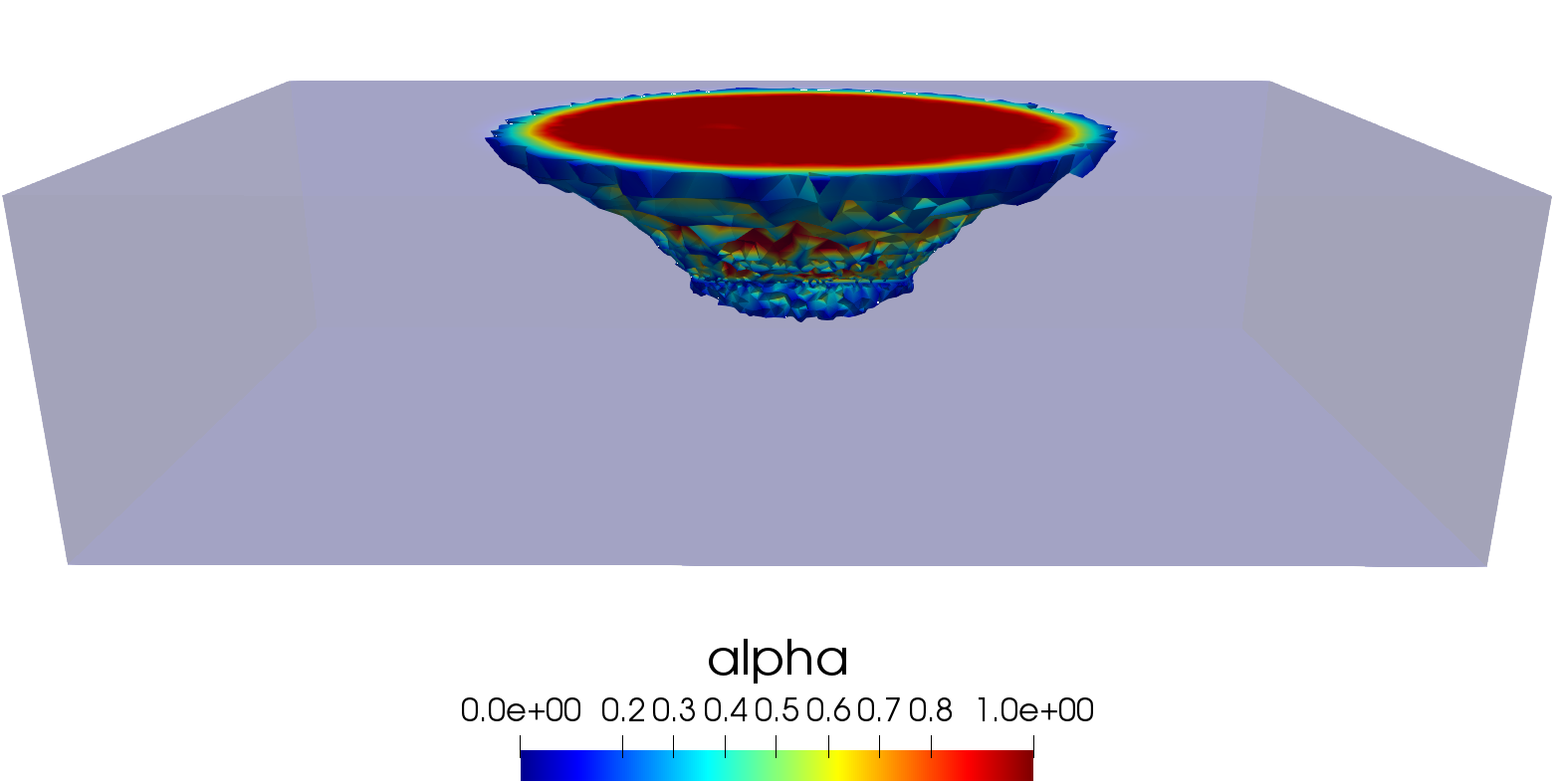}
		\caption{Model 4.}
	\end{subfigure}
	\caption{Damage field distribution in the rock mass for different models in 3D.}\label{Mod1-4_3d}
\end{figure}

For 3D visualization, the figure \ref{Mod1-4_3d} displays the distribution of the damage in the final state for each model. In this view is possible see that the damage is close to zero almost everywhere except around the cavity and in all models the damage reaches the upper boundary with different magnitudes and distributions.

\section{Fast algorithm to solve the block caving process}\label{NewAlg}

\subsection{Errors and new algorithm}

The results, obtained by Agorithm \ref{alg:damage} in the previous section, suffer oscillations in the errors to achieve the convergence. This causes felays in the solving tho problem, causing more time to find the solution to this problem.

\begin{figure}[h!]
	\centering
	\includegraphics[width=\textwidth, height=\textwidth]{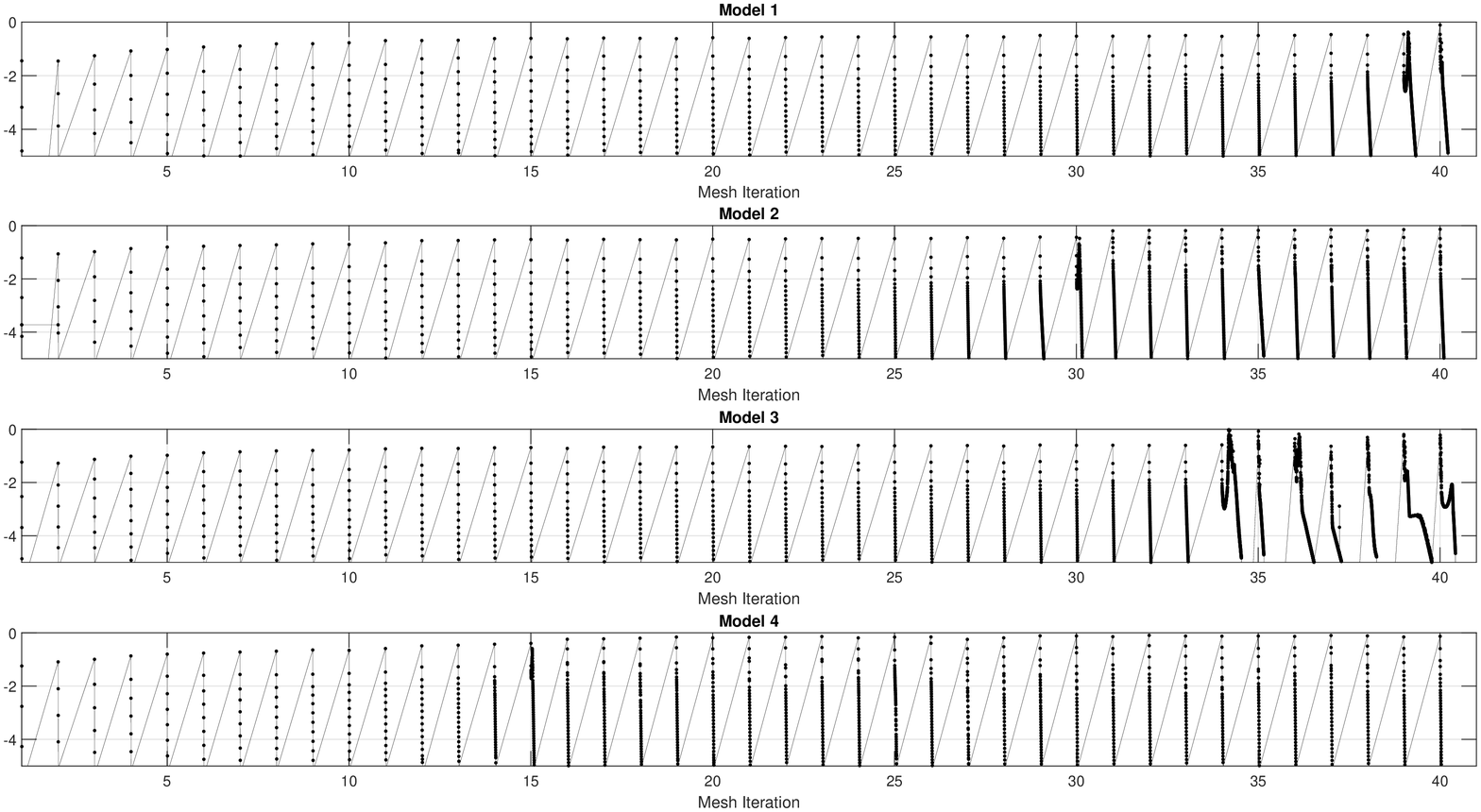}
	\caption{Logarithm of the error  for each cavity advance.}
	\label{errormodels}
\end{figure}

Figure \ref{errormodels} shows the logarithm of the errors in each mesh iteration fot the four models. It is possible to see that, For the Model 1,  the oscillations are obtained in the last time step. In the Model 2 and 3, the oscillations appear around the mesh step 30 and 31 respectively. Finally, in the Model 4, the are no oscillations.

\begin{figure}[h!]
	\centering
	\begin{subfigure}[b]{0.49\textwidth}
		\centering
		\includegraphics[width=\textwidth]{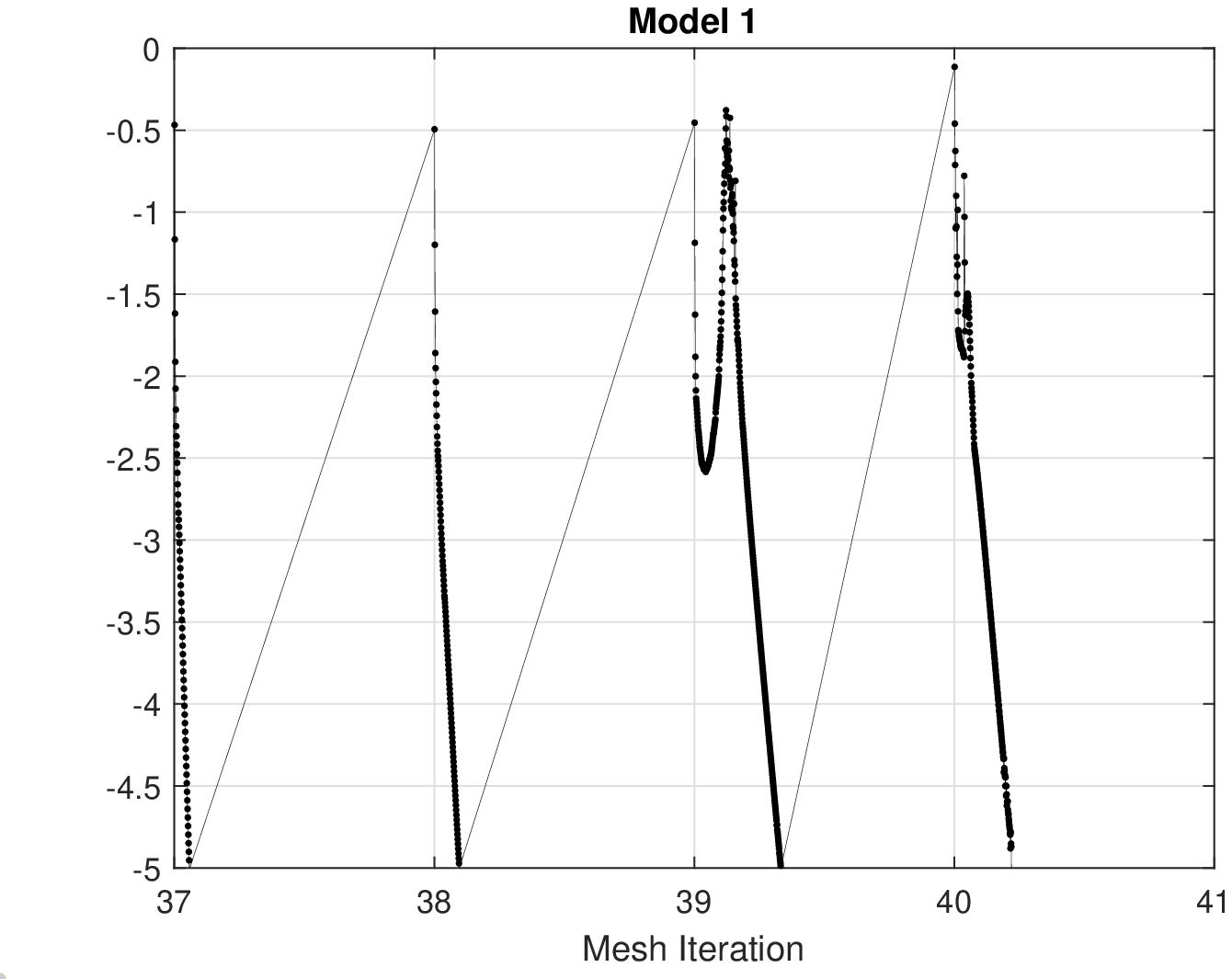}
	\end{subfigure}
	\begin{subfigure}[b]{0.49\textwidth}
		\centering
		\includegraphics[width=\textwidth]{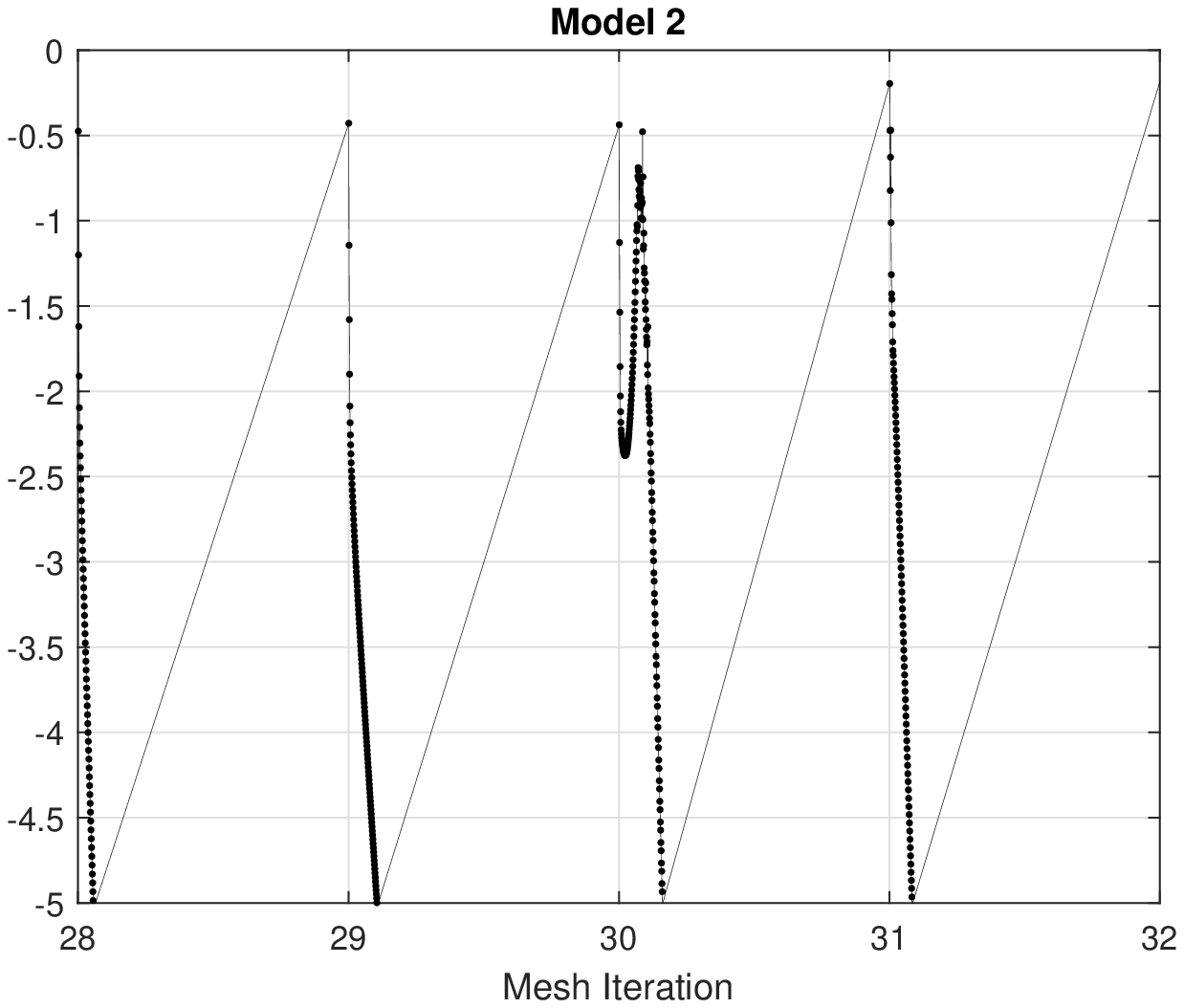}
	\end{subfigure}
	\begin{subfigure}[b]{0.49\textwidth}
		\centering
		\includegraphics[width=\textwidth]{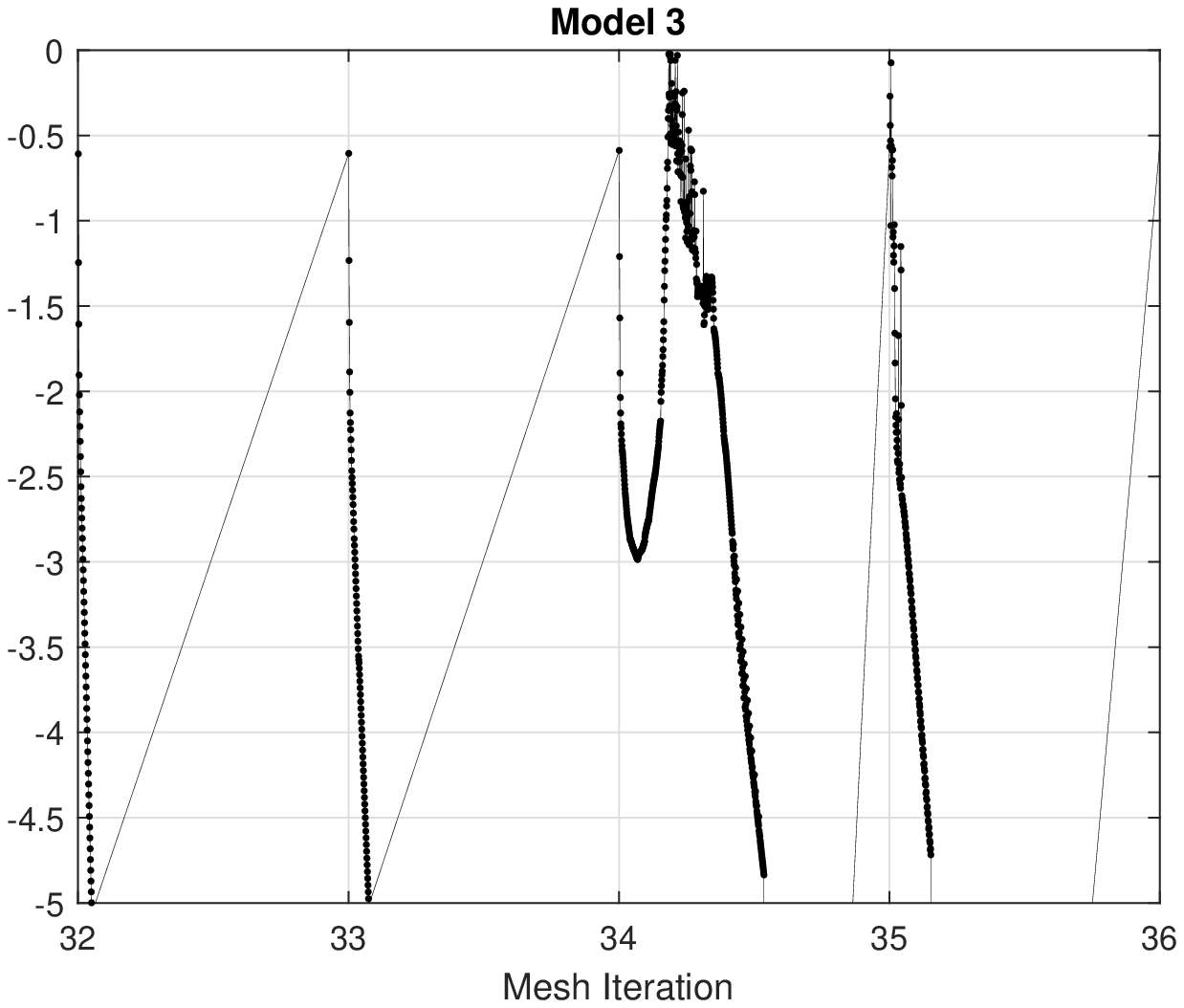}
	\end{subfigure}
	\begin{subfigure}[b]{0.49\textwidth}
		\centering
		\includegraphics[width=\textwidth]{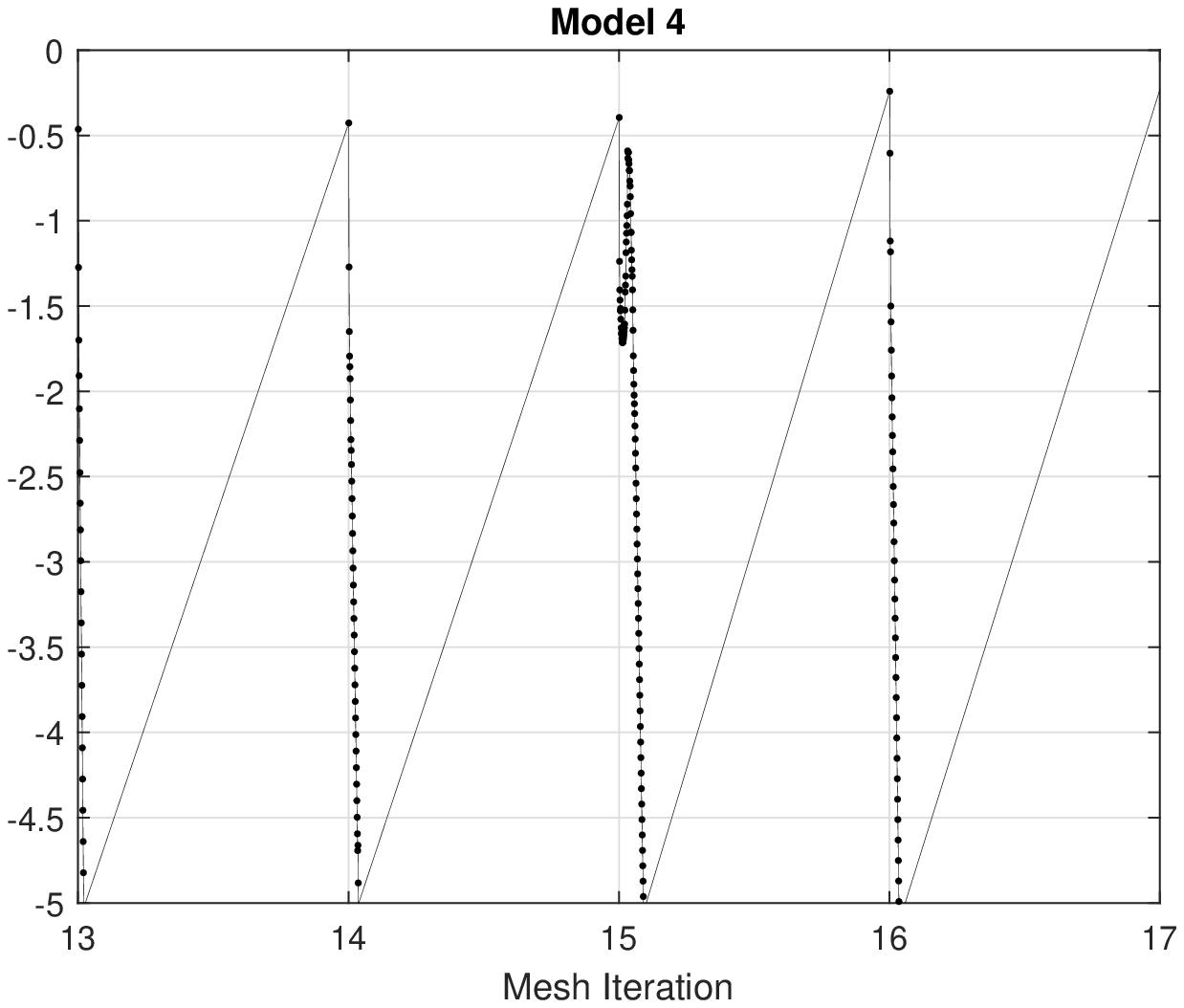}
	\end{subfigure}
	\caption{Zoom of the logarithm of the errors.}
	\label{errormodels_zoom}
\end{figure}

Figure \ref{errormodels_zoom} shows a zoom of the oscillations in the mesh step of these oscillations occur. It can be seen that the oscillations appear in the same mesh step where attained its maximum value (see Figure \ref{alphamaxmods}).

To avoid the oscillations in the solution of the problem and with this improve the time in which the solution is obtained, we propose a modification in the Algorithm \ref{alg:damage}. This modification consists in considering, when the error grows, and arrangement of the form
\begin{equation}
\alpha^{(p)}= C_L \alpha^{(p-1)}+(1-C_L) \alpha^{(p)},
\end{equation}
where this arrangement is a combination between the solution in iteration $(p)$ and the previous solution in the iteration $(p-1)$, with $C_L \in (0,1)$ a constant. This new approach to solve the damage problem is summarized in Algorithm \ref{alg:damage2}.
\begin{algorithm}[h!]
	\begin{algorithmic}[1]
		\RETURN Solution at time step $t_i$.
		\STATE Given $(u_{i-1},\alpha_{i-1})$, the sate at the previous loading step.
		\STATE Set $(u^{(0)},\alpha^{(0)}):=(u_{i-1},\alpha_{i-1})$ and error$^{(0)} = 1.0$ 
		\WHILE {error$^{(p)}>$tolerance}
		\STATE Solve $u^{(p)}$ from (\ref{elasteq}) with $\alpha^{(p-1)}$.
		\STATE Find $\displaystyle \alpha^{(p)}:= \argmin_{\alpha \in \mathcal{D}(\alpha_{i-1})} \overline{\mathcal{P}}(u^{(p)},\alpha)$.
		\STATE error$^{(p)} = \Vert\alpha^{(p-1)} -\alpha^{(p)}\Vert_{\infty}$.
		\WHILE{error$^{(p)}>$error$^{(p-1)}$}
		\STATE $\overline{\alpha}^{(p)}=C_L \alpha^{(p-1)}+(1-C_L)\alpha^{(p)}$.
		\STATE $\alpha^{(p)}=\overline{\alpha}^{(p)}$.
		\STATE error$^{(p)} = \Vert\alpha^{(p-1)} -\alpha^{(p)}\Vert_{\infty}$.
		\ENDWHILE
		\ENDWHILE
		\STATE  Set $(\mathbf{u}_{i}, \alpha_{i}) = (u^{p}, \alpha^{p}).$
	\end{algorithmic}
	\caption{Fast algorithm to solve the damage problem}\label{alg:damage2}
\end{algorithm}

\subsection{numerical results}

\begin{figure}[h!]
	\centering
	\includegraphics[width=\textwidth]{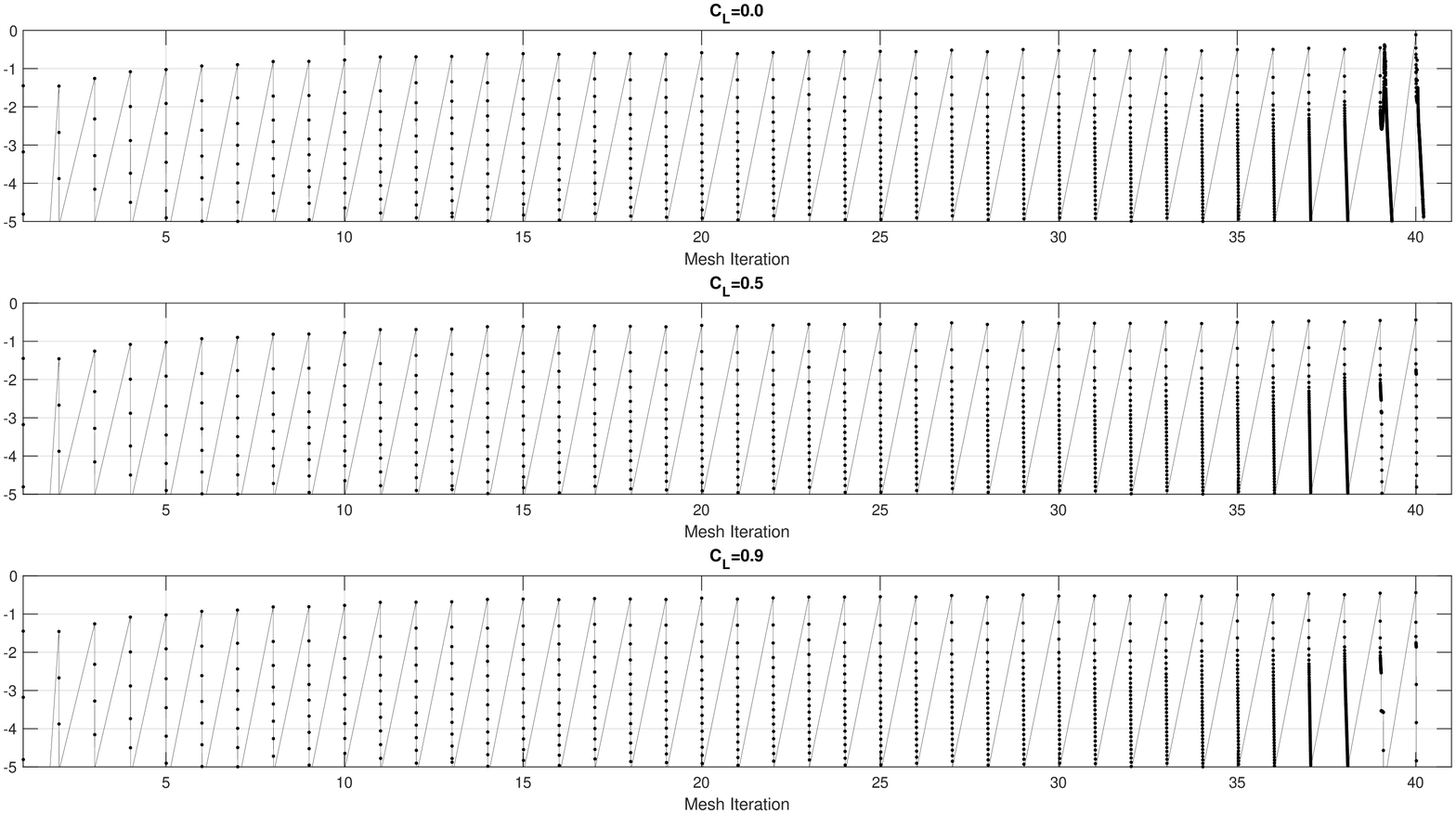} 
	\caption{Logarithm of the error in the Model 1 for different values of $C_L$.}
	\label{errorsimpmod1}
\end{figure}

\begin{figure}[h!]
	\centering
	\includegraphics[width=\textwidth]{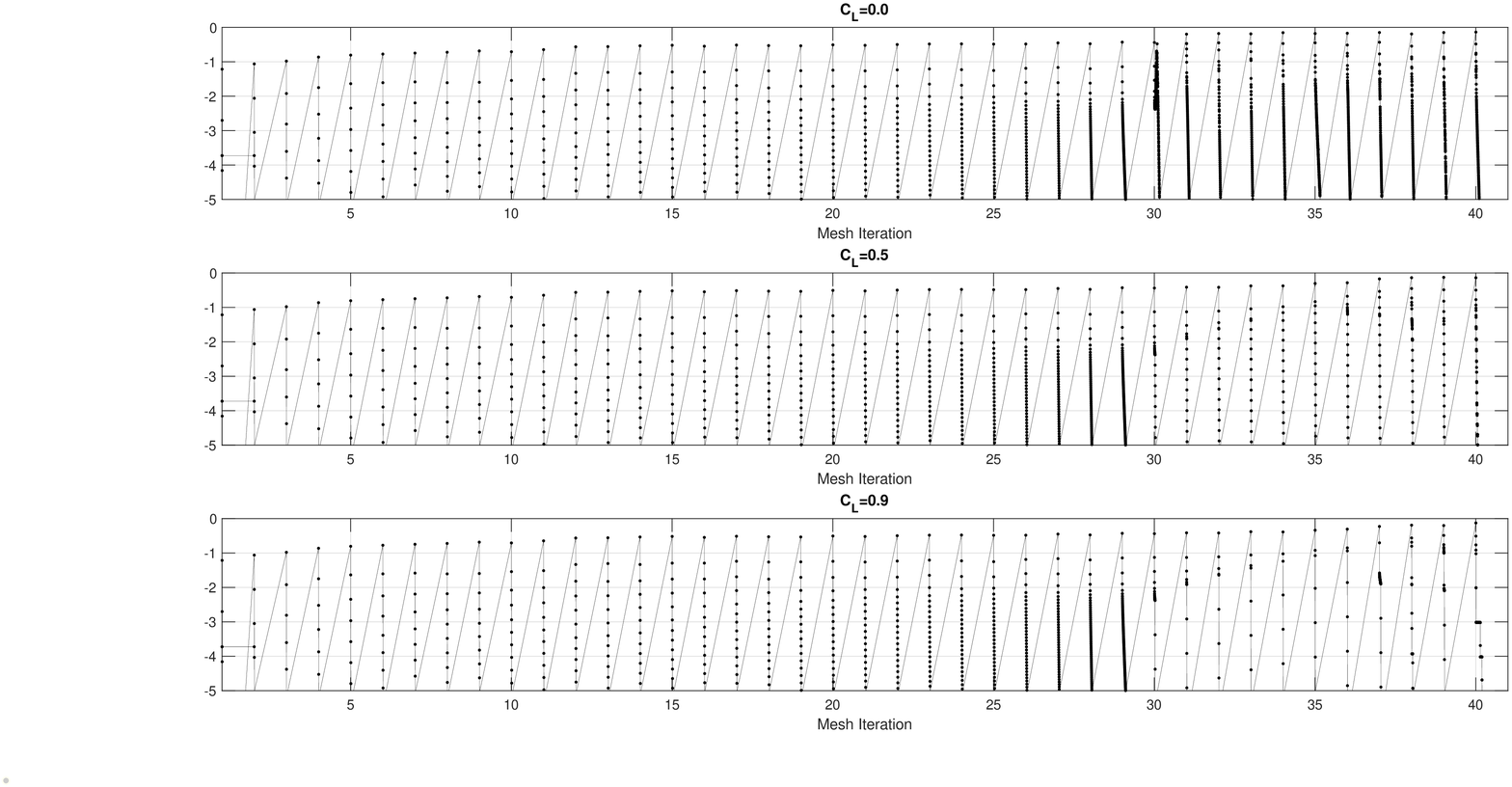} 
	\caption{Logarithm of the error in the Model 2 for different values of $C_L$.}
	\label{errorsimpmod2}
\end{figure}

\begin{figure}[h!]
	\centering
	\includegraphics[width=\textwidth]{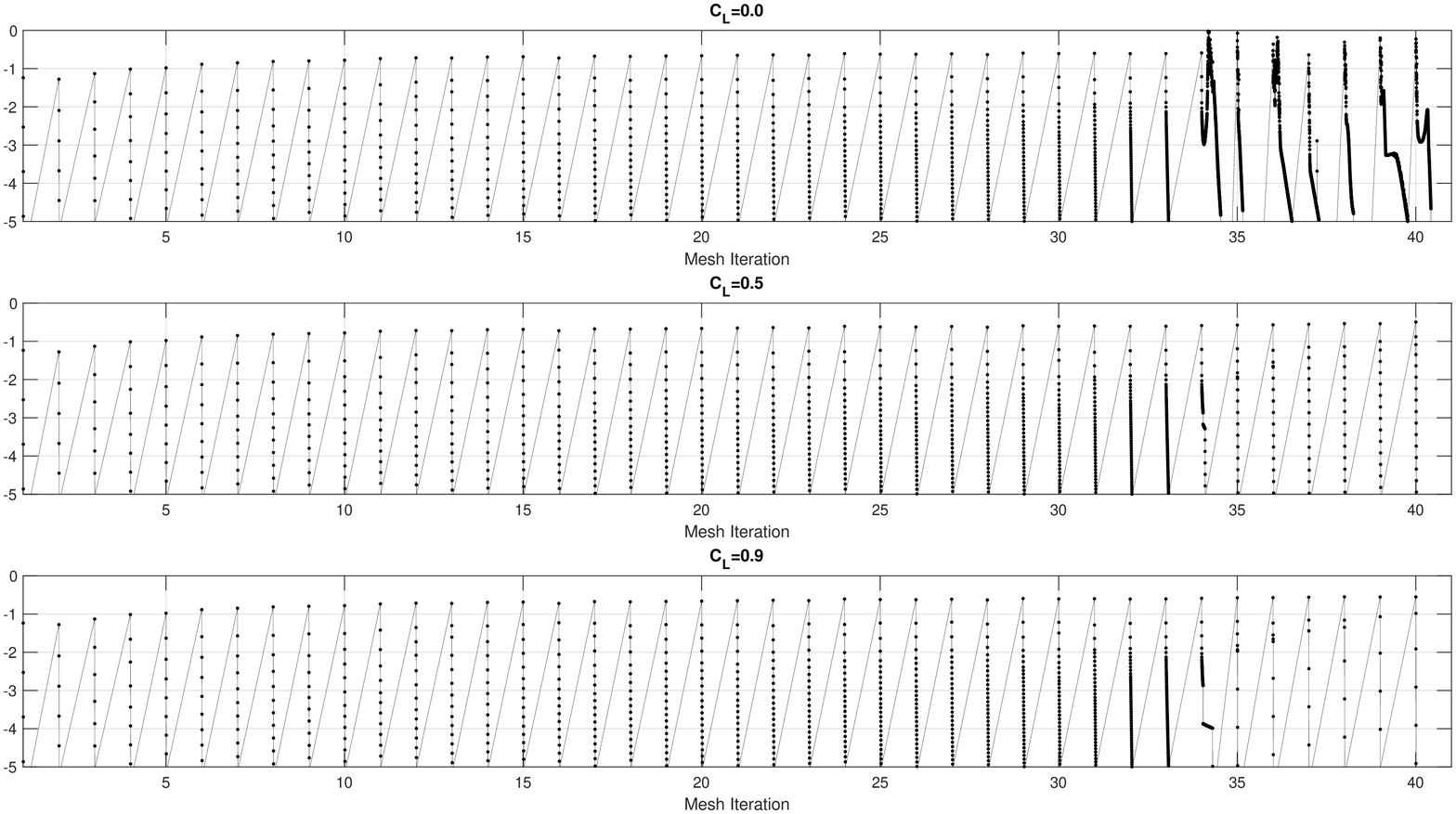} 
	\caption{Logarithm of the error in the Model 3 for different values of $C_L$.}
	\label{errorsimpmod3}
\end{figure}

\begin{figure}[h!]
	\centering
	\includegraphics[width=\textwidth]{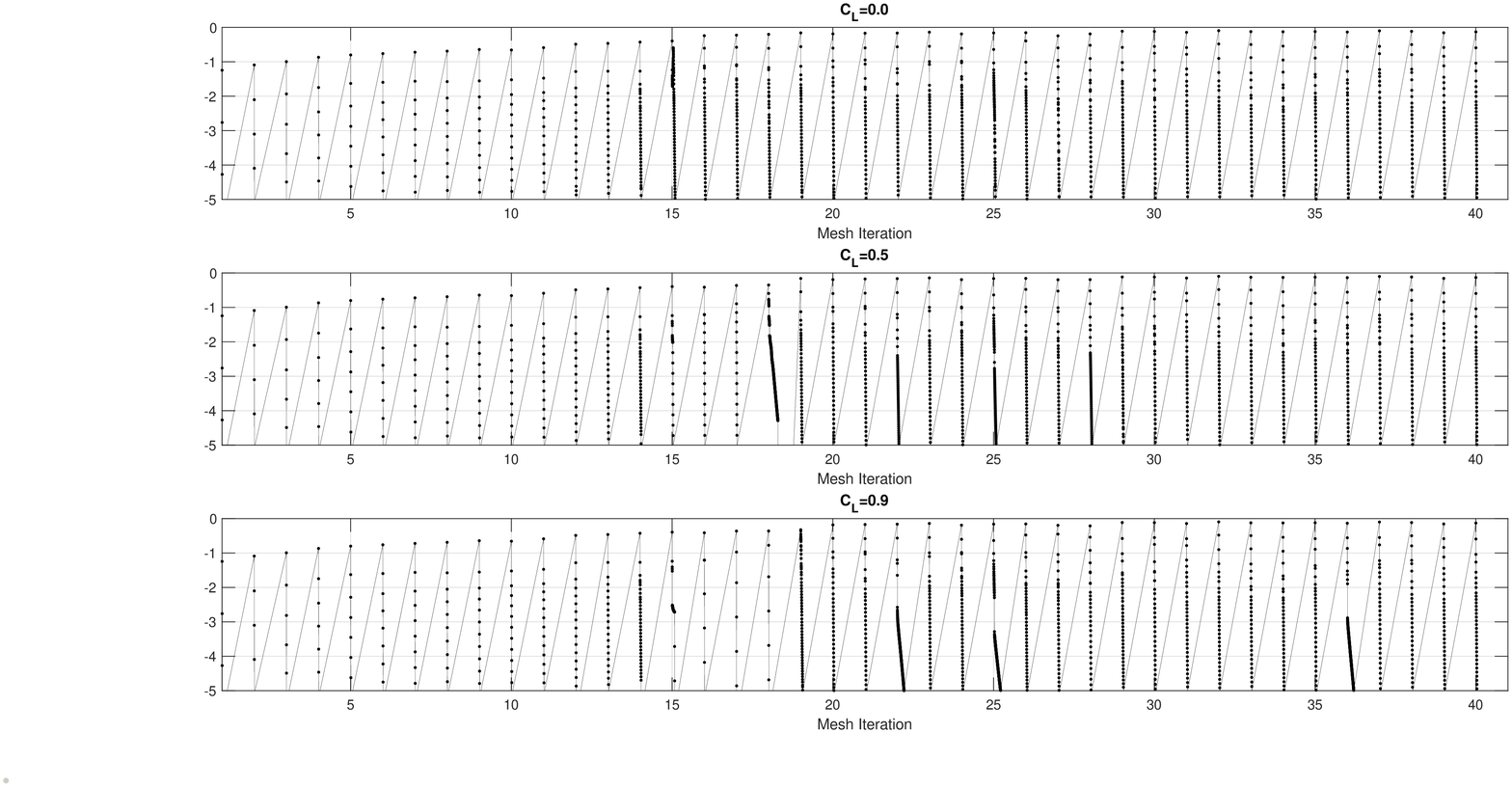} 
	\caption{Logarithm of the error in the Model 4 for different values of $C_L$.}
	\label{errorsimpmod4}
\end{figure}

For this new algorithm, we consider the same previous test, where we apply the Algorithm \ref{alg:damage2} to all models. The Figures \ref{errorsimpmod1}-\ref{errorsimpmod4} display the logarithm of the error for $C_L=0.5$ and $C_L= 0.9$ respectively and for the different models. It is possible to see that the oscillations of the error disappear that with a smaller amount of iterations the convergence of the algorithm is attained, causing with this that the solution is obtained in a lower time.

\begin{figure}[h!]
	\centering
	\begin{subfigure}[b]{0.49\textwidth}
	    \centering
		\includegraphics[width=\textwidth]{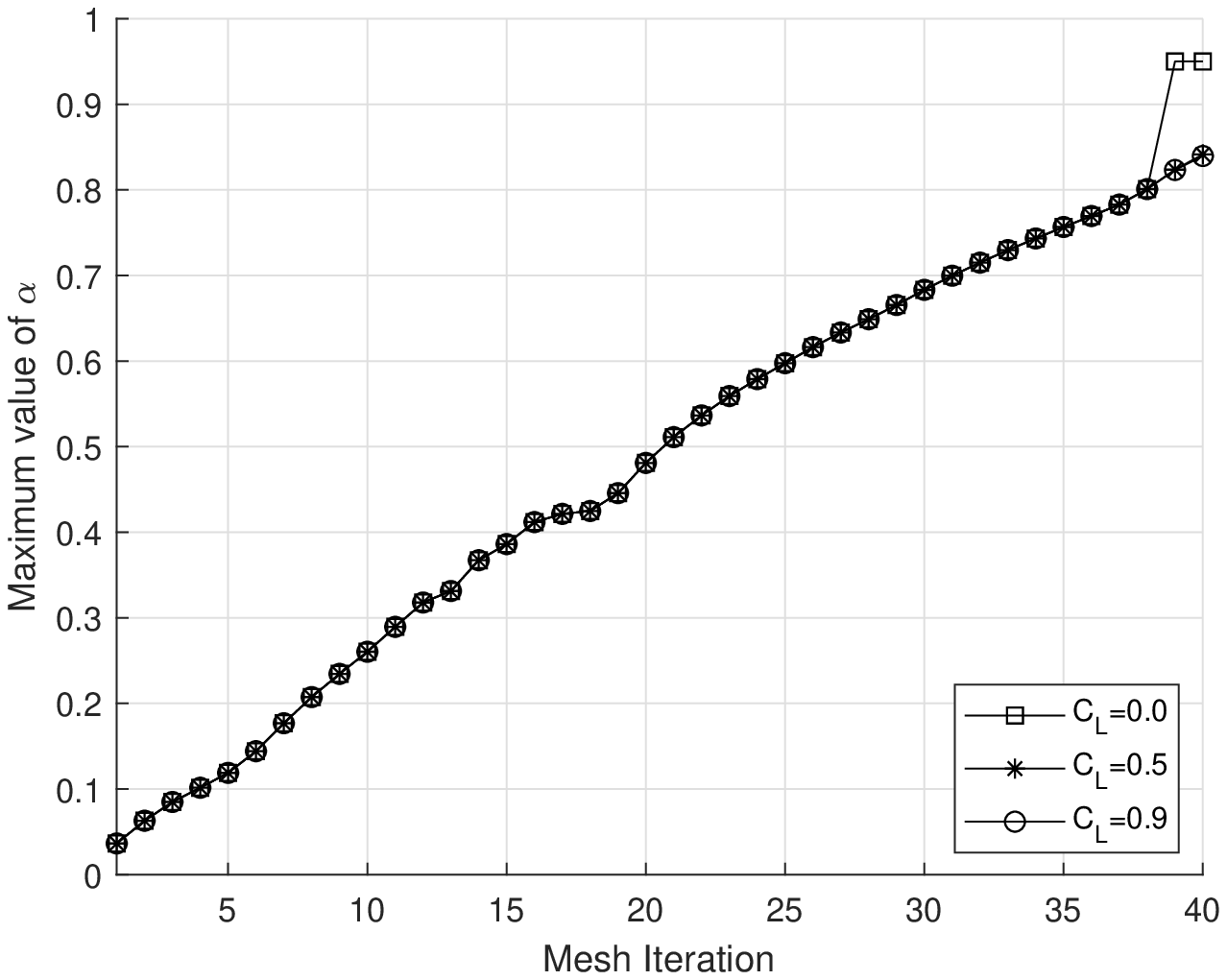}
		\caption{Model 1.}
	\end{subfigure}
	    \centering
		\begin{subfigure}[b]{0.49\textwidth}
		\includegraphics[width=\textwidth]{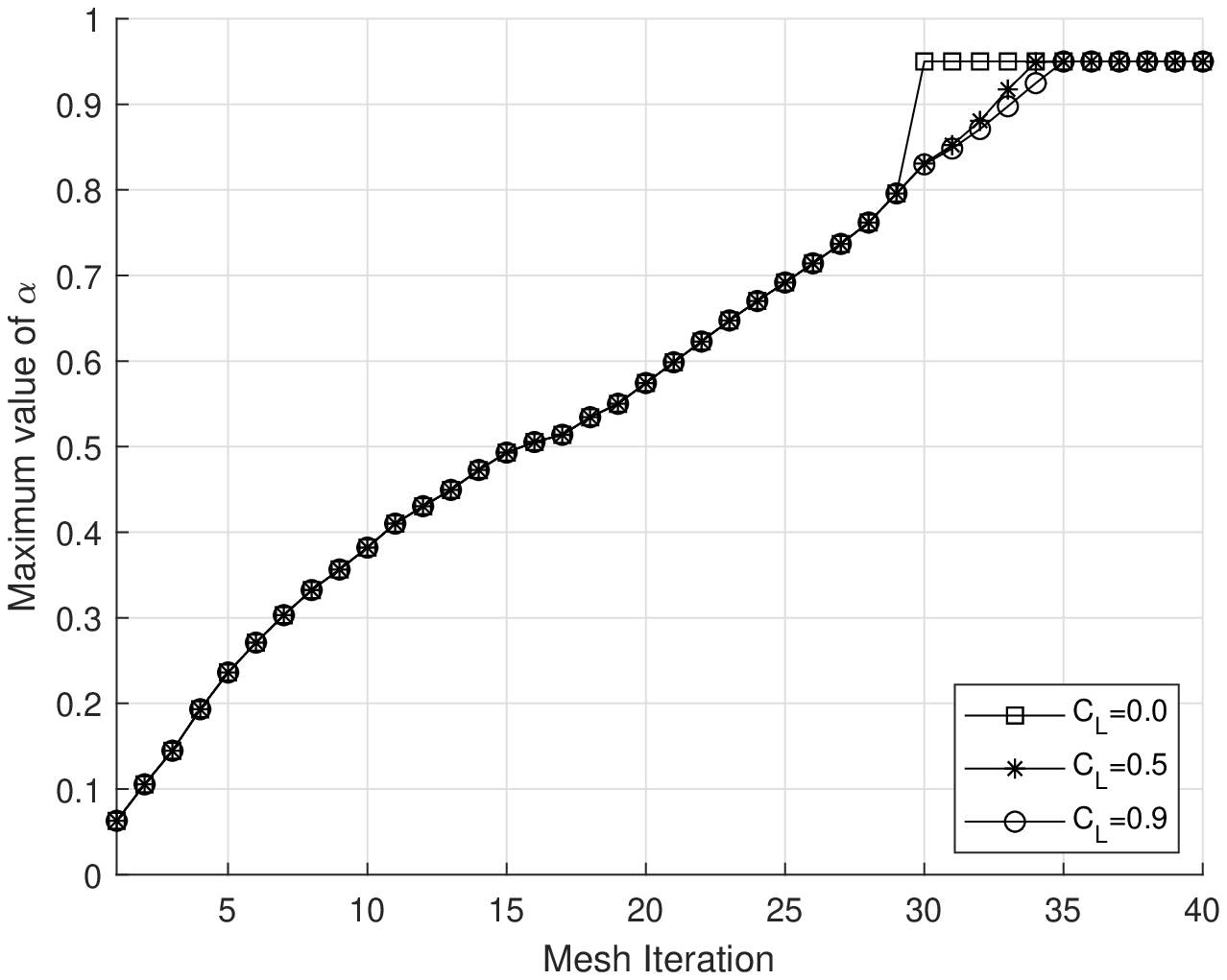}
		\caption{Model 2.}
	\end{subfigure}
	    \centering
	\begin{subfigure}[b]{0.49\textwidth}
		\centering
		\includegraphics[width=\textwidth]{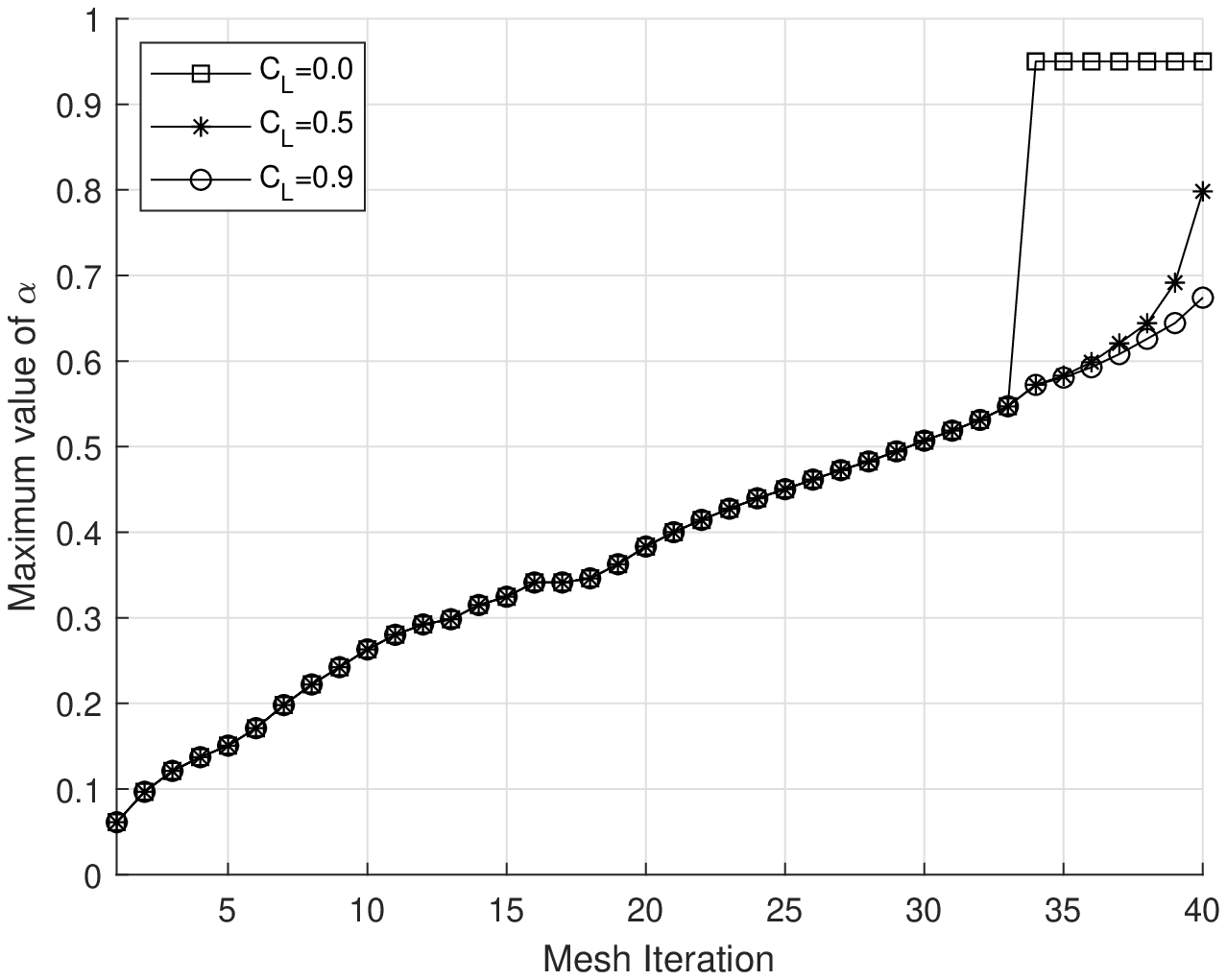}
		\caption{Model 3.}
	\end{subfigure}
	\begin{subfigure}[b]{0.49\textwidth}
		\centering
		\includegraphics[width=\textwidth]{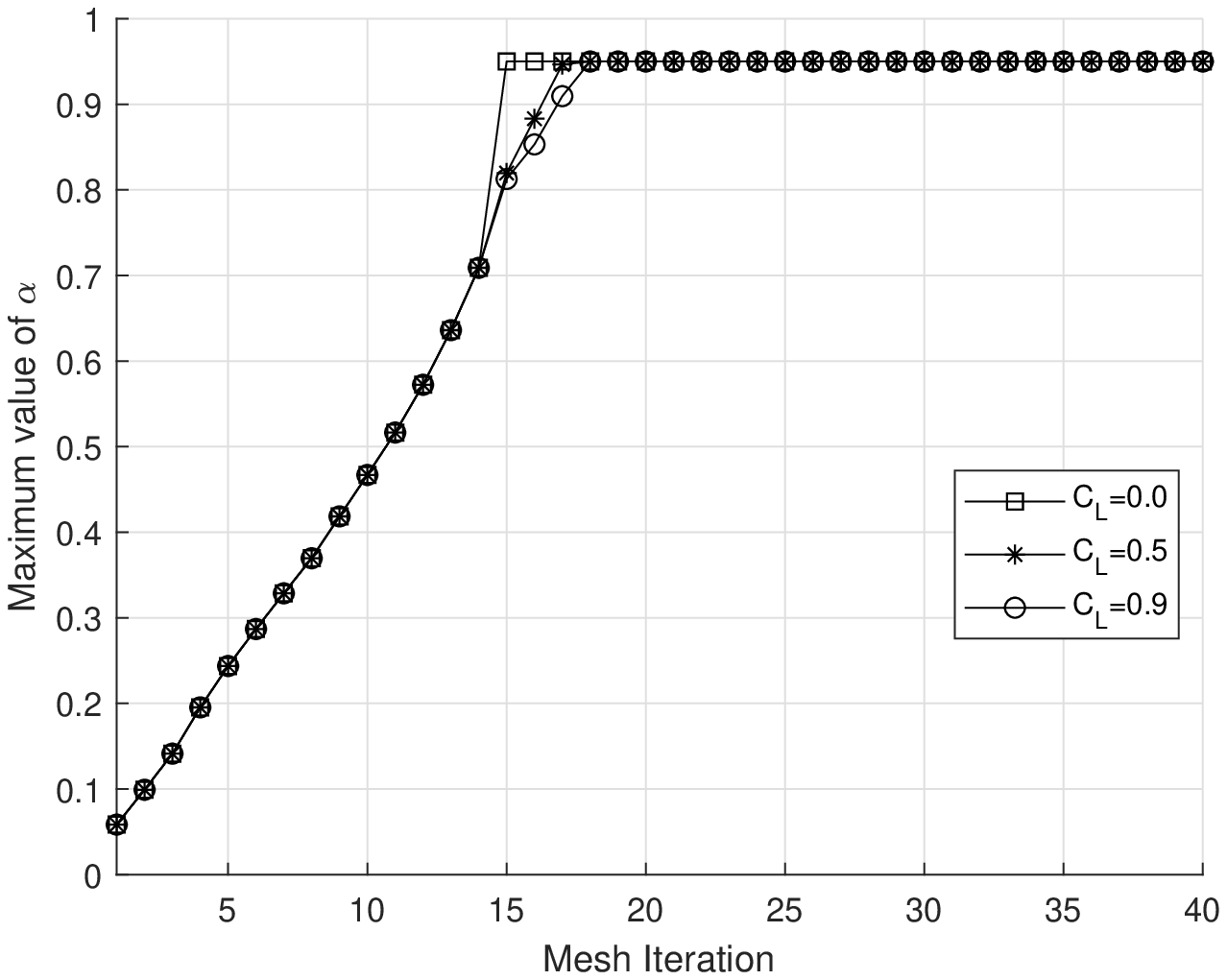}
		\caption{Model 4.}
	\end{subfigure}
	\caption{Maximum damage value for each mesh step for different values of $C_L$.}\label{alphamaxmodsimp}
\end{figure}

Figure \ref{alphamaxmodsimp} summarized the maximum values of $\alpha$ for the different values of $C_L$ and different models in each mesh iteration. The figure displays the evolution of this maximum value and it can be see that the value of the constant $C_L$ is higher. It is possible see that the value of the maximum damage is the same for the different values of $C_L$ in the final mesh step, except in the Model 1 that the maximum value is lower, for this reason we will consider $C_L=0.9$ since this value require less iterations to achieve the convergence. 

\begin{figure}[h!]
	\centering
	\begin{subfigure}[b]{0.49\textwidth}
		\centering
		\includegraphics[width=\textwidth]{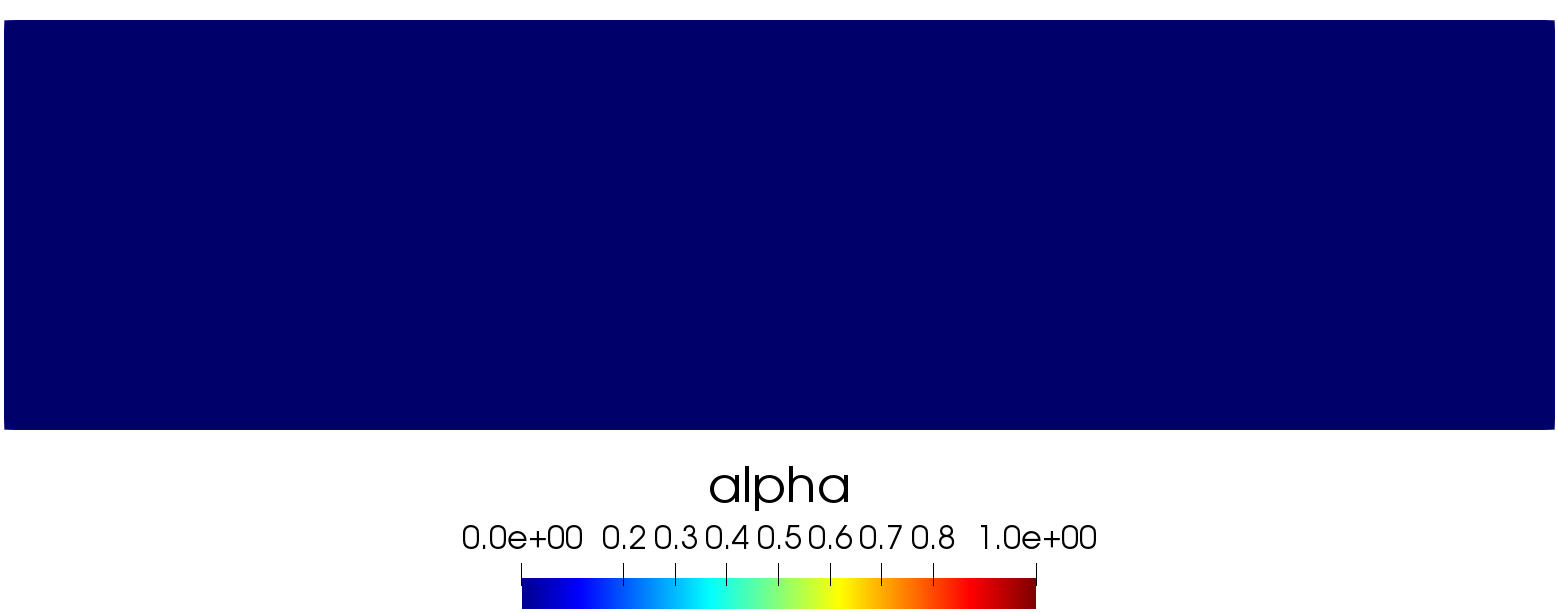}
		\caption{$\Omega(t_{0})$}
	\end{subfigure}
	\begin{subfigure}[b]{0.49\textwidth}
		\centering
		\includegraphics[width=\textwidth]{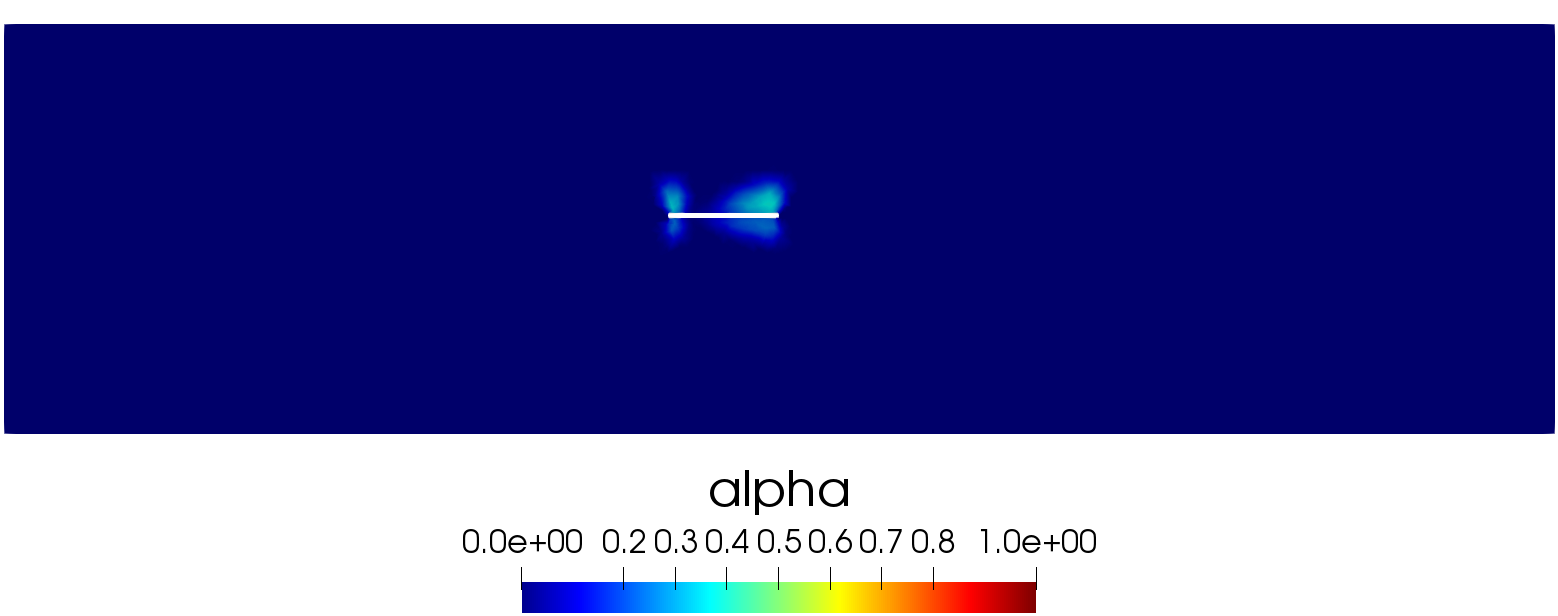}
		\caption{$\Omega(t_{15})$}
	\end{subfigure}
	\begin{subfigure}[b]{0.49\textwidth}
		\centering
		\includegraphics[width=\textwidth]{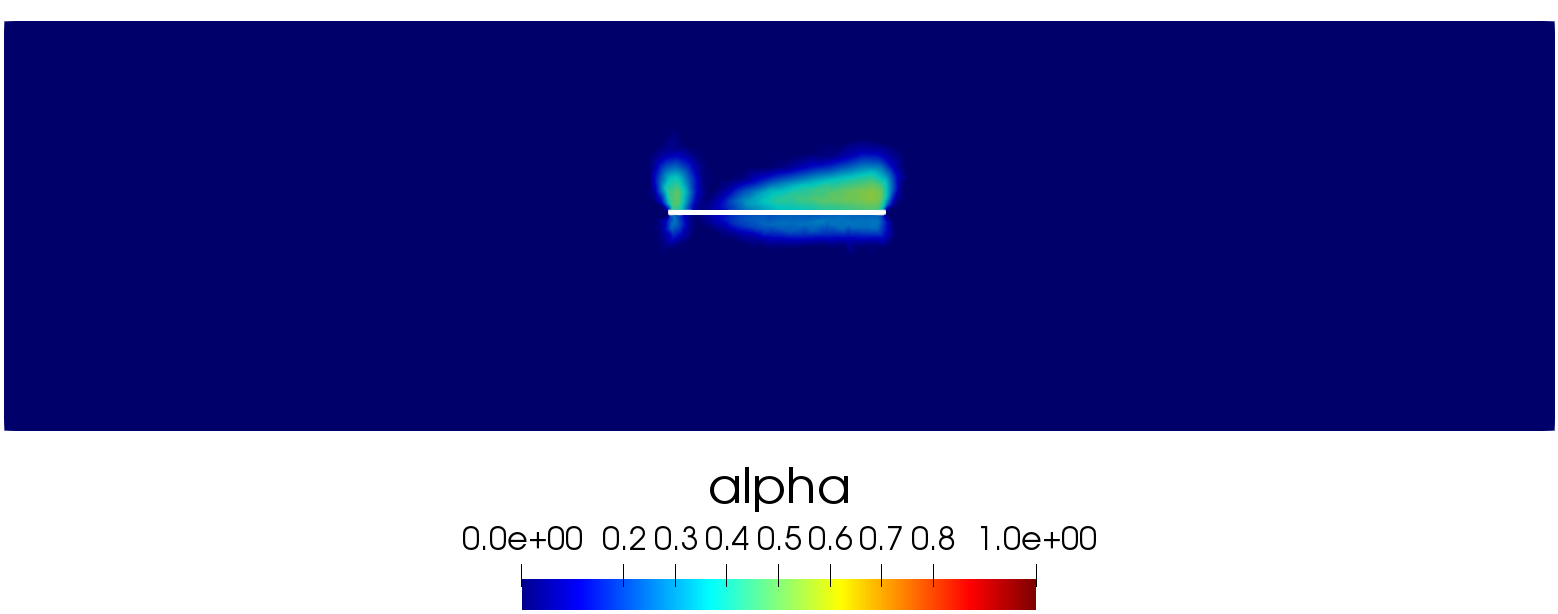}
		\caption{$\Omega(t_{30})$}
	\end{subfigure}	
	\begin{subfigure}[b]{0.49\textwidth}
		\centering
		\includegraphics[width=\textwidth]{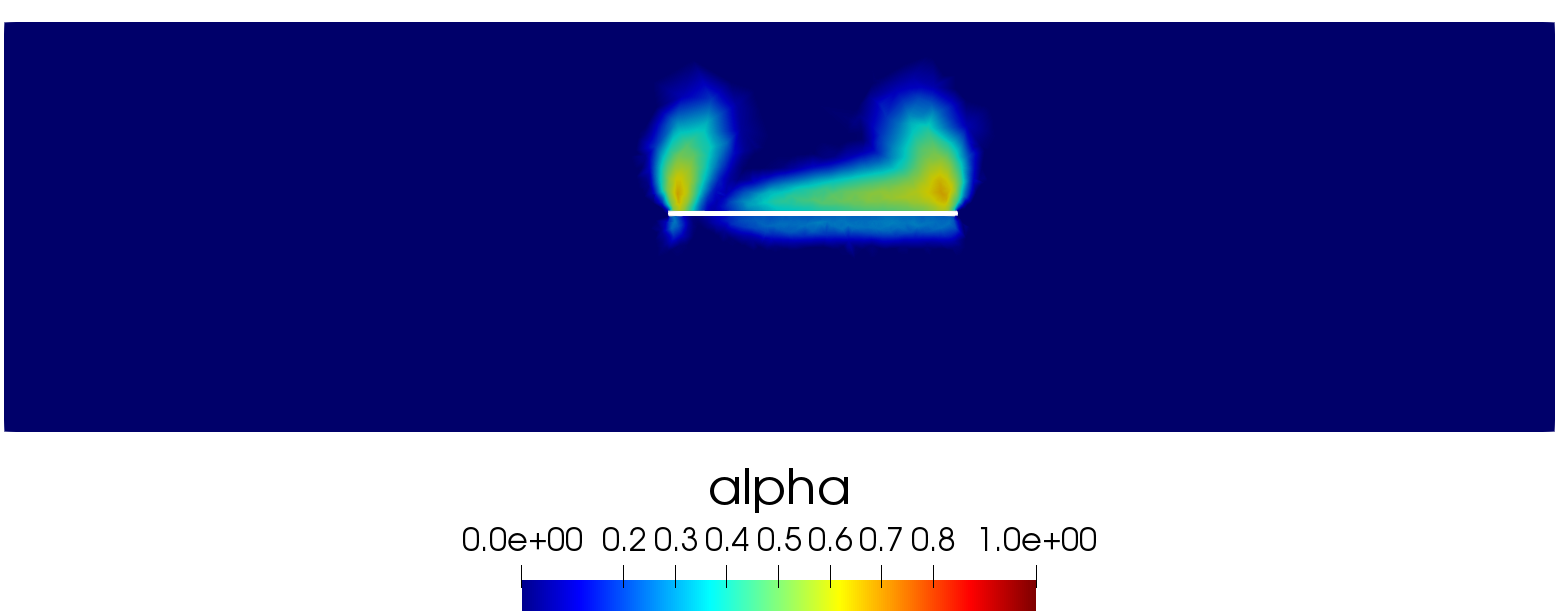}
		\caption{$\Omega(t_{40})$}
	\end{subfigure}
	\caption{Damage field distribution in the rock mass for Model 1 and $C_L=0.9$.}\label{mod1latcutCL09}
\end{figure}

\begin{figure}[h!]
	\centering
	\begin{subfigure}[b]{0.49\textwidth}
		\centering
		\includegraphics[width=\textwidth]{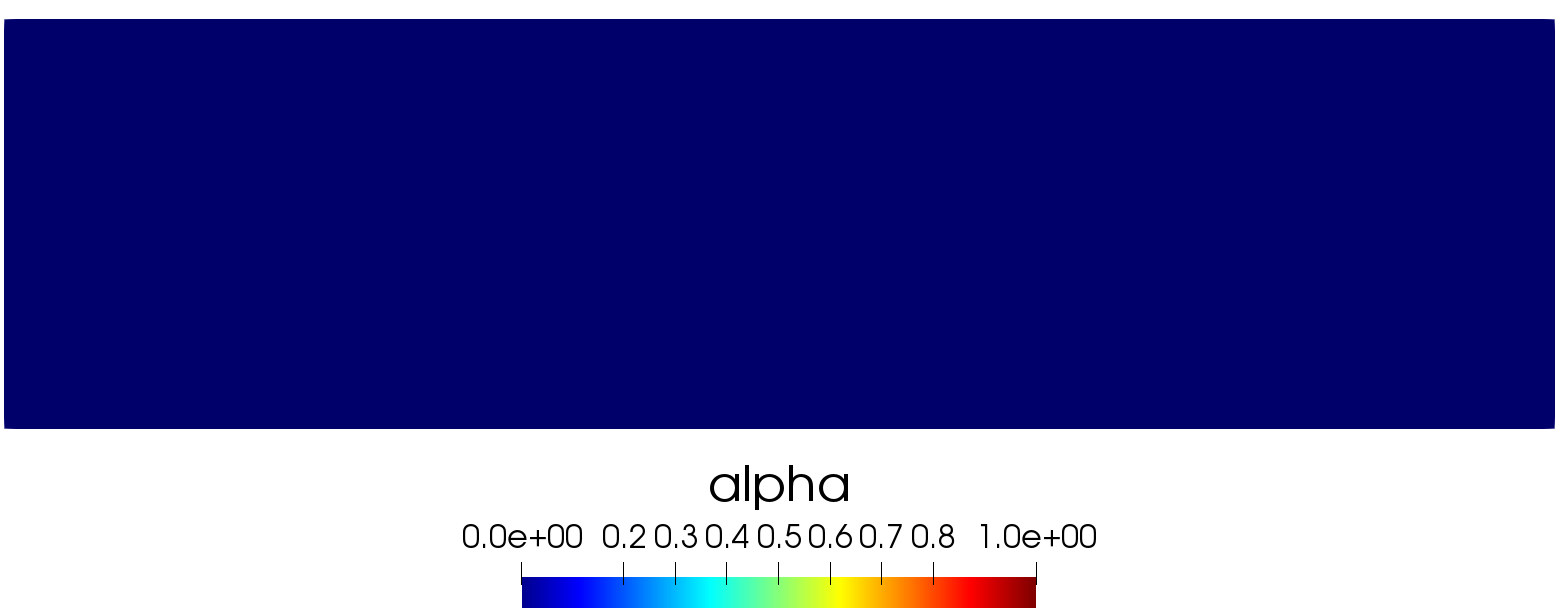}
		\caption{$\Omega(t_{0})$}
	\end{subfigure}
	\begin{subfigure}[b]{0.49\textwidth}
		\centering
		\includegraphics[width=\textwidth]{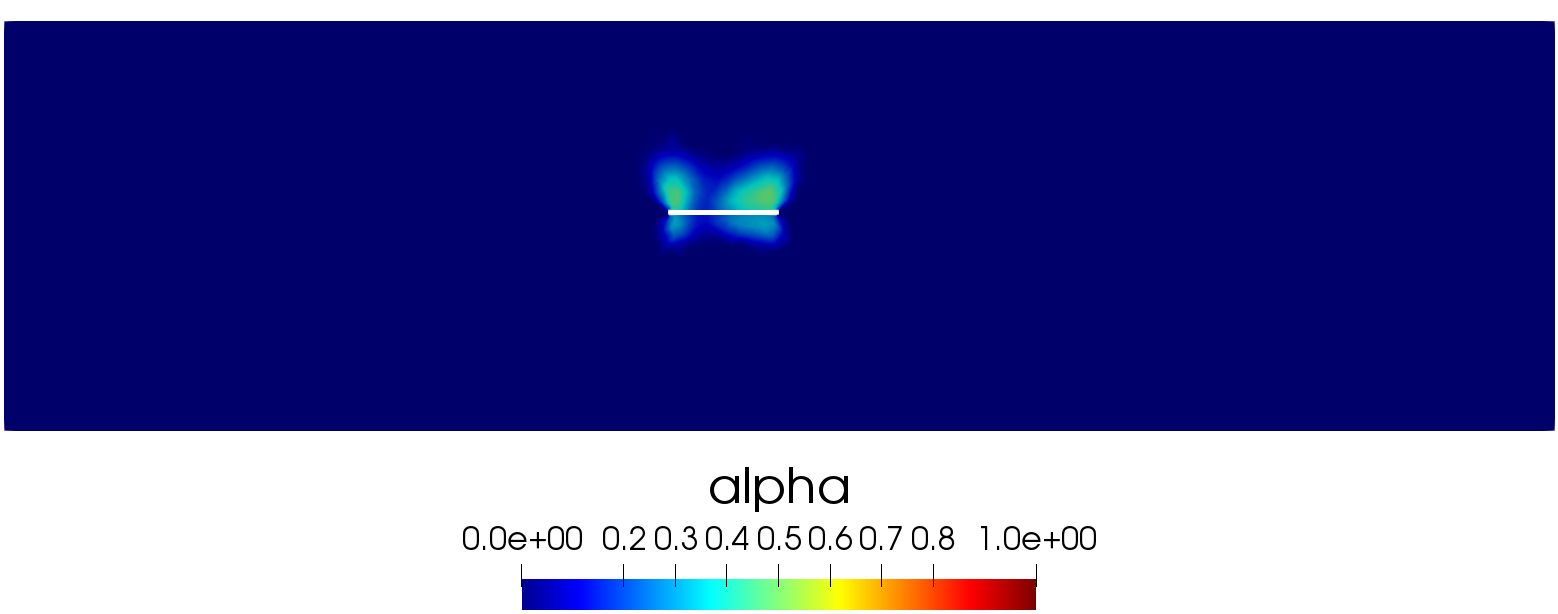}
		\caption{$\Omega(t_{15})$}
	\end{subfigure}
	\begin{subfigure}[b]{0.49\textwidth}
		\centering
		\includegraphics[width=\textwidth]{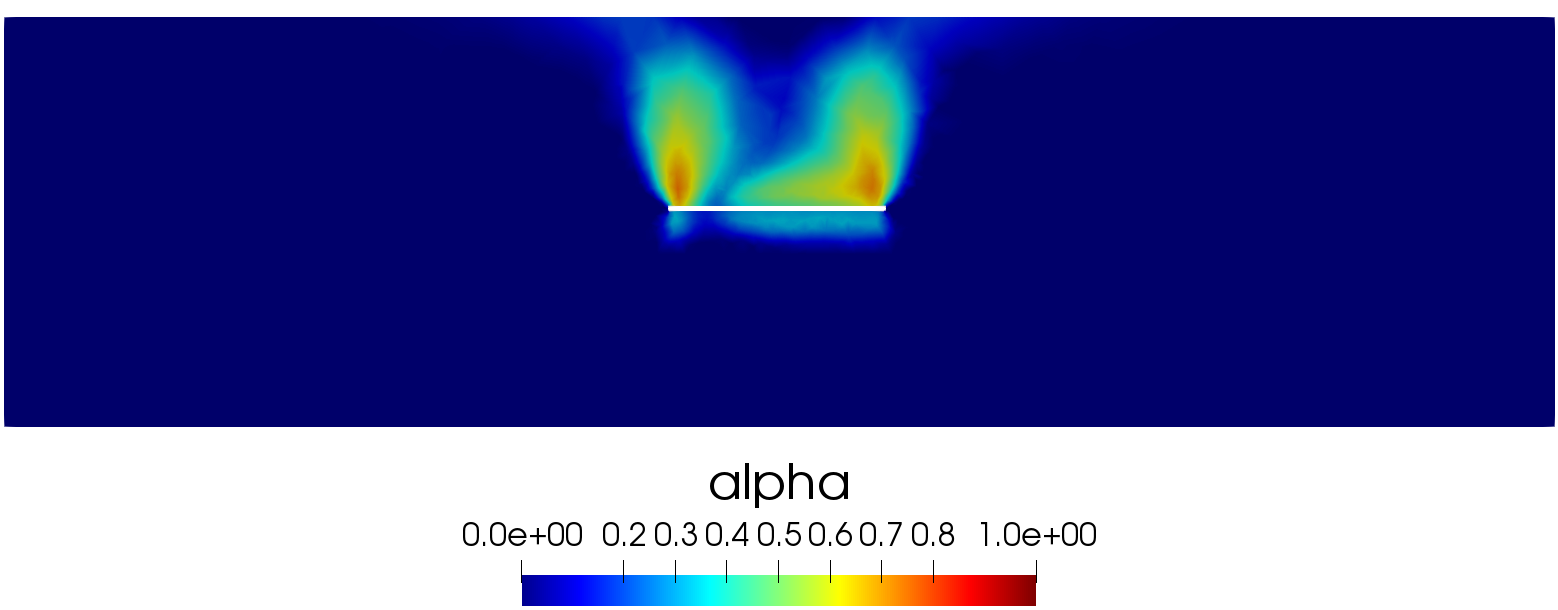}
		\caption{$\Omega(t_{30})$}
	\end{subfigure}	
	\begin{subfigure}[b]{0.49\textwidth}
		\centering
		\includegraphics[width=\textwidth]{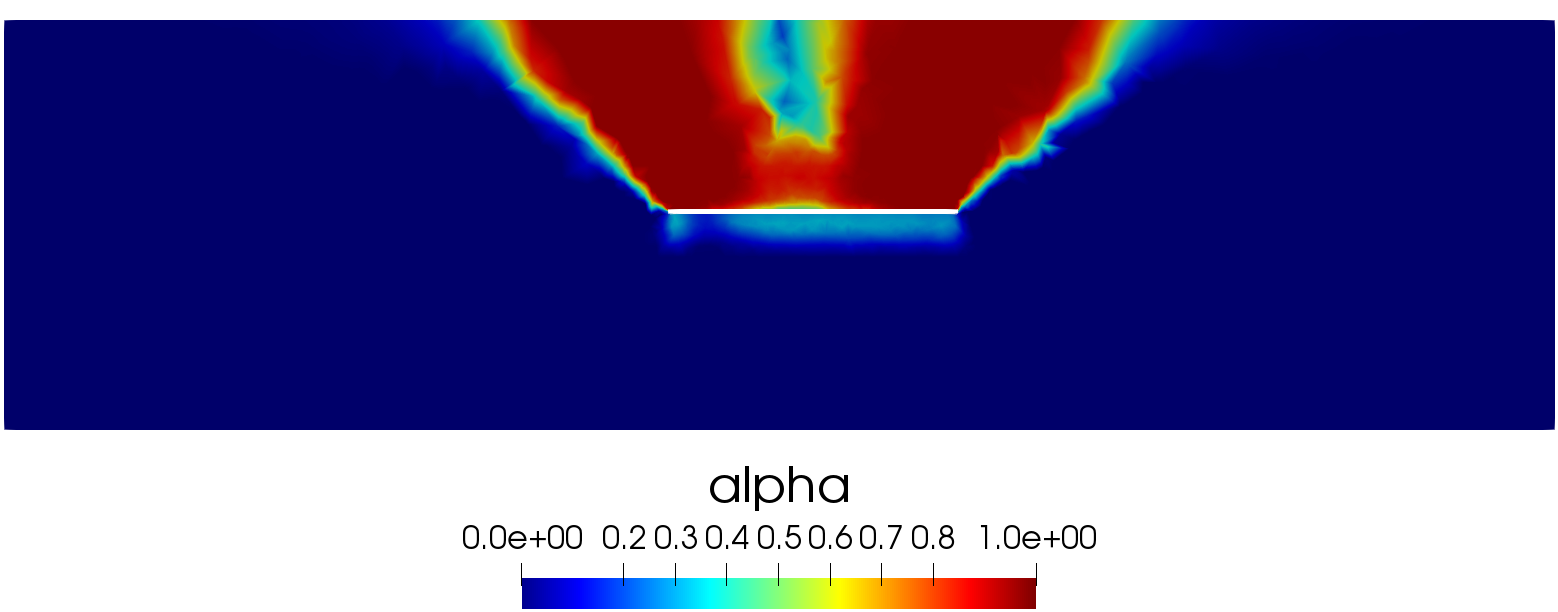}
		\caption{$\Omega(t_{40})$}
	\end{subfigure}
	\caption{Damage field distribution in the rock mass for Model 2 and $C_L=0.9$.}\label{mod2latcutCL09}
\end{figure}

\begin{figure}[h!]
	\centering
	\begin{subfigure}[b]{0.49\textwidth}
		\centering
		\includegraphics[width=\textwidth]{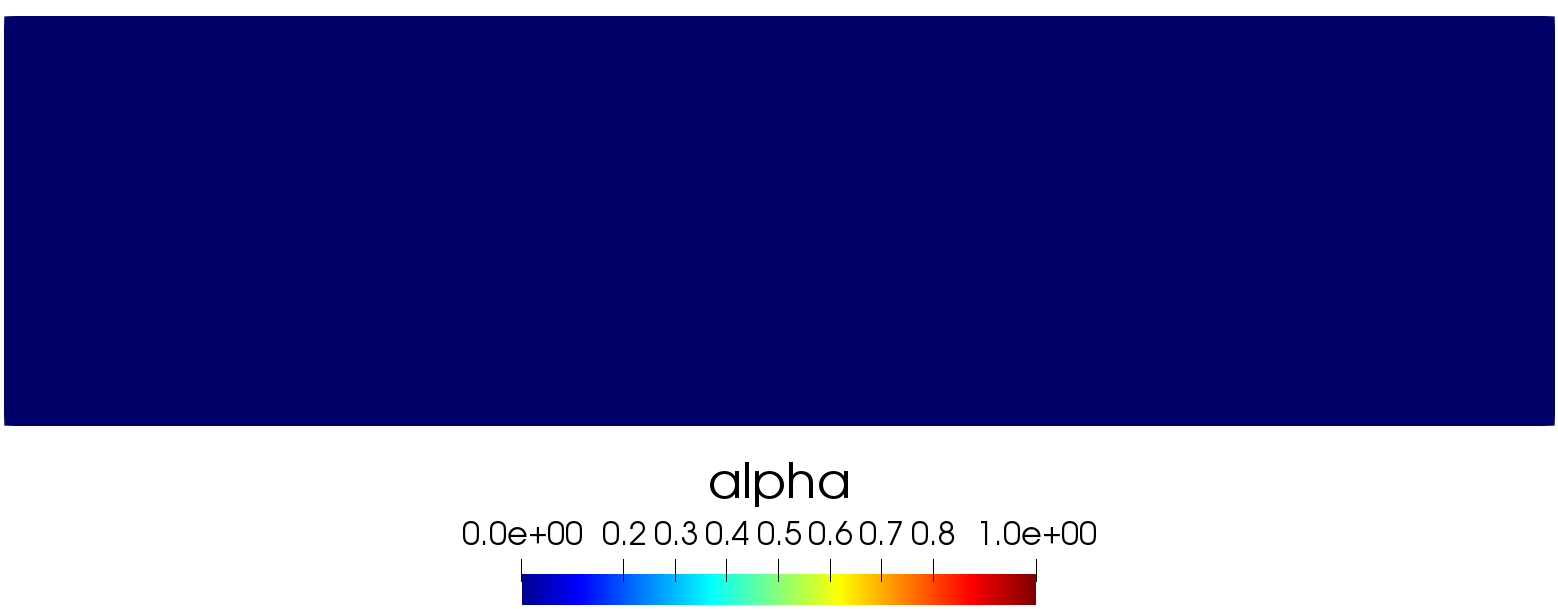}
		\caption{$\Omega(t_{0})$}
	\end{subfigure}
	\begin{subfigure}[b]{0.49\textwidth}
		\centering
		\includegraphics[width=\textwidth]{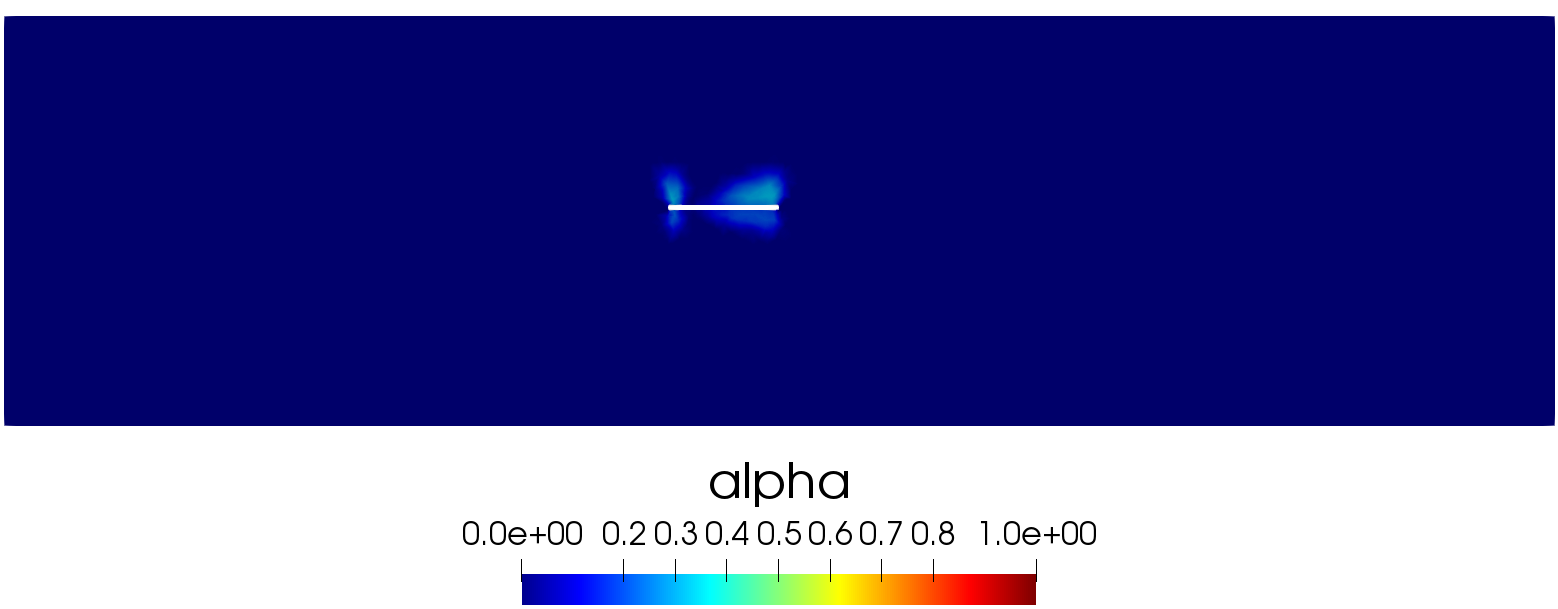}
		\caption{$\Omega(t_{15})$}
	\end{subfigure}
	\begin{subfigure}[b]{0.49\textwidth}
		\centering
		\includegraphics[width=\textwidth]{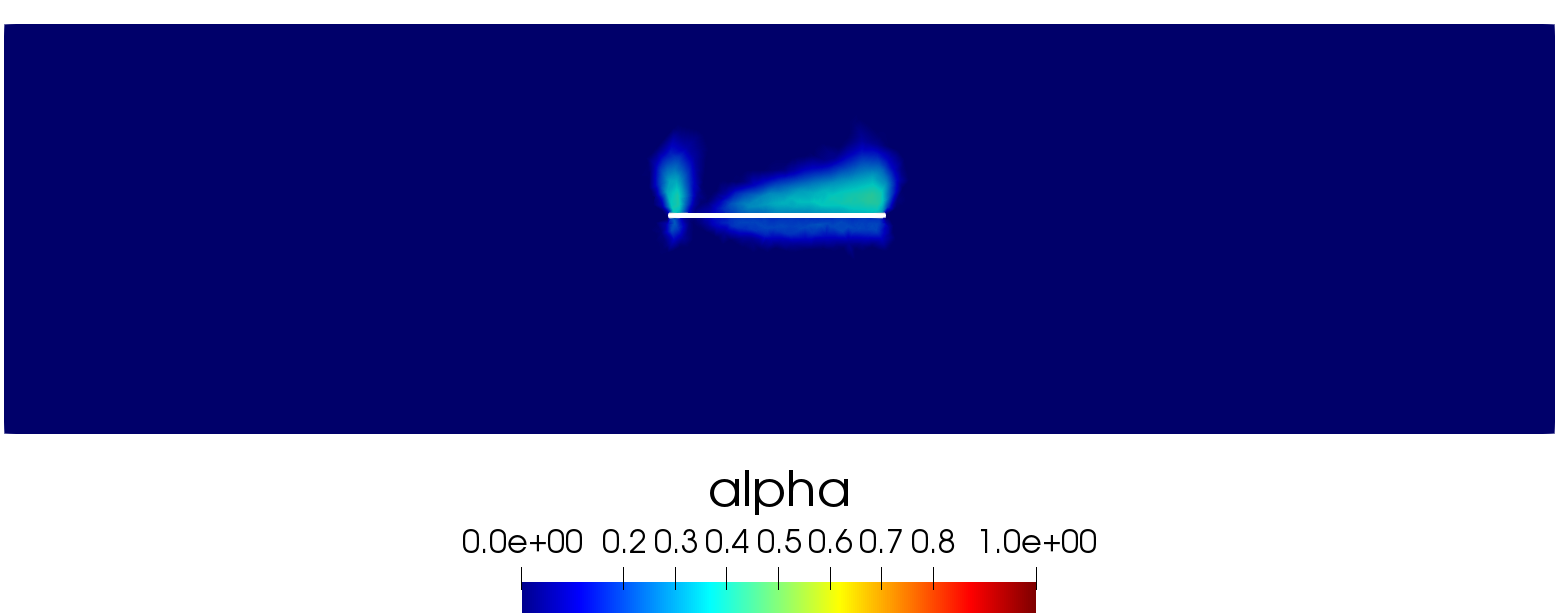}
		\caption{$\Omega(t_{30})$}
	\end{subfigure}	
	\begin{subfigure}[b]{0.49\textwidth}
		\centering
		\includegraphics[width=\textwidth]{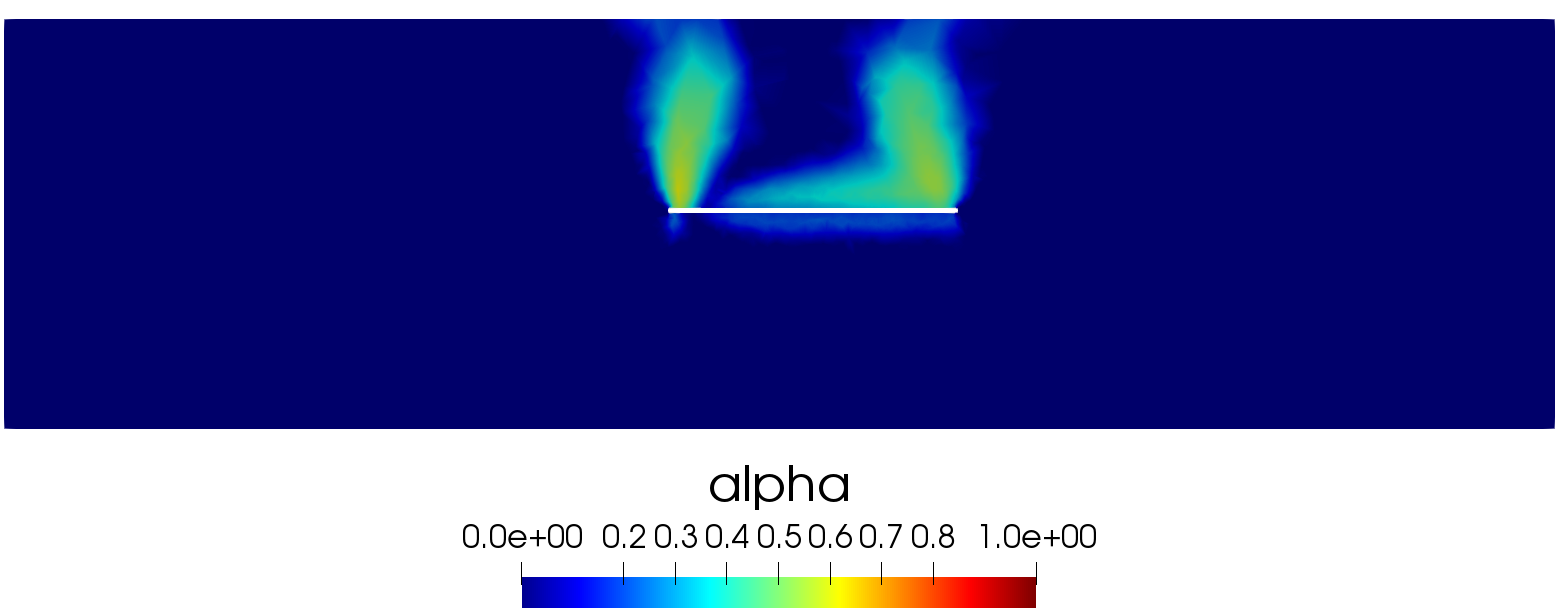}
		\caption{$\Omega(t_{40})$}
	\end{subfigure}
	\caption{Damage field distribution in the rock mass for Model 3 and $C_L=0.9$.}\label{mod3latcutCL09}
\end{figure}

\begin{figure}[h!]
	\centering
	\begin{subfigure}[b]{0.49\textwidth}
		\centering
		\includegraphics[width=\textwidth]{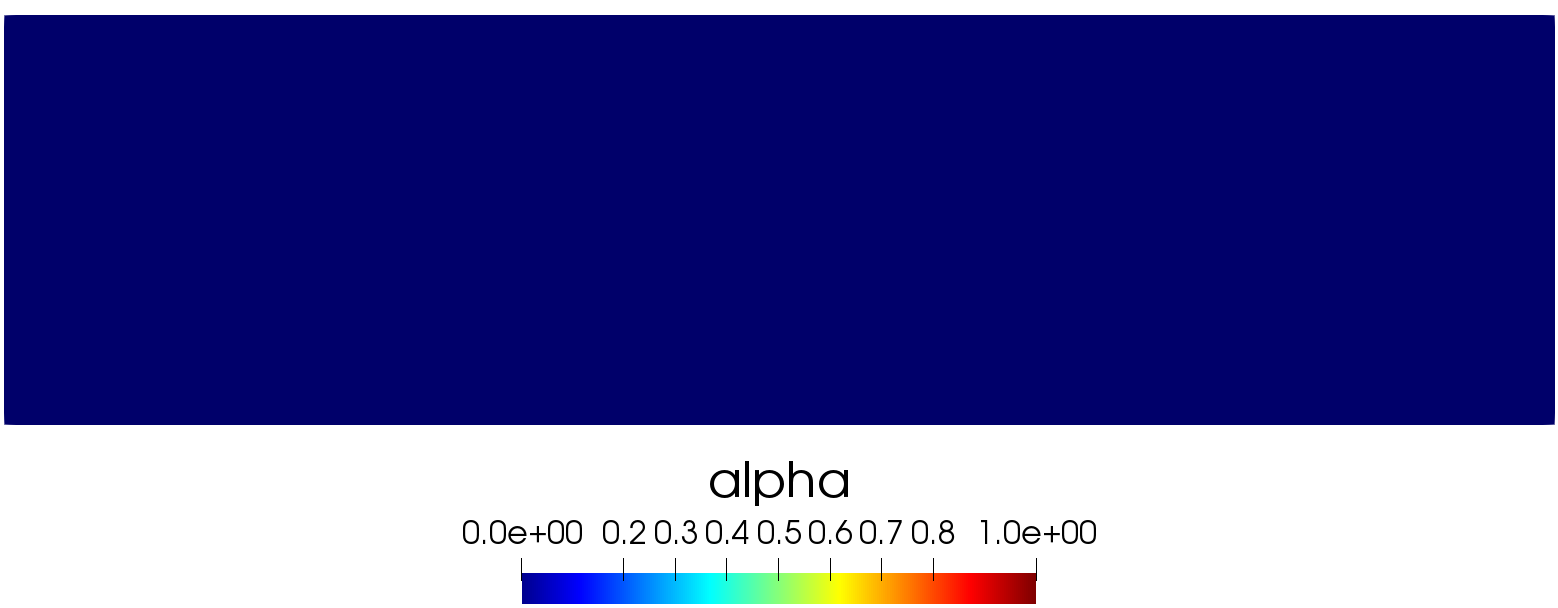}
		\caption{$\Omega(t_{0})$}
	\end{subfigure}
	\begin{subfigure}[b]{0.49\textwidth}
		\centering
		\includegraphics[width=\textwidth]{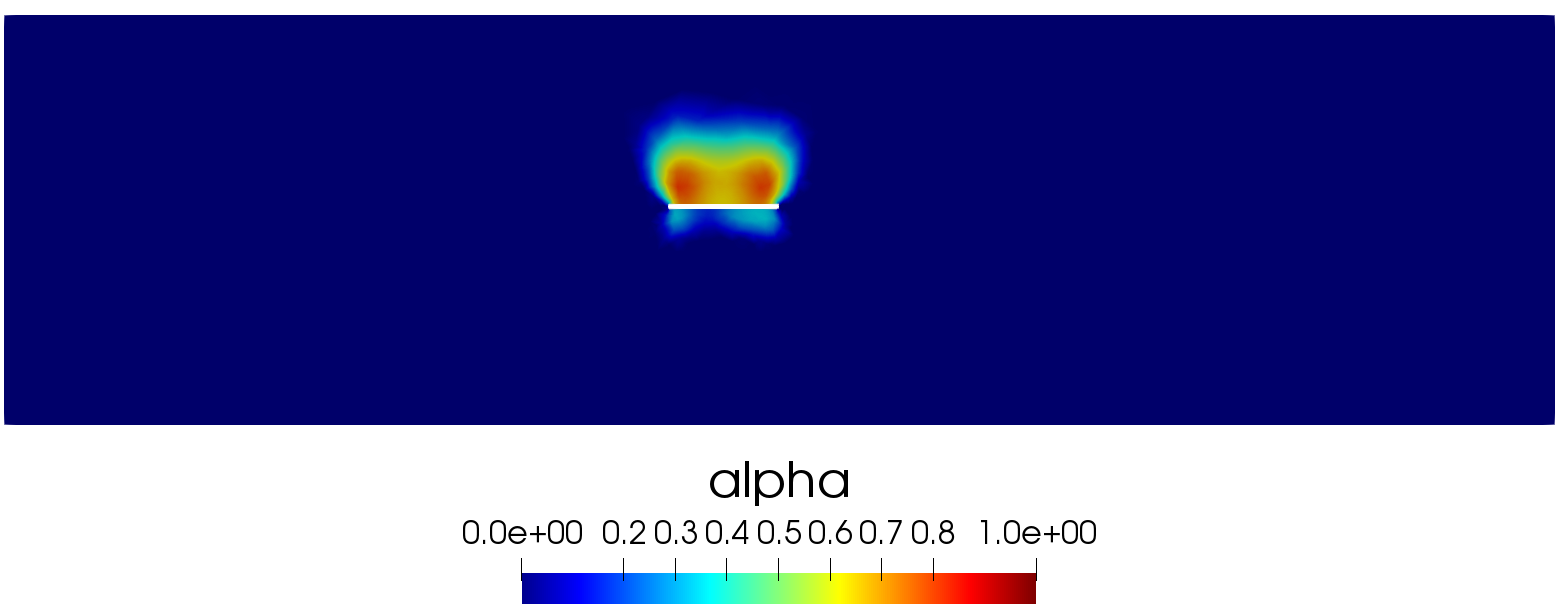}
		\caption{$\Omega(t_{15})$}
	\end{subfigure}
	\begin{subfigure}[b]{0.49\textwidth}
		\centering
		\includegraphics[width=\textwidth]{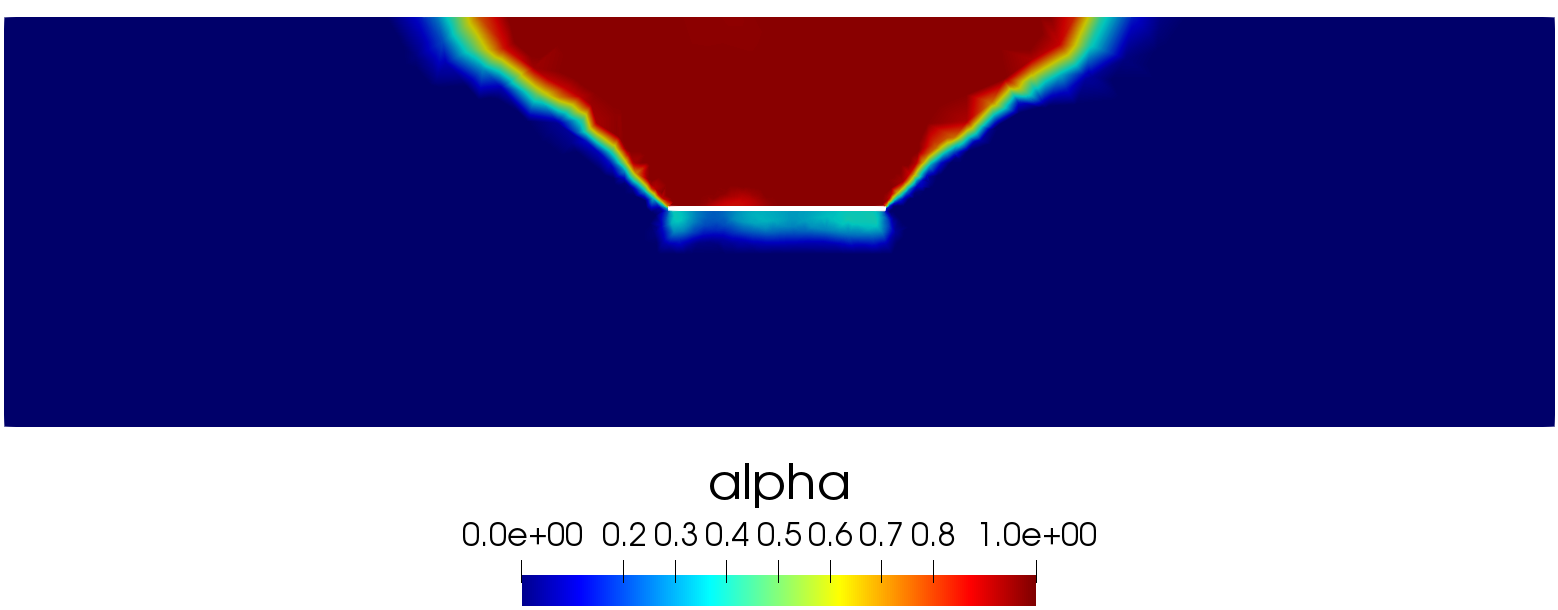}
		\caption{$\Omega(t_{30})$}
	\end{subfigure}	
	\begin{subfigure}[b]{0.49\textwidth}
		\centering
		\includegraphics[width=\textwidth]{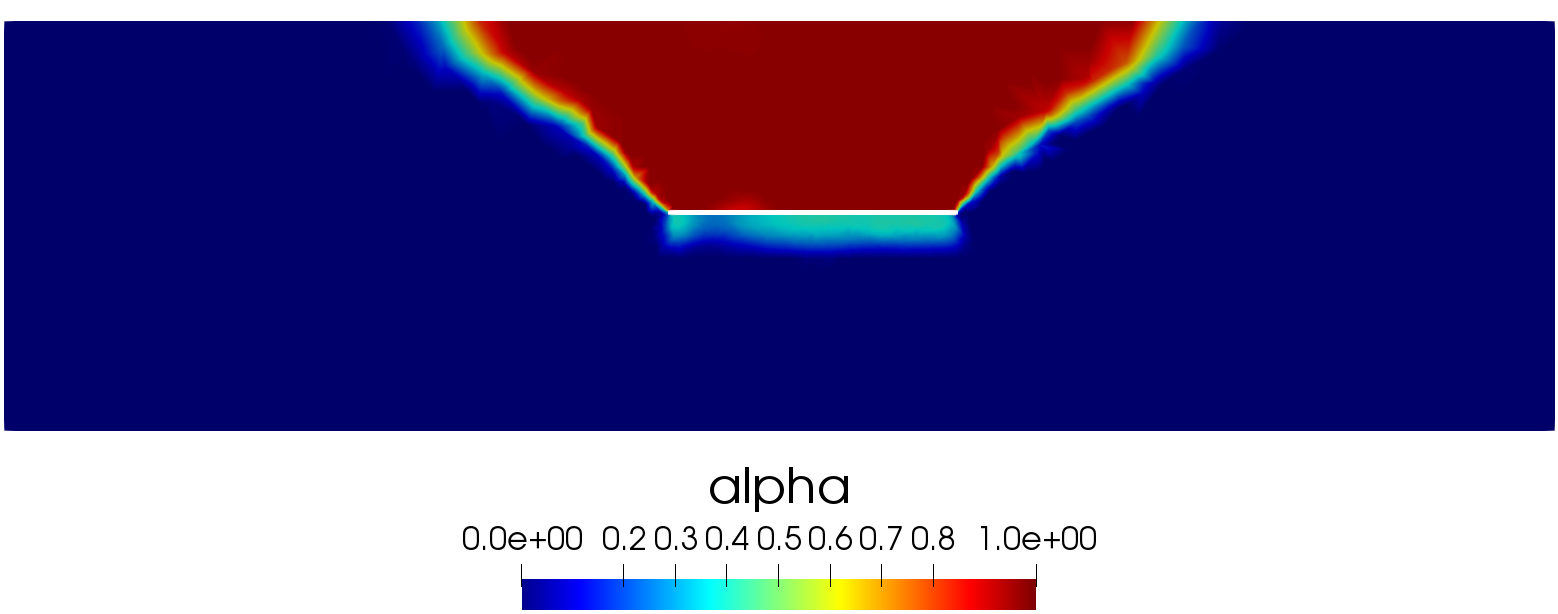}
		\caption{$\Omega(t_{40})$}
	\end{subfigure}
	\caption{Damage field distribution in the rock mass for Model 4 and $C_L=0.9$.}\label{mod4latcutCL09}
\end{figure}

\begin{figure}[h!]
	\centering
	\begin{subfigure}[b]{0.49\textwidth}
	    \centering
		\includegraphics[width=\textwidth]{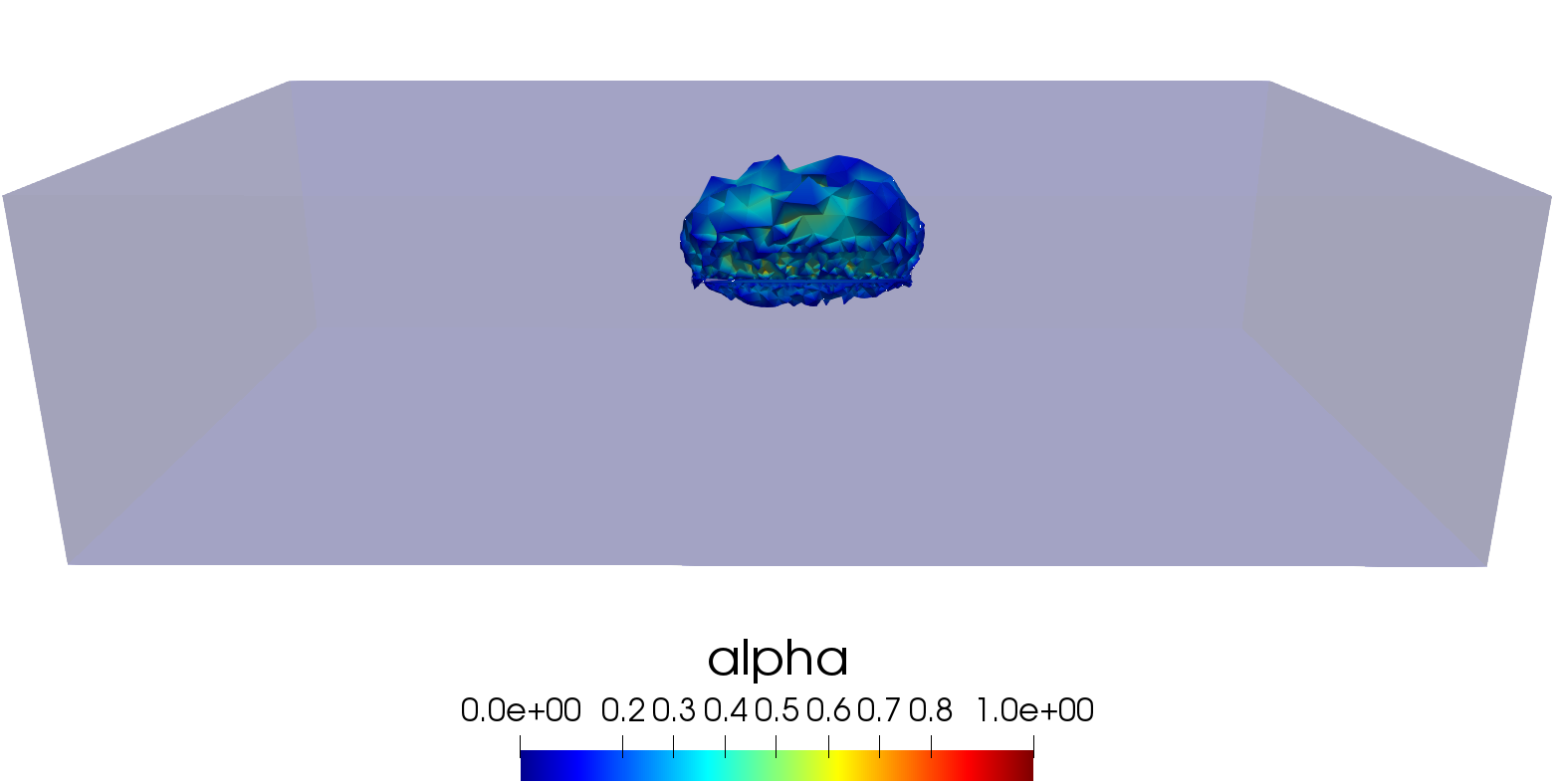}
		\caption{Model 1.}
	\end{subfigure}
	\begin{subfigure}[b]{0.49\textwidth}
		\centering
		\includegraphics[width=\textwidth]{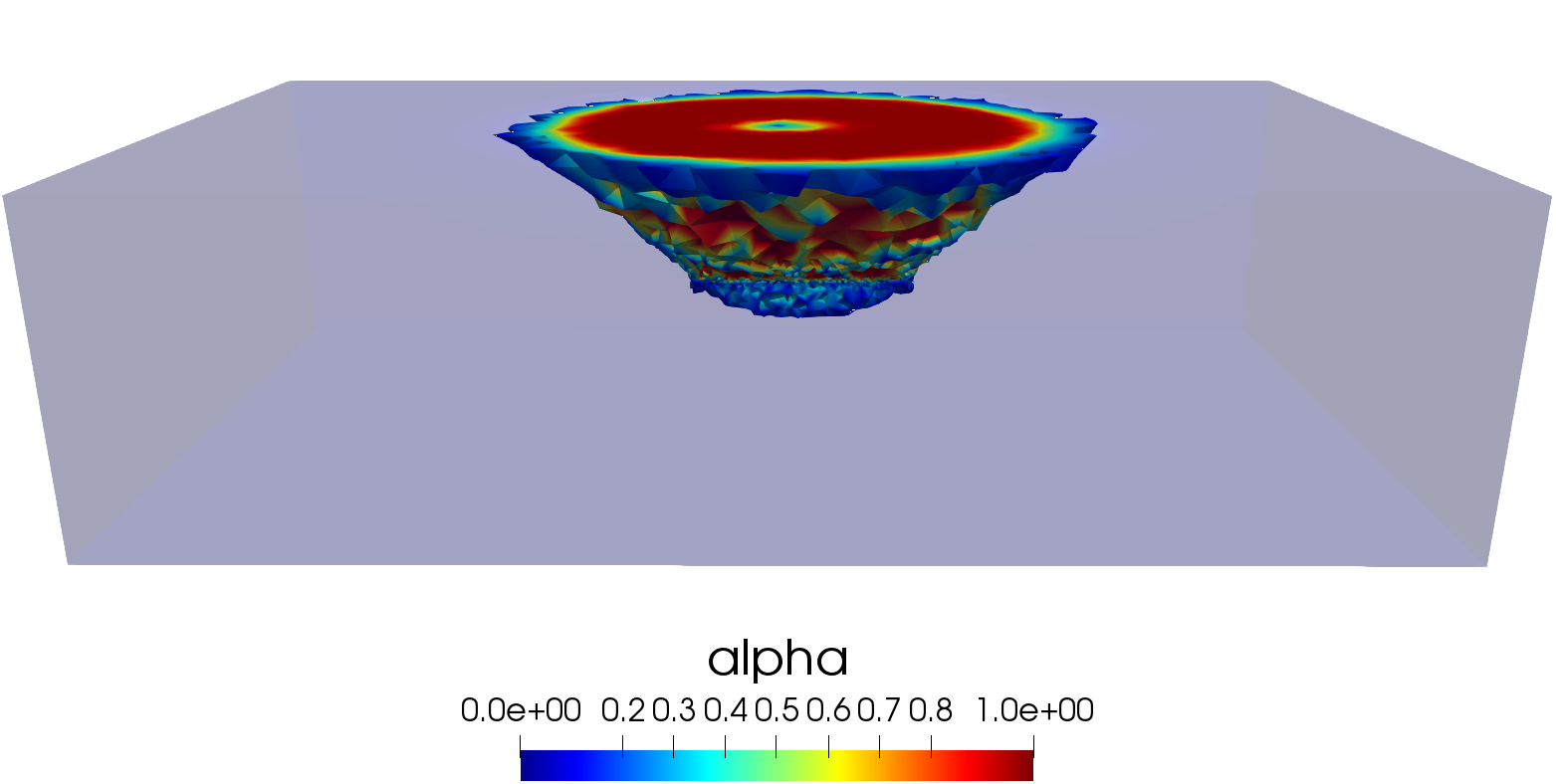}
		\caption{Model 2.}
	\end{subfigure}
		\centering
	\begin{subfigure}[b]{0.49\textwidth}
		\centering
		\includegraphics[width=\textwidth]{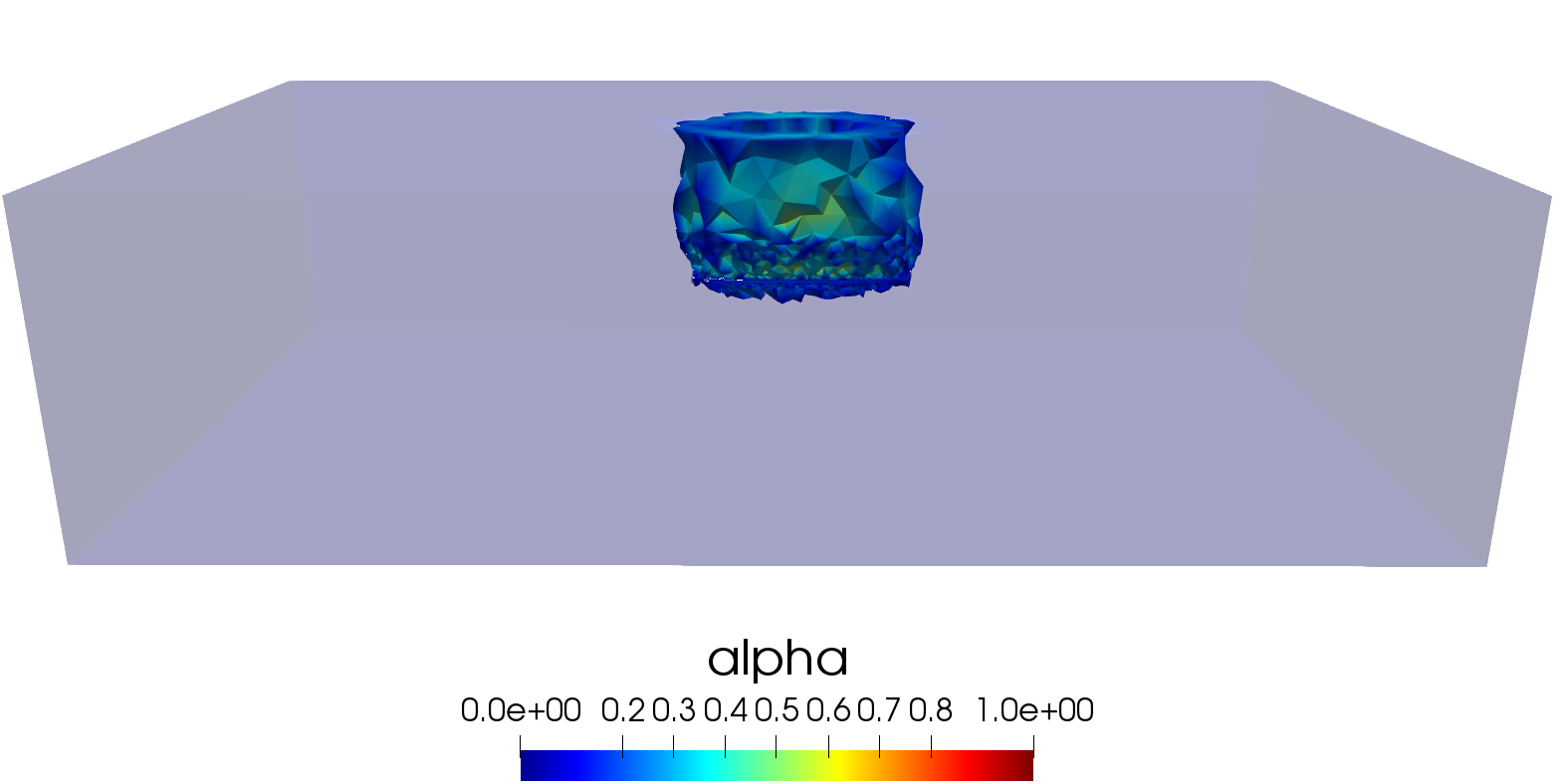}
		\caption{Model 3.}
	\end{subfigure}
	\begin{subfigure}[b]{0.49\textwidth}
		\centering
		\includegraphics[width=\textwidth]{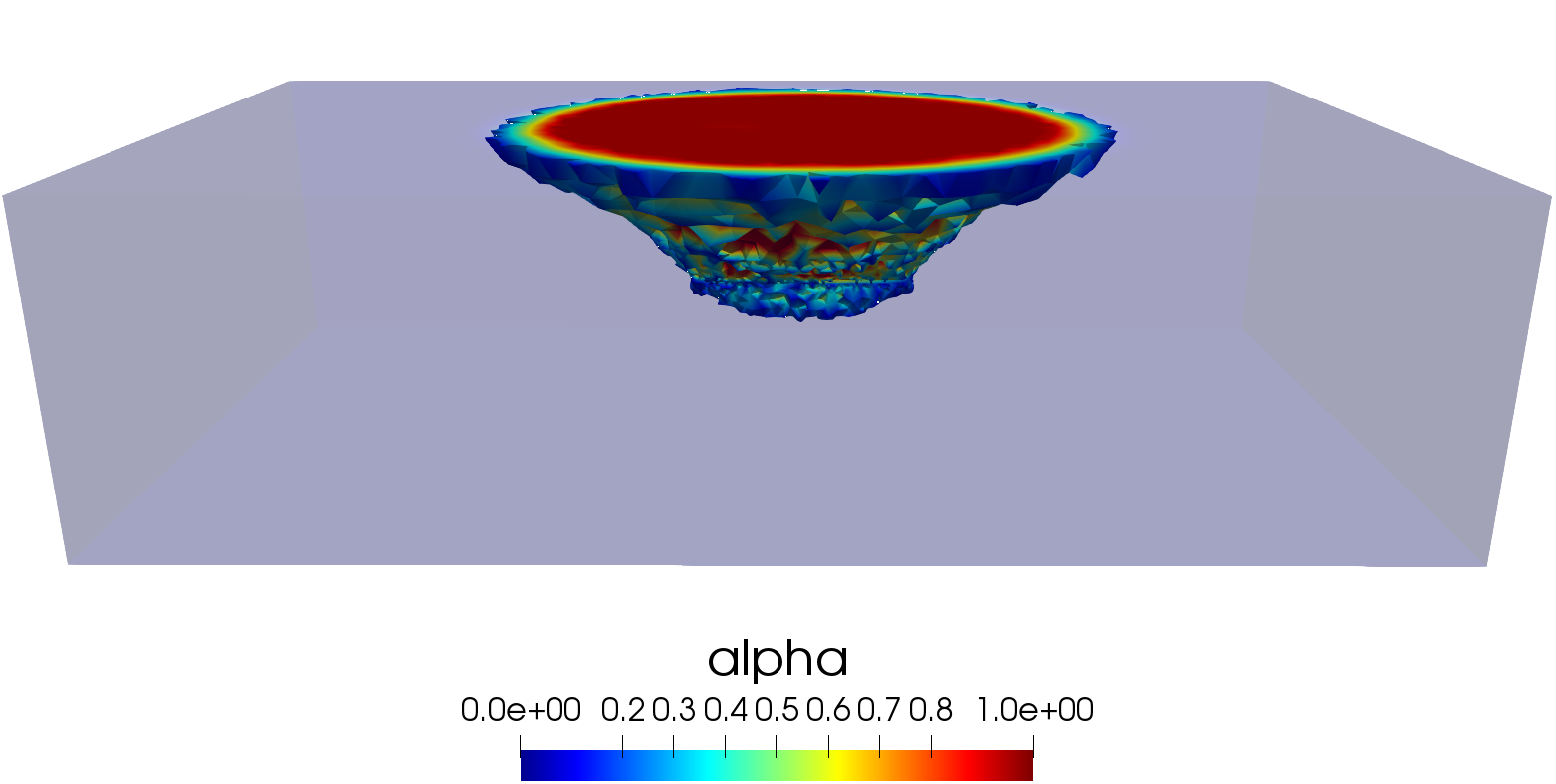}
		\caption{Model 4.}
	\end{subfigure}
	\caption{Damage field distribution in the rock mass in 3D for $C_L=0.9$.}\label{mods_w1-1000_cl-09_40_3D}
\end{figure}

Figures \ref{mod1latcutCL09}-\ref{mod4latcutCL09} display the distribution of the damage when the cavity advance in time using the Algorithm \ref{alg:damage2}. In the same way as using the Algorithm \ref{alg:damage}, the damage is close to zero almost everywhere but with a lower value and not reaching the upper boundary of the domain in the Model 1. In the Models 2 and 3 the material is fully damage over the cavity in the same way as the Algorithm \ref{alg:damage}. The Figure \ref{mods_w1-1000_cl-09_40_3D} displays the distribution of the damage in the final mesh step in 3D.

\section{Conclusions}\label{Con}

We propose a faster algorithm to modeling the block caving process in underground mining based in shear-compression damage model \cite{bonnetier2020shearcompression}. This algorithm is based in the usual alternate algorithm used in damage models \cite{marigo2016overview,bonnetier2020shearcompression}. We show that our new algorithm recover the damage produce by the block caving process using fewer iterations that the usual algorithm to solve the damage problem. In this work, we consider 3D examples to model different damage laws and compare the results produced by the our fast algorithm with the classical alternate algorithm, where where it can be seen that the same results are obtained using fewer iterations, that is, with a lower time.

\FloatBarrier

\bibliography{mybibfile}
\bibliographystyle{amsplain}

\end{document}